\documentclass[absolute]{mymathart}
\usepackage[theorems]{mymathmacros}
\usepackage[usenames,dvipsnames]{color}
\usepackage{subfig}
\usepackage[shortlabels]{enumitem}
\usepackage{hyperref}
\usepackage{amsmath}
\usepackage{amssymb}
\usepackage{graphicx}
\usepackage{epigraph}
\usepackage{amsfonts,amssymb,color}
\usepackage[mathscr]{eucal}
\usepackage{amsmath, amsthm}
\usepackage{mathrsfs}
\usepackage{amsbsy}
\usepackage{float}
\usepackage{verbatim}
\usepackage{amsfonts} 
\usepackage{pifont}
\usepackage{tikz-cd}
\usepackage{import}
\usepackage{xifthen}
\usepackage{pdfpages}
\usepackage{transparent}
\usepackage{dsfont}

\title[On the $n$-matings of polynomials]{On the $n$-matings of polynomials}

\author{Liangang Ma}

\address{Liangang Ma, School of Mathematics and Statistics, Ludong University, Yantai 264025, Shandong, P. R. China.} 
\email{maliangang000@163.com}

\thanks{The work is supported by NSFC-12001056.}  
\newtheorem{theorem}[subsubsection]{Theorem}
\newtheorem{theorem0}[subsection]{Theorem}
\newtheorem{lemma}[subsubsection]{Lemma}
\newtheorem{pro}[subsubsection]{Proposition}

\newtheorem*{Bottcher's Theorem}{B\"ottcher's Theorem}

\newtheorem*{BEZ}{Baer-Epstein-Zieschang Theorem}
\newtheorem{coro}[subsubsection]{Corollary}

\newtheorem{exa}[subsubsection]{Example}
\newtheorem{rem}[subsubsection]{Remark}
\newtheorem{defn}[subsubsection]{Definition}

\numberwithin{equation}{subsection}

\begin{document} 

\begin{abstract}
We introduce the notion of $n$-mating in this work, which includes the classical mating of polynomials as a special case. The new notion brings further links between the polynomial world and the rational world than the classical one, as well as a natural classification of rational maps according to their $n$-unmatability. We classify the hyperbolic $2$-matings according to the (non-)existence of orientation-reversing equators for them. For rational maps admitting orientation-reversing equators, we describe their dynamics via matings of half polynomials. There are diverse types of $n$-matings  from the bicritical family and the degree-$2$ capture family exhibited in our explorations, which demonstrates the ubiquity of them. Finally we consider the postcritical realization programme of rational maps (among the atomic and mating families respectively). The compositive trick is exploited to deal with problems in the programme.                 
\end{abstract}
 
 \maketitle

\tableofcontents

\section{Introduction}\label{sec1}

\numberwithin{equation}{section}

Mating of two quadratic polynomials was introduced by Douady and Hubbard (\cite{Dou}) in the last 80s. The rough idea is to paste together the connected and locally connected  filled Julia sets of two polynomials of the same degree along the landing points of their external rays inversely, which may induce a rational map on the Riemann sphere $\mathbb{P}^1(\mathbb{C})$ in suitable sense. From the reversing point of view one can consider splitting the Julia set of a rational map into Julia sets of two polynomials of the same degree (unmating).  One can refer to Douady-Hubbard (\cite{DH2}) and Milnor (\cite{Mil4}) for the external rays of polynomials and the Mandelbrot set.

There are various notions of matings  subtly distinct from one another in the literature, refer to \cite{MP} by Meyer and Petersen for a survey. One of the main concerns is that in which sense the induced map can be understood as a rational map. For example, in the sense of Thurston equivalence ($\stackrel{Thu}{\approx}$), topological conjugacy ($\stackrel{top}{\approx}$) and conformal conjugacy ($\stackrel{con}{\approx}$). 

Early investigations of the topic focus on the existence of matings of two polynomials of the same degree. Rees (\cite{Ree1, Ree2}), Shishikura (\cite{Shi}) and Tan (\cite{Tan1}) showed that the matings of two postcritically finite quadratic polynomials $\{z^2+c_j\}_{j=1,2}$ with parameters $c_1, c_2$ not in conjugate limbs of the Mandelbrot set exist in the sense of Thurston equivalence, moreover, the induced rational map is unique up to M\"obius conjugations. There is an example of mating  by Pilgrim such that uniqueness of the induced rational map fails in \cite{Mil2}. The importance of the postcritically finite polynomials is that they serve as representatives of hyperbolic polynomials in the same hyperbolic components, while the matability can be transferred to pairs of hyperbolic polynomials by quasi-conformal surgery. For instance, one can refer to Ha\"issinsky and Tan (\cite{HT}). 

Since the Julia sets of non-hyperbolic polynomials may lost local connectivity, the mateability of non-hyperbolic polynomials is subject to it. One is recommended to work of Avila-Lyubich (\cite{AL1, AL2}), Hubbard (\cite{Hub}), McMullen (\cite{Mcm}), Thurston (\cite{Thu})  and Yoccoz (\cite{Yoc}) for the topic on local connectivity of non-hyperbolic polynomials.  In \cite{Pet} Petersen proved that the Julia set of a Siegel quadratic polynomial $e^{2\pi \theta i}z+z^2$ with its rotation number $\theta \in (0,1)$ being a bounded-type irrational is always locally connected. He exploited the 'Julia' sets of some Blaschke functions to model the Julia sets of these polynomials to achieve his goal (see also \cite{Yam} and \cite{PZ}). Upon his result Yampolsky and Zakeri (\cite{YZ}) showed that any two Siegel polynoials $\{e^{2\pi \theta_j i}z+z^2\}_{j=1,2}$ with $\theta_1,\theta_2\in (0,1)$ being bounded-type irrationals and $\theta_1+\theta_2\neq 1$ are conformally matable, while the induced rational map is unique up to M\"obius conjugations.

Yoccoz (\cite{Yoc}) proved local connectivity of Julia sets of the finitely renormalizable polynomials. Aspenberg and Yampolsky considered the matability of non(finitely)-renormalizable quadratic polynomials with the basilica. They showed that a non-renormalizable quadratic polynomial $z^2+c$ without non-repelling periodic orbits is always conformally matable with the basilica polynomial, as long as $c$ is in the Mandelbrot set $\setminus$ $1/2$-limb. The mating is also unique up to M\"obius conjugations. They found bubble rays on the dynamical planes of concerned maps and parabubble rays on the parameter slice so that Yoccoz’ puzzle method applies in their case. 

Most of the above works focus on quadratic matings, while some results on matings of higher-degree polynomials are available. Shishikura and Tan (\cite{ST}) investigated the matability of some cubic polynomials, with some new phenomenons discovered as opposed to the quadratic case. Tan (\cite{Tan2}) explored some family of matings of cubic polynomials which result in Newton maps. A conjecture of her on the topic is obtained by Aspenberg and Roesch in \cite{AR}. 

Frequently there are rational slices involved in the matability of polynomials. While the target rational slices contain the matings we aim for, they also contain other interesting maps, for example, the captures in \cite{AR, AY, ST}. These mysterious maps typically can not be interpreted as matings in ones' expectations, while our notion of $n$-mating will break through these barriers.  

In the following we switch to the opposite viewpoint of mating: unmating rational maps into pairs of polynomials of the same degree. The unmatability of hyperbolic rational maps depends on an important notion given by    
Wittner in his thesis (\cite{Wit}): equator. He proved that if a postcritically finite rational map admits an equator, then the map is the formal mating (a proto-type mating) of a pair of polynomials in the sense of Thurston equivalence. Meyer strengthened his result by giving necessary and sufficient criteria  for the unmatability of postcritically finite hyperbolic rational maps. He showed that a postcritically finite hyperbolic rational map is (arises as) a mating if and only if  it admits an equator, in the sense of topological conjugacy.   

Milnor constructed an explicit example of mating of dendrites in \cite{Mil2}, whose Julia set spectacularly appears as some kind of space-filling (Peano, \emph{cf.} \cite{Mey3}) curve. Meyer also gave some sufficient (unfortunately not necessary)  criterion  for a postcritically finite rational map with extreme Julia set to arise as a mating:  if a postcritically finite rational map with its Julia set being $\mathbb{P}^1(\mathbb{C})$ admits a pseudo-equator, then it is a mating. It follows from his criterion that we better take iterates in considering the unmatability of non-hyperbolic rational maps, as any sufficiently high iterate of a postcritically finite rational map with its Julia set being $\mathbb{P}^1(\mathbb{C})$ admits some pseudo-equator (\cite{Mey2}).

On inspiration of all the above works, we raise the notion of $n$-mating of polynomials (Definition \ref{def2}). The classical mating corresponds to $1$-mating under the notion, while it essentially generalises the classical one. Roughly speaking, a rational map $R$ is an $n$-mating for some integer $n\in\mathbb{N}_+$ if $R^n$ is a classical mating while any of its lower iterates is not a classical mating.  Our first result deals with the unmatability of hyperbolic postcritically finite rational maps with orientation-reversing (OR) equators (Definition \ref{def4}).    

\begin{theorem0}\label{thm1} 
Let $R:\mathbb{P}^1(\mathbb{C})\rightarrow\mathbb{P}^1(\mathbb{C})$ be a hyperbolic postcritically finite rational map. If $R$ admits an orientation-reversing equator, then $R^2$ arises as a $1$-mating. 
\end{theorem0} 

The OR (pseudo-)equator is discovered by Meyer (\cite{Mey1}). Be alert that we do not claim $R$ arises as a $2$-mating in Theorem \ref{thm1}. There is some time that we think a rational map can not admit an equator and an OR equator simultaneously, surprisingly, this is wrong. We discover rational maps with both equators and OR equators: the hermaphroditic matings, in the bicritical family (\cite{Mil3}). Thus these  maps admit at least two dual types of unmating in fact. This means a hyperbolic rational map with an OR equator may arise as either a $1$-mating or a $2$-mating. An interesting query is that whether all $2$-matings arise due to OR equators. Surprisingly, this is also wrong. We discover pure type $2$-matings in the capture family (\cite{AY}) in $V_2$ (\cite{Ree3}). These discoveries provide a clear classification of $2$-matings: the OR matings and non-OR (primitive) $2$-matings (see Section \ref{sec3}).

We then give a description of dynamics of the hyperbolic postcritically finite rational maps with OR equators via matings of pairs of half polynomials ($\amalg$, Definition \ref{def5}).  
\begin{theorem0}\label{thm5}
For any hyperbolic postcritically finite rational map $R: \mathbb{P}^1(\mathbb{C})\rightarrow\mathbb{P}^1(\mathbb{C})$ admitting an OR equator, there exists a pair of half polynomials $\{\mathcal{P}_{wb}, \mathcal{P}_{bw}\}$, such that 
\begin{equation}\label{eq33}
R\stackrel{Thu}{\approx} \mathcal{P}_{wb}\amalg \mathcal{P}_{bw}
\end{equation}
on $\mathbb{P}^1(\mathbb{C})$.
\end{theorem0}

The half polynomials (Definition \ref{def1}) differ from the typical polynomials in the sense that we distinguish their domain and range as different Riemann spheres. 

We point out that the $n$-matings with folds $n\geq 2$ exist extensively in various families of rational maps well studied in the literature. We explore them in the bicritical family and the $V_2$ family in the following.

The investigation of bicritical rational maps was initiated by Milnor (\cite{Mil3}). In fact his focus is laid on the moduli space (up to M\"obius conjugacies). One can find some interesting topological descriptions on the  moduli spaces (up to M\"obius conjugacies) of degree-$d$ rational maps in DeMarco's work \cite{Dem} and of degree-$d$ polynomials in DeMarco-McMullen's work \cite{DM}. For most recent developments on bicritical rational maps one can refer to Koch-Lindsey-Sharland (\cite{KLS}). There are diverse interests on the family with attractive dynamics in it, as to our concern, we discover $4$-matings, hermaphroditic matings and OR matings in the bicritical family.  

\begin{theorem0}\label{thm8}
\begin{enumerate}[(A).]
\item For $k\geq 2$, any bicritical map $R_{a,b,c,d}$ with $(a,b,c,d) \in \Omega_+^{k}$ is a $4$-mating. 

\item For $k\geq 3$, any bicritical map $R_{a,b,c,d}$ with $(a,b,c,d) \in \Omega_-^{k}$ is a hermaphroditic mating. 

\item Any bicritical map $R_{a,b,c,d}$ with $(a,b,c,d) \in \Omega_-^{2}$ is an orientation-reversing mating. 
\end{enumerate}
\end{theorem0}

\begin{figure}[h]
\centering
\includegraphics[scale=1]{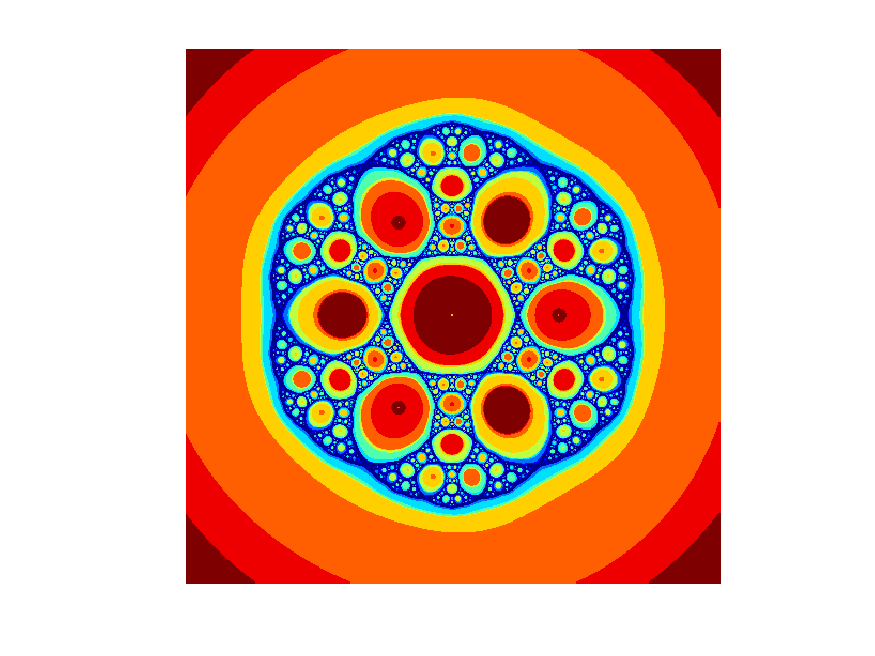}
\caption{Dynamics of $R_{\omega_{+,3}}$}
\label{fig21}
\end{figure}

Theorem \ref{thm8} results in new interesting interpretation of dynamics of certain maps in the bicritical family. For example, the map 
$$R_{\omega_{+,3}}=\cfrac{z^3-1}{z^3+1}$$
with parameter in $\Omega_+^{3}$ arises as a $4$-mating. Its Fatou and Julia sets ($\mathcal{F}(R_{\omega_{+,3}})$ and $\mathcal{J}(R_{\omega_{+,3}})$) are depicted in Figure \ref{fig21}.

The map does not arise as a classical mating, neither does $R_{\omega_{+,3}}^2$ nor $R_{\omega_{+,3}}^3$. However, since $R_{\omega_{+,3}}^4$ arises as the classical mating and $\mathcal{J}(R_{\omega_{+,3}})=\mathcal{J}(R_{\omega_{+,3}}^4)$ ($\mathcal{F}(R_{\omega_{+,3}})=\mathcal{F}(R_{\omega_{+,3}}^4)$), we can still interpret its Julia set as the Julia set of mating of two degree-$81$ polynomials.  

We also explore the captures (refer to \cite{Luo} by Luo) belonging to the degree-$2$ rational family well-studied by Aspenberg and Yampolsky (\cite{AY}). These captures are known not to arise as classical matings  (Theorem \ref{lem9}). However, when the $n$-unmatability is considered, things change dramatically for the odd-generation captures. 

\begin{theorem0}\label{thm7}
The postcritically finite captures of even generations in the family $\Big\{R_a=\cfrac{a}{z^2+2z}\Big\}_{a\in \mathbb{C}\setminus\{0\}}$ are all atomic rational maps, while the postcritically finite captures of odd generations  are all primitive $2$-matings.
\end{theorem0}

The appearance of primitive $2$-matings in the family astounds us greatly. Existence of $n$-matings with $2$ and higher folds is obviously not confined in the above families. For example, we can imagine that $n$-matings will appear in the families studied by Shishikura-Tan in \cite{ST} and by Aspenberg-Roesch in \cite{AR}.  

Rational maps are traditionally classified into the hyperbolic family and non-hyperbolic family, according to whether all Fatou components of the maps are attractive. The non-hyperbolic rational maps can be once more classified into some sub-families, according to the distinguished features of some Fatou components of the rational maps, for example, the parabolic family, the Siegel family, the Herman family, $\cdots$. One can refer to Sullivan's Nonwandering Theorem (\cite{Sul}). These classifications set basic regulations on the dynamical behaviours of the maps in these kingdoms.  McMullen and Sullivan (\cite{MS}) developed the theory of holomorphic dynamical system  to investigate dynamics of rational maps, Kleinian groups, polynomial-like maps and holomorphic correspondences simultaneously. The unmating of holomorphic correspondences into polynomials and modular groups was carried out by Bullett-Penrose (\cite{BP}) and Bullett-Lomonaco (\cite{BL}). It was shown by Lee-Lyubich-Makarov-Mukherjee (\cite{LLMM1, LLMM2}) that any geometrically finite Schwarz reflection with respect to the cardioid or a circle is the mating of some anti-holomorphic polynomial with some reflection map.

Now we can give a classification of the rational maps according to a new criterion: the unmatability of the maps. Basically, all rational maps are classified into two disjoint families, the mating family and the atomic family. A map in the mating family is an $n$-mating for some integer $n\in\mathbb{N}_+$, while a map in the atomic family does not arise as any $n$-mating for any $n\in\mathbb{N}_+$. Maps in the former family are automatically classified into sub-families according to their folds such that they arise as $n$-matings. We are not sure how much this classification will be linked to further dynamics of maps in corresponding families, however, certainly it will. For example, there are distinctions when we consider the postcritical realization programme of (postcritically finite) rational maps, where the mating technique can be applied (refer to \cite[Theorem 5.6]{DKM}). Let $\bar{\mathbb{Q}}$ be the algebraic closure of $\mathbb{Q}$, we consider its projection on the Riemann sphere $\mathbb{P}^1(\bar{\mathbb{Q}})$. 

\begin{theorem0}\label{thm9}
Let $X\subset\mathbb{P}^1(\bar{\mathbb{Q}})$ be a finite set satisfying $\#X\geq 2$, then 
\begin{enumerate}[(A).]
\item There exists a rigid postcritically finite rational map arising as a $1$-mating such that its postcritical set coincides with $X$.

\item In case $\#X=2$, there does not exist any atomic postcritically finite rational map whose postcritical set coincides with $X$. In case $\#X=3$, there exists an atomic postcritically finite rational map such that its postcritical set coincides with $X$.
\end{enumerate}
\end{theorem0}

Obviously the mating family reaches further than the atomic family on realizing prescribed postcritical sets in Theorem \ref{thm9}. For results on realizing maps on prescribed postcritical sets of three points by atoms and matings, see Theorem \ref{thm10}.   

The postcritical realization programme is initiated by DeMarco, Koch and McMullen (\cite{DKM}). There are basically two aims of different level in the programme. In its general form, the first one is to decide which  countable (finite) sets in $\mathbb{P}^1(\bar{\mathbb{Q}})$ (or even $\mathbb{P}^1(\mathbb{C})$) can be realized (Definition \ref{def6}) by (postcritically finite) rational maps. DeMarco, Koch and McMullen had shown that any finite set in $\mathbb{P}^1(\bar{\mathbb{Q}})$ can be realized, with the employment of the Belyi polynomials (\cite{Bel}) in their construction. The second level on realizing maps (Definition \ref{def7}) on countable (finite) sets in $\mathbb{P}^1(\bar{\mathbb{Q}})$ (or even $\mathbb{P}^1(\mathbb{C})$) by (postcritically finite) rational maps is more difficult, see \cite[Question 1.2]{DKM}. If the prescribed set contains $2$ points, not every map can be realized, while every map can be realized if the set contains $3$ points according to \cite{DKM}. There are typically $k^k$ distinct maps on a set $X$ of $k$ points, which makes the realization of them a tremendous work. However, according to knowledge on the algebraic structure (refer to Theorem \ref{thm12}, \ref{thm13}, \ref{thm15}) on the collections of the (periodic and strictly pre-periodic) maps,  we can greatly simplify their question. The following result illustrates some simplification of \cite[Question 1.2]{DKM} in case $\#X=4$.

\begin{theorem0}\label{thm16}
Let $X=\{p_1, p_2, p_3, p_4\}\subset \mathbb{P}^1(\bar{\mathbb{Q}})$ with $4$ distinct points. Assuming we can realize the three maps 
\begin{center}
\begin{tikzpicture}
    \node (a) at (0,0) {$p_1$};
    \node (b) at (1,0) {$p_2$};
    \node (c) at (2,0) {$p_3$};
    \node (d) at (3,0) {$p_4,$};

    \draw[->] (a) to node {} (b);
    \draw[->] (b) to node {} (c);
    \draw[->] (c) to node {} (d);
    \draw[->] (d) to [in=-15, out=195, looseness=1] node {} (a);
\end{tikzpicture}
\begin{tikzpicture}
    \node (a) at (4,0) {$p_1$};
    \node (b) at (5,0) {$p_2$};
    \node (c) at (6,0) {$p_3$};
    \node (d) at (7,0) {$p_4$};

    \draw[->] (a) to node {} (b);
    \draw[->] (b) to node {} (c);
    \draw[->] (c) to [in=-20, out=200, looseness=1] node {} (a);
    \draw[->] (d) edge[loop right] node {,} (d);
\end{tikzpicture}
\begin{tikzpicture}
    \node (a) at (8,0) {$p_1$};
    \node (b) at (9,0) {$p_2$};
    \node (c) at (10,0) {$p_3$};
    \node (d) at (11,0) {$p_4$};

    \draw[->] (a) to node {} (b);
    \draw[->] (b) to node {} (c);
    \draw[->] (c) to node {} (d);
    \draw[->] (d) to [in=-20, out=200, looseness=1] node {} (b);
\end{tikzpicture}
\end{center}
on $X=\{p_1, p_2, p_3, p_4\}$  by some rational maps $R_{\mathfrak{P}|1}, R_{\mathfrak{P}||A1}, R_{\mathfrak{S}|A1}$ respectively with their critical-value sets $V(R_{\mathfrak{P}|1})=V(R_{\mathfrak{P}||A1})=V(R_{\mathfrak{S}|A1})=X$, then all the maps on $X$ are realizable by rational maps on $\mathbb{P}^1(\mathbb{C})$.
\end{theorem0}

The requirement $V(R_{\mathfrak{P}|1})=V(R_{\mathfrak{P}||A1})=V(R_{\mathfrak{S}|A1})=X$ in Theorem \ref{thm16} can be discounted, see for example Corollary \ref{cor8}.  Instead of achieving the postcritical realization programme via the Teichm\"uller theory, we use the compositive trick, which relies on established knowledge of postcritically finite rational maps realizing certain maps. The technique can be applied to realization of maps on sets of more than $4$ points. One can refer to Theorem  \ref{thm12} and Theorem \ref{thm14} for the realization of periodic maps and strictly pre-periodic maps  respectively.   

All the polynomials which give rise to the matings in our work can be identified via the Hubbard trees (see \cite{BFH} by Bielefeld-Fisher-Hubbard) available from maps on the postcritical sets partitioned by the (pseudo-)equators of the concerned rational maps (in case the Julia sets of the matings being $\mathbb{P}^1(\mathbb{C})$, see \cite[Section 10]{Mey1}). Attentive readers may carry out the programme. 

Some results may have extensions onto (expanding) Thurston maps (\emph{cf}. \cite{BM} by Bonk-Meyer) or (coarse expanding) conformal maps (\emph{cf}. \cite{HP} by Ha\"issinsky-Pilgrim), however, we take main interest in rational maps in this work. Moreover, we wonder whether some counterpart notions or results apply to, for instance, unmating of holomorphic correspondences.

The structure of the work is as following. In Section \ref{sec2} we introduce the main notations and concepts employed throughout the work, including the $n$-matings, (OR) equators, etc. Some results on the behaviour of (OR) equators are presented here. Section \ref{sec3} is devoted to the classification of the hyperbolic $2$-matings, according to whether they admit OR equators or not. Theorem \ref{thm1} is proved here. In Section \ref{sec4} we introduce the notion of half polynomials and employ matings of them to interpret the dynamics of rational maps with OR equators (Theorem \ref{thm5}). In the following two sections we show that $n$-matings are widespread in the rational world. In Section \ref{sec6}  we explore the $n$-unmatability of postcritically finite rational maps in the bicritical family. We will demonstrate that there are $4$-matings, hermaphroditic matings and OR matings in the family (Theorem \ref{thm8}). In Section \ref{sec5}  we explore the $n$-unmatability of postcritically finite capture maps in the $V_2$ family. We will demonstrate that there are atomic maps and primitive $2$-matings appearing alternatively in the family (Theorem \ref{thm7}).  The readers will experience a gallery of (OR) equators of sample maps and their pre-images in the bicritical family and the capture family throughout the two sections. In Section \ref{sec7} we carry out the postcritical realization programme of atomic maps and matings respectively (Theorem \ref{thm9}, \ref{thm10}). In the final section we introduce the compositive trick to deal with the postcritical realization programme. We will study the algebraic structure built on different collections of maps on prescribed sets,  together with the compositive trick this results in simplification of the postcritical realization programme (Theorem \ref{thm16}).

\numberwithin{equation}{subsection}

\section{Basic notations and some preliminary results}\label{sec2}

\subsection{Notations and concepts}
Let $\mathbb{C}$ be the complex plane and $\mathbb{P}^1(\mathbb{C})$ be the Riemann sphere. Let $\bar{\mathbb{Q}}$ be the collection of all algebraic numbers and $\mathbb{P}^1(\bar{\mathbb{Q}})$ be its projection onto $\mathbb{P}^1(\mathbb{C})$. Let 
$$\hat{\mathbb{C}}=\mathbb{C}\cup\infty S^1=\mathbb{C}\cup\{\infty e^{2\pi it}: 0\leq t<1\}$$ 
be the \emph{exploded complex plane}.

We typically use $\partial A, A^o, \bar{A}$ to mean the \emph{boundary}, \emph{interior} and \emph{closure} of some set $A$ respectively. The \emph{disjoint union} of two sets $A$ and $B$ is $A\sqcup B$. For a set $A$, a \emph{partition} of $A$ means a collection of non-empty mutually disjoint subsets $\{A_j\}_{j\in \Lambda}$ with $\Lambda$ being some index set, such that 
\begin{center}
$A=\cup_{j\in\Lambda} A_j$.
\end{center}

Let $\mathds{1}_A: A\rightarrow A$ be the \emph{identity map} on $A$ (we sometimes simply write $\mathds{1}$ if the set on which it is defined is clear). For a map $f: A\rightarrow A$ and some subset $B\subset A$, let $f|_{B}: B\rightarrow B$ be its \emph{restriction} on $B$.  

We typically use $R: \mathbb{P}^1(\mathbb{C})\rightarrow \mathbb{P}^1(\mathbb{C})$ to mean a rational map of degree $d$ on the Riemann sphere. Sometimes it can also be understood as $R: \hat{\mathbb{C}}\rightarrow \hat{\mathbb{C}}$. Especially, in case $R$ is a degree $d$ polynomial, we mean that 
\begin{equation}\label{eq15}
R(\infty e^{2\pi it})=\infty e^{2\pi idt}
\end{equation}
for all $t\in I=[0,1]$ on $\infty S^1$.  The collections of critical points and critical values of a rational map $R$ are denoted by $C(R)$ and $V(R)$ respectively. The \emph{postcritical set} $P(R)$ is
$$P(R)=\{R^k(c): c\in C(R), k\geq 1\},$$
which is sub-preserved (Definition \ref{def10}) by the map $R$. The three sets are crucial in exploring the dynamics of the map $R$. If  the cardinality $\#P(R)<\infty$ we call $R$ \emph{postcritically finite}. These notations apply to branched coverings on the (topological) $2$-sphere $S^2$, with the branching points serving as critical points.  The \emph{Julia} and \emph{Fatou} sets of a rational map $R$ are $\mathcal{J}(R)$ and $\mathcal{F}(R)$ respectively. A polynomial is typically denoted by $P: \mathbb{P}^1(\mathbb{C})\rightarrow \mathbb{P}^1(\mathbb{C})$, with its filled Julia set being $\mathcal{K}(P)$.

We say a set $B\subset A$ is \emph{completely invariant} under a map $f: A\rightarrow A$ if $$f(B)=B=f^{-1}(B).$$ 
Now we introduce some relationships between maps. For two branched coverings $f,g: S^2\rightarrow S^2$, we write $f\stackrel{top}{\approx} g$ to mean $f$ is \emph{topologically conjugate} to $g$. This means there exists some homeomorphism $\kappa: S^2\rightarrow S^2$, such that 
$$\kappa\circ f=g\circ\kappa$$
on $S^2$. In case there are Riemannian metrics on $S^2$, the relationship escalates into a \emph{conformal conjugacy} denoted by $f\stackrel{con}{\approx} g$ if the homeomorphism $\kappa$ is conformal on $S^2$ (or an open dense subset of $S^2$). 

For two postcritically finite branched coverings $f,g: S^2\rightarrow S^2$, we write $f\stackrel{Thu}{\approx} g$ to mean $f$ is \emph{Thurston equivalent} to $g$. This means there exist two orientation-preserving homeomorphisms $\kappa_1, \kappa_2: S^2\rightarrow S^2$ isotopic to each other \emph{rel.} $P(f)$, such that
$$\kappa_1\circ f=g\circ\kappa_2$$
on $S^2$.

Obviously these relationships are subject to the following order structure,    
\begin{center}
$\stackrel{con}{\approx} \Rightarrow \stackrel{top}{\approx} \Rightarrow \stackrel{Thu}{\approx}$ 
\end{center}
for postcritically finite maps.

We will use the B\"ottcher's Theorem (refer to \cite[Theorem 9.1]{Mil1}) in the following. It describes the behaviour of some holomorphic map around superattracting fixed points (in case of their existence) by $z^d$ with local degree $d\geq 2$. Let $\mathbb{D}\subset\mathbb{C}$ be the closed unit disc.  The polynomial version is enough for our aim in this work. 
\begin{Bottcher's Theorem}
Let $P: \mathbb{P}^1(\mathbb{C})\rightarrow \mathbb{P}^1(\mathbb{C})$ be a polynomial of degree $d\geq 2$ with connected Julia set. There exists a conformal map $\phi: \mathbb{P}^1(\mathbb{C})\setminus\mathcal{K}(P)\rightarrow\mathbb{P}^1(\mathbb{C})\setminus\mathbb{D}$, such that
\begin{center}
$\phi\circ P\circ\phi^{-1}(z)=z^d$
\end{center}
for $z\in\mathbb{P}^1(\mathbb{C})\setminus\mathbb{D}$. The map $\phi$ is unique up to multiplication by a $(d-1)$-th root of $1$.
\end{Bottcher's Theorem}

We call the map $\phi$ the  \emph{B\"ottcher's coordinate} of $P$ in the following. The coordinate extends continuously to the boundary $\mathcal{J}(P)$ if $\mathcal{J}(P)$ ($\mathcal{K}(P)$) is locally connected, according to Carath\'eodory (\cite[Theorem 17.14]{Mil1}).

\subsection{The $n$-mating}

To introduce the $n$-mating of polynomials, we need the classical mating. We start from some proto-type classical mating-the \emph{formal mating} of polynomials. For $\hat{\mathbb{C}}_w=\hat{\mathbb{C}}_b=\hat{\mathbb{C}}$, let 
\begin{center}
$P_w: \hat{\mathbb{C}}_w\rightarrow\hat{\mathbb{C}}_w$ and $P_b: \hat{\mathbb{C}}_b\rightarrow\hat{\mathbb{C}}_b$ 
\end{center}
be two monic polynomials of the same degree $d\geq 2$ with connected and locally connected Julia sets. (One can interchange the index 'w' and 'b' throughout the section, which imposes no essential affect on the results.) Let $\phi_w$ and $\phi_b$ be their B\"ottcher's coordinates respectively (the coordinates extend to the boundaries). Let $\sim_F$ be the equivalent relationship on $\hat{\mathbb{C}}_w\sqcup\hat{\mathbb{C}}_b$ such that 
\begin{center}
$\infty e^{2\pi it}\sim_F \infty e^{-2\pi it}$
\end{center}
for $\infty e^{2\pi it}\in\hat{\mathbb{C}}_w, \infty e^{-2\pi it}\in\hat{\mathbb{C}}_b$ and $t\in I$.

\begin{defn}
The formal mating of the two  degree-$d$ monic polynomials $P_w$ and $P_b$ is defined to be the map
\begin{center}
$P_w$ $\mathrm{\rotatebox[origin=c]{90}{$\vDash$}}_\mathrm{F}$ $P_b: (\hat{\mathbb{C}}_w\sqcup\hat{\mathbb{C}}_b)/ \sim_F\rightarrow (\hat{\mathbb{C}}_w\sqcup\hat{\mathbb{C}}_b)/ \sim_F$
\end{center} 
with $(P_w$ $\mathrm{\rotatebox[origin=c]{90}{$\vDash$}}_\mathrm{F}$ $P_b)|_{\hat{\mathbb{C}}_w}=P_w$ and $(P_w$ $\mathrm{\rotatebox[origin=c]{90}{$\vDash$}}_\mathrm{F}$ $P_b)|_{\hat{\mathbb{C}}_b}=P_b$.
\end{defn}
Note that the map $P_w$ $\mathrm{\rotatebox[origin=c]{90}{$\vDash$}}_\mathrm{F}$ $P_b$ is well-defined  on $\partial \hat{\mathbb{C}}_w=\partial \hat{\mathbb{C}}_b$ according to (\ref{eq15}).  The classical mating can be induced from the  
formal mating by deeper quotient on the sphere $\hat{\mathbb{C}}_w\sqcup\hat{\mathbb{C}}_b$. To do this, assume $\mathcal{J}(P_w)$ and $\mathcal{J}(P_b)$ are connected and locally connected now.  Let $\sim$ be the equivalent relationship on $\hat{\mathbb{C}}_w\sqcup\hat{\mathbb{C}}_b$ such that
\begin{center}
$\phi_w^{-1}(r_w e^{2\pi it})\sim \phi_b^{-1}(r_b e^{-2\pi it})$
\end{center}
for $\phi_w^{-1}(r_w e^{2\pi it})\in\hat{\mathbb{C}}_w,\ \phi_b^{-1}(r_b e^{-2\pi it})\in\hat{\mathbb{C}}_b$ and $r_w, r_b \in[1,\infty], t\in I$. We can judge whether the quotient space $(\hat{\mathbb{C}}_w\sqcup\hat{\mathbb{C}}_b)/ \sim$ is a sphere by the Moore's Theorem (\cite{Moo}). If $(\hat{\mathbb{C}}_w\sqcup\hat{\mathbb{C}}_b)/ \sim$ is a sphere and the degenerate  formal (topological) mating 
\begin{center}
$P_w$ $\rotatebox[origin=c]{90}{$\vDash$}$ $P_b=(P_w$ $\mathrm{\rotatebox[origin=c]{90}{$\vDash$}}_\mathrm{F}$ $P_b)|_{(\hat{\mathbb{C}}_w\sqcup\hat{\mathbb{C}}_b)/ \sim}: (\hat{\mathbb{C}}_w\sqcup\hat{\mathbb{C}}_b)/ \sim\rightarrow (\hat{\mathbb{C}}_w\sqcup\hat{\mathbb{C}}_b)/ \sim$
\end{center}
can be viewed as a rational map in some suitable sense, we get the classical mating ($1$-mating) of $P_w$ and $P_b$.

\begin{defn}\label{def3}
For two monic polynomials $P_w$ and $P_b$ of the same degree $d\geq 2$ with connected and locally connected Julia sets, assume $(\hat{\mathbb{C}}_w\sqcup\hat{\mathbb{C}}_b)/ \sim$ is a topological sphere. Now if there exists some rational map $R:\mathbb{P}^1(\mathbb{C})\rightarrow\mathbb{P}^1(\mathbb{C})$ such that 
\begin{center}
$R\stackrel{top}{\approx}P_w$ $\rotatebox[origin=c]{90}{$\vDash$}$ $P_b$,
\end{center}
we say $P_w$ and $P_b$ are $1$-matable, while $R$ is $1$-unmatable.  We will also say $R$ arises as (or is) a $1$-mating of $P_w$ and $P_b$.
\end{defn}

One can also use $\stackrel{con}{\approx}$ or $\stackrel{Thu}{\approx}$ instead of $\stackrel{top}{\approx}$ in Definition \ref{def3}, which may lead to a stronger or weaker concept in some situations. Upon the classical mating we introduce the $n$\emph{-mating} of polynomials.

\begin{defn}\label{def2}
For some integer $n\geq 1$, a rational map $R: \mathbb{P}^1(\mathbb{C})\rightarrow\mathbb{P}^1(\mathbb{C})$ of degree $d\geq 2$ is said to be  (or arise as) an $n$-mating  if the following two conditions are satisfied.  
\begin{enumerate}[(1).]
\item There exists two monic polynomials $P_w$ and $P_b$ of degree $d^n$, such that $R^n$ arises as the $1$-mating of $P_w$ and $P_b$.

\item In case $n\geq 2$, for any $1\leq m<n$, $R^m$ does not arise as a $1$-mating.
\end{enumerate}

\end{defn}

The integer $n$ is called the \emph{fold} of the (un)mating, which will be denoted by $F(R)$ in the following.  Thus the higher-fold matings are supposed to be of essential distinctions on their dynamics contrasting the classical matings, which never degenerate into the classical ones. From now on we will simply say that a rational map $R$ arises as (or is) a \emph{mating} if the exists some $n\geq 1$ such that $R$ arises as an $n$-mating.  The following definition is dual to Definition  \ref{def2}.

\begin{defn}
For two monic polynomials $P_w$ and $P_b$ of degree $d\geq 2$, if there exists some rational map $R:\mathbb{P}^1(\mathbb{C})\rightarrow\mathbb{P}^1(\mathbb{C})$ which arises as an $n$-mating of $P_w$ and $P_b$ for some $n\in\mathbb{N}_+$, then $P_w$ and $P_b$ are called $n$-matable. 
\end{defn}

Below is a travail example of $n$-mating.
\begin{exa}
The two postcritically finite polynomials $z^4$ and $z^4-2z^2$ are both $1$-matable and $2$-matable since
\begin{center}
$(z^2-1)^{\circ 2}=(z^2-1)\circ(z^2-1)$ $\stackrel{top}{\approx}$ $z^4$\rotatebox[origin=c]{90}{$\vDash$}$(z^4-2z^2)$.
\end{center}
\end{exa} 

The folds of matings of polynomials admit some order structure. 
\begin{rem}
Any two $(n_1n_2)$-matable polynomials $P_w$ and $P_b$ for two positive integers $n_1, n_2$ must be both $n_1$-matable and $n_2$-matable.
\end{rem}
If two monic polynomials $P_w$ and $P_b$ of degree $d\geq 2$ are matable, it would be interesting to ask the largest fold $n$ such that they are $n$-matable. It is necessary a logarithmic integer of $d$ with respect to some integer base.

\subsection{Equators and orientation-reversing equators}

To investigate the $1$-unmatability of rational maps, Wittner (\cite{Wit}) introduced the notion of \emph{equator} for rational maps, which can be employed to split the dynamics of the rational map into dynamics of two polynomials. One can also refer to \cite[Section 4]{Mey1}. A Jordan curve $\Xi\subset \mathbb{P}^1(\mathbb{C})$ is said to be \emph{circular} to a set $A\subset \mathbb{P}^1(\mathbb{C})$ if $A$ is contained in either one of the two connected components of $\mathbb{P}^1(\mathbb{C})\setminus \Xi$.
\begin{defn}\label{def8}
For a rational map $R:\mathbb{P}^1(\mathbb{C})\rightarrow\mathbb{P}^1(\mathbb{C})$, a Jordan curve $\Xi\subset \mathbb{P}^1(\mathbb{C})\setminus P(R)$ is called an \emph{equator} of $R$, if $R^{-1}(\Xi)$ is connected and orientation-preserving isotopic to $\Xi$ \emph{rel.} $\mathcal{P}(R)$ in $\mathbb{P}^1(\mathbb{C})$.    
\end{defn}
We traditionally use $\{H_t: \mathbb{P}^1(\mathbb{C})\rightarrow\mathbb{P}^1(\mathbb{C})\}_{t\in I}$ to denote the isotopy from $\Xi$ to $R^{-1}(\Xi)$. The orientation-preserving requirement in Definition \ref{def8} means that the two orientations induced by the map $R$ as well as the end of the isotopy $H_1$ on $\Xi$ and $R^{-1}(\Xi)$ coincide with each other. 

The two closed components of $\mathbb{P}^1(\mathbb{C})$ separated by a Jordan curve $\Xi$ are denoted by $U_w^0$ and $U_b^0$ respectively, while their pre-images under $R$ are denoted by $U_w^1$ and $U_b^1$ respectively. This means that 
\begin{equation}\label{eq3}
R(U_w^1)=U_w^0, R(U_b^1)=U_b^0.
\end{equation}
In case $\Xi$ is an equator of $R$, $U_w^1, U_b^1$ are two connected components and we have 
\begin{equation}\label{eq4}
H_1(U_w^0)=U_w^1, H_1(U_b^0)=U_b^1.
\end{equation}

\begin{figure}[H]
\centering
\def\svgwidth{12cm}
\begingroup%
  \makeatletter%
  \providecommand\color[2][]{%
    \errmessage{(Inkscape) Color is used for the text in Inkscape, but the package 'color.sty' is not loaded}%
    \renewcommand\color[2][]{}%
  }%
  \providecommand\transparent[1]{%
    \errmessage{(Inkscape) Transparency is used (non-zero) for the text in Inkscape, but the package 'transparent.sty' is not loaded}%
    \renewcommand\transparent[1]{}%
  }%
  \providecommand\rotatebox[2]{#2}%
  \newcommand*\fsize{\dimexpr\f@size pt\relax}%
  \newcommand*\lineheight[1]{\fontsize{\fsize}{#1\fsize}\selectfont}%
  \ifx\svgwidth\undefined%
    \setlength{\unitlength}{1102.97507174bp}%
    \ifx\svgscale\undefined%
      \relax%
    \else%
      \setlength{\unitlength}{\unitlength * \real{\svgscale}}%
    \fi%
  \else%
    \setlength{\unitlength}{\svgwidth}%
  \fi%
  \global\let\svgwidth\undefined%
  \global\let\svgscale\undefined%
  \makeatother%
  \begin{picture}(1,0.59585398)%
    \lineheight{1}%
    \setlength\tabcolsep{0pt}%
    \put(0,0){\includegraphics[width=\unitlength,page=1]{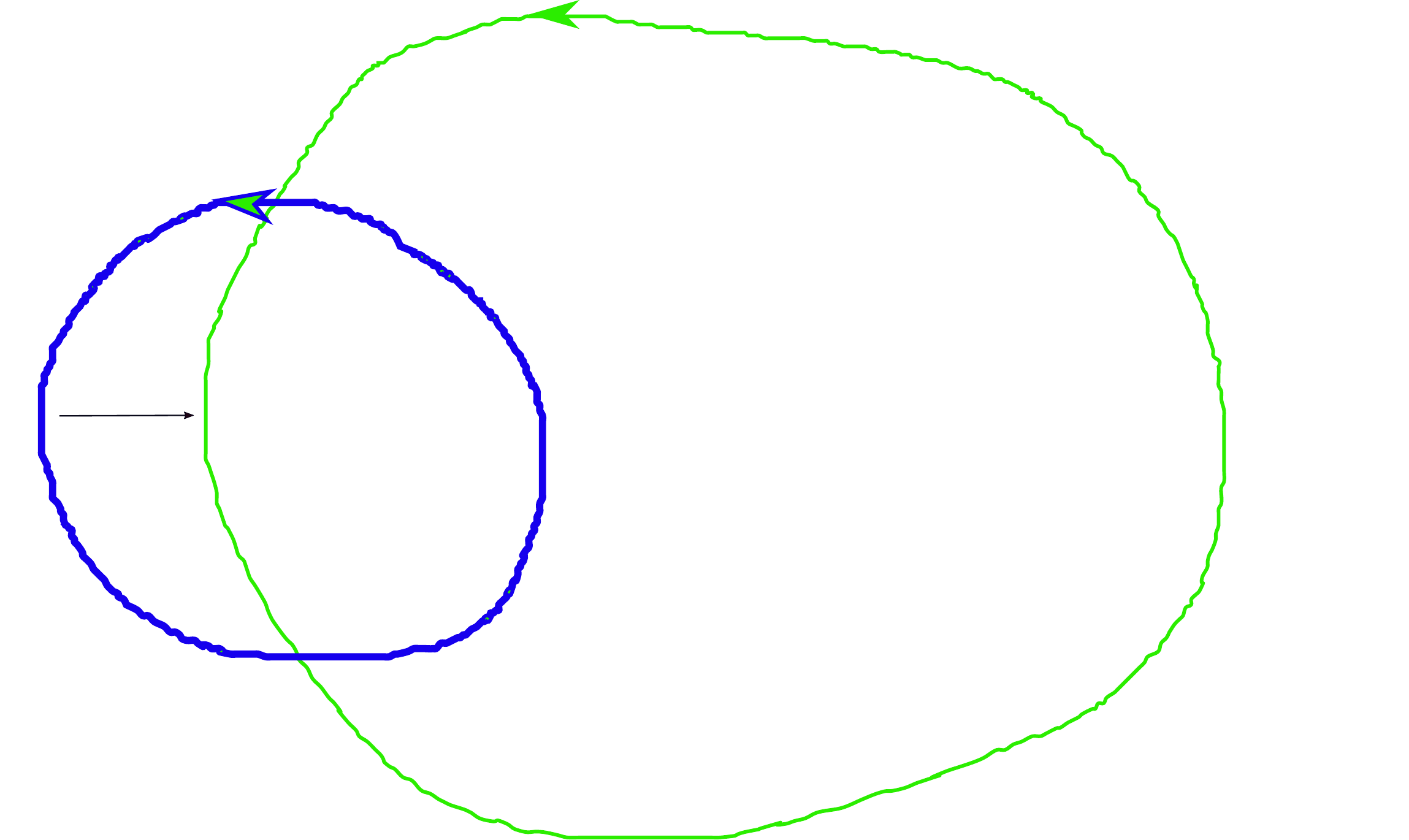}}%
    \put(-0.00205853,0.29441041){\color[rgb]{0,0,0}\makebox(0,0)[lt]{\lineheight{0}\smash{\begin{tabular}[t]{l}$\Xi$\end{tabular}}}}%
    \put(0.1520145,0.29461629){\color[rgb]{0,0,0}\makebox(0,0)[lt]{\lineheight{0}\smash{\begin{tabular}[t]{l}$R^{-1}(\Xi)$\end{tabular}}}}%
    \put(0.39684483,0.28663924){\color[rgb]{0,0,0}\makebox(0,0)[lt]{\lineheight{0}\smash{\begin{tabular}[t]{l}$U_b^0$\end{tabular}}}}%
    \put(0.33147275,0.29052486){\color[rgb]{0,0,0}\makebox(0,0)[lt]{\lineheight{0}\smash{\begin{tabular}[t]{l}$U_w^0$\end{tabular}}}}%
    \put(0.87351326,0.28663916){\color[rgb]{0,0,0}\makebox(0,0)[lt]{\lineheight{0}\smash{\begin{tabular}[t]{l}$U_b^1$\end{tabular}}}}%
    \put(0.81911488,0.2866392){\color[rgb]{0,0,0}\makebox(0,0)[lt]{\lineheight{0}\smash{\begin{tabular}[t]{l}$U_w^1$\end{tabular}}}}%
    \put(0.07239768,0.31772397){\color[rgb]{0,0,0}\makebox(0,0)[lt]{\lineheight{0}\smash{\begin{tabular}[t]{l}$H_t$\end{tabular}}}}%
  \end{picture}%
\endgroup%

\caption{Topological dynamics of a sample map $R$ with an equator $\Xi$}
\label{fig2}
\end{figure}
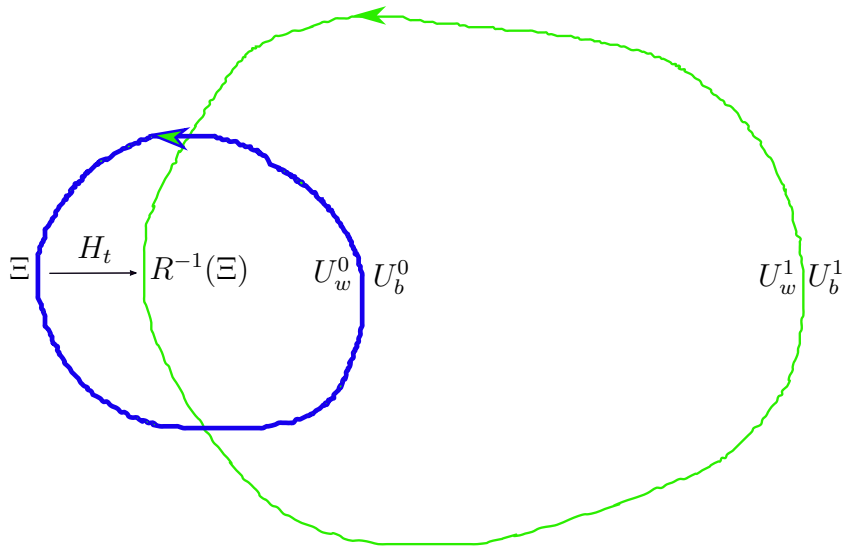

An exhibition of the topological dynamics of some sample rational map $R$ with an equator $\Xi$ is presented in Figure  \ref{fig2}.

Throughout the paper the readers will intermittently see pictures of Jordan curves and their pre-images on the dynamical planes of some maps. We may mark them in some pictures, but not always. In case of no markings we always colour an original Jordan curve $\Xi$ (may not be an equator or an OR equator) blue, its first pre-image  $R^{-1}(\Xi)$ green, its second pre-image  $R^{-2}(\Xi)$ black, its third pre-image  $R^{-3}(\Xi)$ yellow and its fourth pre-image  $R^{-4}(\Xi)$ purple on the the dynamical plane of a map $R$.

Let $P_w(R)=U_w^0\cap P(R)$ and $P_b(R)=U_b^0\cap P(R)$ be the two subsets of $P(R)$ separated by a Jordan curve $\Xi\subset \mathbb{P}^1(\mathbb{C})\setminus P(R)$. Then $\{P_w(R), P_b(R)\}$ is a (non-travail) partition of $P(R)$ if neither of them is empty. The two subsets are \emph{immune} to each other in the following sense if $\Xi$ is an equator. 

\begin{lemma}\label{lem5}
For a rational map $R$ with an equator $\Xi$, the two sets $P_w(R)$ and $P_b(R)$ are immune to each other under $R$, that is, 
\begin{equation}\label{eq36}
R(P_w(R))\subset P_w(R)\subset R^{-1}(P_w(R))
\end{equation}
and
\begin{equation}\label{eq37}
R(P_b(R))\subset P_b(R)\subset R^{-1}(P_b(R)).
\end{equation}
\end{lemma}
\begin{proof}
First note that $U_w^0\cap U_w^1\neq \emptyset$. It follows from (\ref{eq4}) that  
$$P_w(R)\subset U_w^0\cap U_w^1.$$
Considering (\ref{eq3}) this implies  $P_w(R)\subset R^{-1}(P_w(R))$, which gives rise to (\ref{eq36}). The relationship (\ref{eq37}) follows in a similar way. 
\end{proof}

The simple but important property will be used frequently in searching equators for rational maps in the following. Note that the $\subset$ can not be replaced by $=$ in (\ref{eq36}) and (\ref{eq37}), since there might be critical but not postcritical points in $U_w^1\setminus U_w^0$. In case 
\begin{center}
$R(P_w(R))=P_w(R)=R^{-1}(P_w(R))$ and $R(P_b(R))=P_b(R)=R^{-1}(P_b(R))$
\end{center}
hold we say the two sets $P_w(R)$ and $P_b(R)$ are \emph{strictly immune} to each other.

Now we introduce the \emph{orientation-reversing (OR) equators} for rational maps.  Meyer (\cite[Section 12]{Mey1}) first discovered the existence of some orientation-reversing (pseudo-)equators for some rational maps in unmating certain rational maps.    
\begin{defn}\label{def4}
For a rational map $R: \mathbb{P}^1(\mathbb{C})\rightarrow\mathbb{P}^1(\mathbb{C})$, a Jordan curve $\Xi\subset \mathbb{P}^1(\mathbb{C})\setminus P(R)$ is called an \emph{OR equator} of $R$, if $R^{-1}(\Xi)$ is connected and orientation-reversing isotopic to $\Xi$ \emph{rel.} $P(R)$.    
\end{defn}
This means that the two orientations induced by the map $R$ as well as the end of the isotopy $H_1$ for the OR equator on $\Xi$ and $R^{-1}(\Xi)$ opposite each other. There are usually apriori properties for the OR equators dual to those of equators, but not always. Let $\Xi$ be a Jordan curve along with the sets $U_w^0, U_b^0$, $U_w^1, U_b^1$ satisfying (\ref{eq3}). Now if  $\Xi$ is an OR equator with the isotopy $\{H_t: \mathbb{P}^1(\mathbb{C})\rightarrow\mathbb{P}^1(\mathbb{C})\}_{t\in I}$ between $\Xi$ and $R^{-1}(\Xi)$, then the following equations dual to (\ref{eq4}) hold,
\begin{equation}\label{eq14}
H_1(U_w^0)=U_b^1, H_1(U_b^0)=U_w^1.
\end{equation}

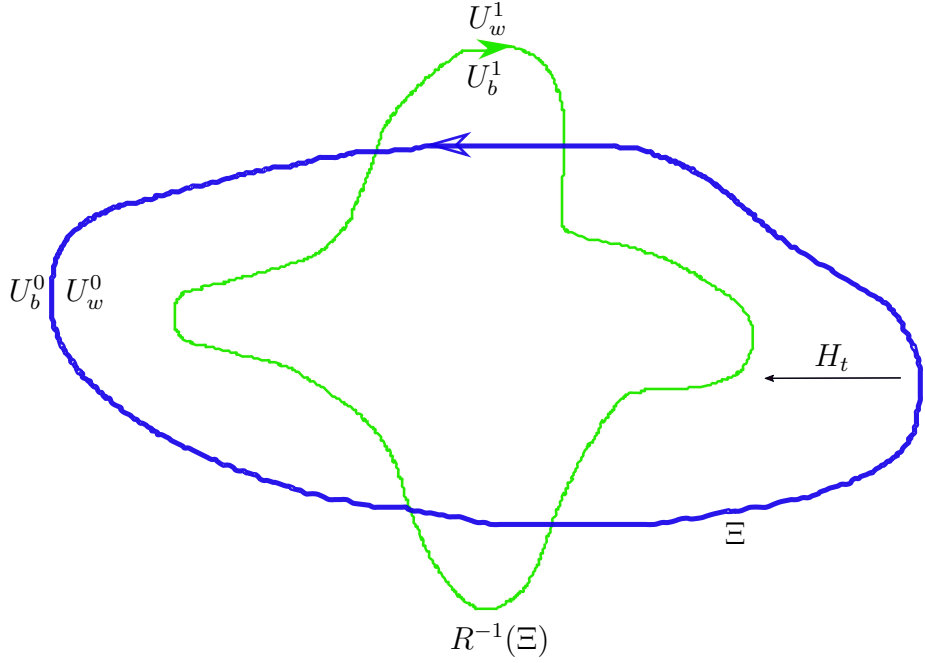
\begin{figure}[ht]
\centering
\def\svgwidth{12cm}
\begingroup%
  \makeatletter%
  \providecommand\color[2][]{%
    \errmessage{(Inkscape) Color is used for the text in Inkscape, but the package 'color.sty' is not loaded}%
    \renewcommand\color[2][]{}%
  }%
  \providecommand\transparent[1]{%
    \errmessage{(Inkscape) Transparency is used (non-zero) for the text in Inkscape, but the package 'transparent.sty' is not loaded}%
    \renewcommand\transparent[1]{}%
  }%
  \providecommand\rotatebox[2]{#2}%
  \newcommand*\fsize{\dimexpr\f@size pt\relax}%
  \newcommand*\lineheight[1]{\fontsize{\fsize}{#1\fsize}\selectfont}%
  \ifx\svgwidth\undefined%
    \setlength{\unitlength}{1014.69063699bp}%
    \ifx\svgscale\undefined%
      \relax%
    \else%
      \setlength{\unitlength}{\unitlength * \real{\svgscale}}%
    \fi%
  \else%
    \setlength{\unitlength}{\svgwidth}%
  \fi%
  \global\let\svgwidth\undefined%
  \global\let\svgscale\undefined%
  \makeatother%
  \begin{picture}(1,0.71482143)%
    \lineheight{1}%
    \setlength\tabcolsep{0pt}%
    \put(0,0){\includegraphics[width=\unitlength,page=1]{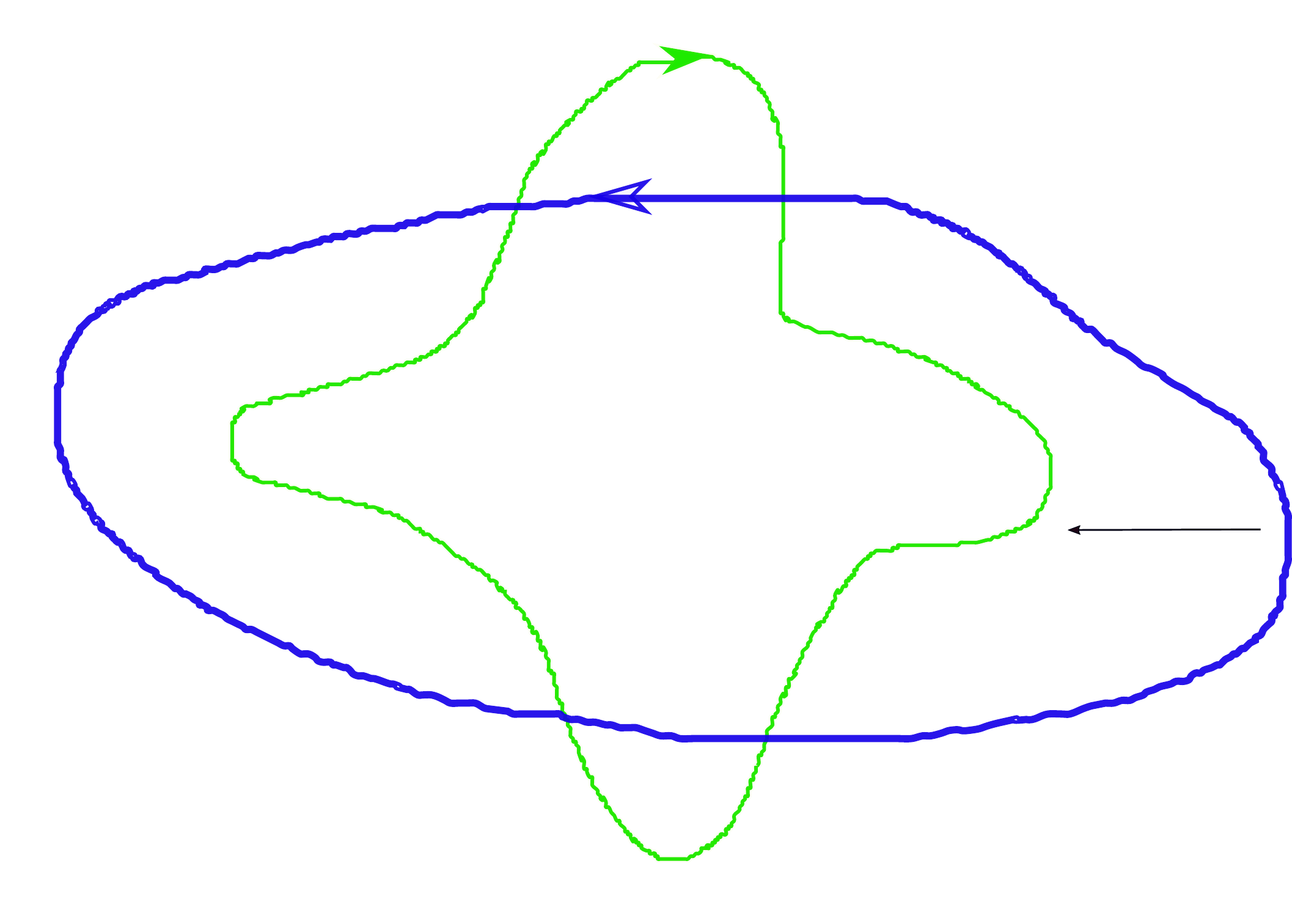}}%
    \put(0.78267956,0.12187483){\color[rgb]{0,0,0}\makebox(0,0)[lt]{\lineheight{0}\smash{\begin{tabular}[t]{l}$\Xi$\end{tabular}}}}%
    \put(0.48133095,0.00594778){\color[rgb]{0,0,0}\makebox(0,0)[lt]{\lineheight{0}\smash{\begin{tabular}[t]{l}$R^{-1}(\Xi)$\end{tabular}}}}%
    \put(-0.00223764,0.39007766){\color[rgb]{0,0,0}\makebox(0,0)[lt]{\lineheight{0}\smash{\begin{tabular}[t]{l}$U_b^0$\end{tabular}}}}%
    \put(0.05974801,0.3900777){\color[rgb]{0,0,0}\makebox(0,0)[lt]{\lineheight{0}\smash{\begin{tabular}[t]{l}$U_w^0$\end{tabular}}}}%
    \put(0.49900925,0.62660286){\color[rgb]{0,0,0}\makebox(0,0)[lt]{\lineheight{0}\smash{\begin{tabular}[t]{l}$U_b^1$\end{tabular}}}}%
    \put(0.50112108,0.69206973){\color[rgb]{0,0,0}\makebox(0,0)[lt]{\lineheight{0}\smash{\begin{tabular}[t]{l}$U_w^1$\end{tabular}}}}%
    \put(0.88050862,0.31827539){\color[rgb]{0,0,0}\makebox(0,0)[lt]{\lineheight{0}\smash{\begin{tabular}[t]{l}$H_t$\end{tabular}}}}%
  \end{picture}%
\endgroup%

\caption{Topological dynamics of a sample map $R$ with an OR equator-I}
\label{fig3}
\end{figure}

We illustrates the topological dynamics of some sample rational map $R$ with an OR equator $\Xi$ in Figure  \ref{fig3}.

In case $\Xi$ is an OR equator for $R$, the two sets in the partition $\{P_w(R), P_b(R)\}$ of $P(R)$ satisfy the following \emph{swapping property} under $R$, dual to Lemma \ref{lem5}.

\begin{lemma}\label{lem7}
The two sets $P_w(R)$ and $P_b(R)$ are swapped to each other under $R$, that is, 
\begin{equation}\label{eq1}
R(P_w(R))\subset P_b(R)\subset R^{-1}(P_w(R))
\end{equation}
and
\begin{equation}\label{eq2}
R(P_b(R))\subset P_w(R)\subset R^{-1}(P_b(R)).
\end{equation}
\end{lemma}
\begin{proof}
These formulas follow from similar arguments as in the proof of Lemma \ref{lem5}, in virtue of (\ref{eq14}).  
\end{proof}

In case the equalities hold in (\ref{eq1}) and (\ref{eq2}), we say $P_w(R)$ and $P_b(R)$ are \emph{strictly swapping} to each other.

For a hyperbolic postcritically finite rational map $R$ with an equator $\Xi$, let 
$$\tilde{R}_w=H_1\circ R: U_w^1/R^{-1}(\Xi)\rightarrow U_w^1/R^{-1}(\Xi)$$
be the white map on the quotient space, which is Thurston equivalent to some polynomial $P_w$. Let 
$$\tilde{R}_b=H_1\circ R: U_b^1/R^{-1}(\Xi)\rightarrow U_b^1/R^{-1}(\Xi)$$ 
be the black map, which is Thurston equivalent to some polynomial $P_b$.
For a hyperbolic postcritically finite rational map $R$ with an OR equator $\Xi$, let 
$$\tilde{R}_{wb}=H_1\circ R: U_w^1/R^{-1}(\Xi)\rightarrow U_b^1/R^{-1}(\Xi)$$
be the white-black map, let 
$$\tilde{R}_{bw}=H_1\circ R: U_b^1/R^{-1}(\Xi)\rightarrow U_w^1/R^{-1}(\Xi)$$ 
be the black-white map. These maps are degenerates of the maps $T_{wb}, T_{bw}$ in Section \ref{sec4}. The following result describes the backward limit behaviour of equators and OR equators of rational maps. 

\begin{lemma}\label{lem6}
Let $R$ be a hyperbolic postcritically finite rational map with an equator (or OR equator) $\Xi$. Considering the backward iterations of $\Xi$ under the map $R$, the limit set 
$$\Xi_\infty=\lim_{k\rightarrow\infty} R^{-k}(\Xi)$$ 
satisfies
\begin{equation}\label{eq7}
\Xi_\infty=\mathcal{J}(R),
\end{equation}
in which 
\begin{center}
$
\begin{array}{ll}
& \lim_{k\rightarrow\infty} R^{-k}(\Xi)\vspace{3mm}\\
=&\big\{z\in\mathbb{P}^1(\mathbb{C}): \mbox{ there exists a sequence } \{z_k\in R^{-k}(\Xi)\}_{k\in\mathbb{N}} \mbox{ such that } \lim_{k\rightarrow\infty}z_k=z\big\}.
\end{array}
$
\end{center}
\end{lemma}
\begin{proof}
First note that the restriction $R: R^{-k}(\Xi)\rightarrow R^{-k+1}(\Xi)$ is always a degree-$d$ covering map for any $k\geq 0$. We claim that 
\begin{equation}\label{eq6}
\Xi_\infty\cap \mathcal{J}(R)\neq \emptyset.
\end{equation}

We justify (\ref{eq6}) in two cases. If $\Xi\cap \mathcal{J}(R)\neq \emptyset$, say, $z_0\in \Xi\cap \mathcal{J}(R)$, then $R^{-k}(z_0)\subset R^{-k}(\Xi)\cap \mathcal{J}(R)$ for any $k\geq 1$, which results in at least one limit point necessarily in  $\mathcal{J}(R)$ since $\mathcal{J}(R)$ is perfect.

Now we consider the case $\Xi\cap \mathcal{J}(R)= \emptyset$, so $\Xi$ is confined in some single Fatou component of $R$. In case $\Xi$ is an equator, this happens if and only if $\#P_w(R)=1$ (or $\#P_b(R)=1$), that is, there is a superattracting completely invariant fixed point $z_\infty$ of $R$ in $U_w^0$. Moreover, $\Xi$ must be circular to the superattracting fixed point and contained in the interior of the basin $U_\infty$ of $z_\infty$. In case $\Xi$ is an OR equator, this happens if and only if $\#P_w(R)=\#P_b(R)=1$, that is, there is a superattracting completely invariant fixed point $z_\infty$ of $R^2$ in $U_w^0$.  In both cases one can see that 
$$\lim_{k\rightarrow\infty} R^{-k}(\Xi)=\partial U_\infty$$ 
in virtue of the B\"ottcher's Theorem, which gives rise to (\ref{eq6}).

Note that $\Xi_\infty$ is a completely invariant set, that is,
$$R(\Xi_\infty)=\Xi_\infty=R^{-1}(\Xi_\infty).$$
Obviously $\Xi_\infty\neq \mathbb{P}^1(\mathbb{C})$, together with the claim above, we get (\ref{eq7}).
\end{proof}

In fact Lemma \ref{lem6} holds for equators or OR equators of any postcritically finite rational map. The limit set corresponds to the boundary of the filled Julia set $\partial \mathcal{K}_b$ of the degree-$d$ polynomial $P_b$ Thurston equivalent to the black map $\tilde{R}_b=H_1\circ R: U_b^1/R^{-1}(\Xi)\rightarrow U_b^1/R^{-1}(\Xi)$ in case $\Xi$ is an equator. Lemma \ref{lem6} implies the following topological description of the distribution of an equator (or OR equator) with respect to its pre-image for rational maps. 

\begin{pro}\label{pro8}
For a postcritically finite hyperbolic rational map $R$ with an equator $\Xi$, we have
\begin{equation}\label{eq8}
U_w^0\setminus U_w^1\neq\emptyset
\end{equation}
and
\begin{equation}\label{eq9}
U_w^1\setminus U_w^0\neq\emptyset
\end{equation}
hold simultaneously unless $R$ is a Thurston polynomial.  For a postcritically finite hyperbolic rational map $R$ with an OR equator $\Xi$, we have (\ref{eq8}) and (\ref{eq9})
hold simultaneously unless $R$ is conformally conjugate to $\cfrac{1}{z^2}$.
\end{pro}

\begin{proof}
We deal with the equator case first.  Suppose $R$ is not a Thurston polynomial (so it does not admit any completely invariant superattracting fixed point), we will show that (\ref{eq8}) and (\ref{eq9}) must hold simultaneously. Note that in this case since $\Xi$ is circular to $P_w(R)$ and $P_b(R)$ subjecting to the immune property, it must be that 
\begin{equation}\label{eq10}
(U_w^0)^o\cap \mathcal{J}(R)\neq\emptyset
\end{equation}
and 
\begin{equation}\label{eq11}
(U_b^0)^o\cap \mathcal{J}(R)\neq\emptyset
\end{equation}
hold simultaneously. Now if either of (\ref{eq8}) and (\ref{eq9}) fails, without loss of generality we suppose $U_w^1\subset U_w^0$. Then inductively one can see that 
$$R^{-k}(\Xi)\subset U_w^1\subset U_w^0$$
for any $k\geq 1$,  which results in that $\Xi_\infty\subset U_w^0$.
Now combining (\ref{eq7}), (\ref{eq10}) and (\ref{eq11}) together we get some contradiction. 

In case $R$ admitting an OR equator, if $R$ does not admit a completely invariant superattracting periodic cycle of period $2$ (so $R^2$ does not admit a completely invariant superattracting fixed point), by similar steps, one can deduce that $\Xi_\infty$ will be confined in some subsets ($U_w^0$ or $U_b^0$) which contradicting (\ref{eq7}). If $R$ admits a completely invariant superattracting periodic cycle of period $2$ and some OR equator $\Xi$ satisfying  (\ref{eq8}) and (\ref{eq9}) simultaneously, let $\{z_0, z_\infty\}$ be the periodic cycle with $z_0=R^2(z_\infty)=R^{-2}(z_\infty)$. Now let $U_0, U_\infty$ be the immediate basin of $z_0, z_\infty$ respectively. In this case without loss of generality we can assume $\Xi\subset U_0$ and $R^{-1}(\Xi)\subset U_\infty$. Note that there are no critical or postcritical points in the annulus bounded by $\Xi$ and $R^{-1}(\Xi)$, so there are not other Fatou components besides $U_0, U_\infty$ for $R$, which forces $R$ to be conformally conjugate to $\cfrac{1}{z^2}$. 
\end{proof}

Proposition \ref{pro8} denies the hope that one can isotope an equator (OR equator) of some rational map $R$ to some equator (OR equator) which confines its pre-image in general cases \emph{rel.} $P(R)$. This indicates complexity of the unmating process.

\section{Classification of hyperbolic $2$-matings}\label{sec3}

\subsection{OR equators switch to equators of the second iterate}

We mean to prove Theorem \ref{thm1} in this subsection. We first point out that equators and OR equators do not resist each other, which means that there exists some rational map $R$ admitting an equator and an OR equator simultaneously.  A hyperbolic rational map admitting an equator and an OR equator simultaneously is called a \emph{hermaphroditic mating}. One can find such maps in the cubic bicritical family, see for instance Proposition \ref{pro5}.

We will use the simplified version of the Baer-Epstein-Zieschang Theorem on certain subsets of the Riemann sphere in the following, see \cite{Bae, Eps, Zie} as well as \cite[Appendix]{Bus} for versions on two dimensional compact connected manifolds with piecewise smooth boundaries.

\begin{BEZ}
Let $M\subset \mathbb{P}^1(\mathbb{C})$ be a compact connected manifold with smooth boundary. Let $\Xi_1, \Xi_2\subset M^o$ be two non-travail Jordan curves homotopic to each other in $M$. Then there exists some isotopy $\{G_t\}_{0\leq t\leq 1}$ from $\Xi_1$ to $\Xi_2$ \emph{rel.} $\mathbb{P}^1(\mathbb{C})\setminus M^o$ in $\mathbb{P}^1(\mathbb{C})$.
\end{BEZ}

\begin{proof}
Considering the two non-travail Jordan curves $\Xi_1$ and $\Xi_2$ homotopic to each other in the smooth compact connected two dimensional manifold $M$ with smooth boundary, apply \cite[Theorem A.3]{Bus} in our case, we get a homeomorphism $G_1^*: M\rightarrow M$ such that
\begin{itemize}
\item $G_1^*(\Xi_1)=\Xi_2$.

\item $\mathds{1}_{M}$ is isotopic to $G_1^*$ through an isotopy $\{G_t^*\}_{t\in [0,1]}$ on $M$.

\item The restriction $G_t^*|_{\partial M}=\mathds{1}_{\partial M}$ for any $t\in [0,1]$.
\end{itemize}
Then the extension
$$G_t(z)=\left\{
\begin{array}{lll}
G_t^*(z) & for & z\in M^o,\\
z & for & z\in \mathbb{P}^1(\mathbb{C})\setminus M^o
\end{array}
\right.$$
gives an isotopy from $\Xi_1$ to $\Xi_2$ \emph{rel.} $\mathbb{P}^1(\mathbb{C})\setminus M^o$ in $\mathbb{P}^1(\mathbb{C})$.
\end{proof}

The following result is crucial in our proof of Theorem \ref{thm1}.

\begin{lemma}\label{lem4}
Let $R:\mathbb{P}^1(\mathbb{C})\rightarrow\mathbb{P}^1(\mathbb{C})$ be a rational map.  If $\Xi$ is an OR equator of $R$, then it is an equator of $R^2$.
\end{lemma}
\begin{proof}
As $\Xi$ is an OR equator of $R$, there exists an isotopy $\{F_t\}_{t\in[0,1]}$ from $\Xi$ to $R^{-1}(\Xi)$ \emph{rel.} $P(R)$ in $\mathbb{P}^1(\mathbb{C})$. Since $\Xi\cap P(R)=\emptyset$, we have $R^{-1}(\Xi)\cap P(R)=\emptyset$, which forces
\begin{center}
$R^{-2}(\Xi)\cap C(R)=\emptyset$.
\end{center}
This guarantees $R^{-2}$ is a covering map without branching points from $R^{-2}(\Xi)$ to $\Xi$. We use the notation $U_w^0, U_b^0, U_w^1, U_b^1$ with respect to the OR equator $\Xi$ satisfying (\ref{eq3}) and $P_w(R), P_b(R)$ satisfying (\ref{eq1}), (\ref{eq2}).

We claim that $R^{-2}(\Xi)$ is connected. Suppose this is not the case,  $R^{-2}(\Xi)$ splits into $m$ connected components $\{\Xi_j\}_{j=1}^m$ for some integer $m\geq 2$. Obviously $\Xi_{j_1}\cap \Xi_{j_2}=\emptyset$, moreover, they can not be isotopic to each other \emph{rel.} $P(R)$ in $\mathbb{P}^1(\mathbb{C})$ for any $1\leq j_1<j_2\leq m$. Then $\Xi_j$ can not be circular to $P_w(R)$ and $P_b(R)$, in other words,  $\Xi_j$ must be circular to a set containing points in both $P_w(R)$ and $P_b(R)$ simultaneously  for some $1\leq j\leq m$ (in case $\Xi_j$ is circular to $P_w(R)\cup P_b(R)$ it is a travail Jordan curve).  In any case these conditions  ruin the swapping property which the partition $\{P_w(R), P_b(R)\}$ satisfies since $\Xi$ is an OR equator. 

We denote the two closed components of $\mathbb{P}^1(\mathbb{C})$ separated by $R^{-2}(\Xi)$ by $U_w^2$ and $U_b^2$ respectively, along with the relationships
\begin{center}
$R(U_w^2)=U_w^1, R(U_b^2)=U_b^1$.
\end{center}
$R^{-2}(\Xi)$ must be circular to $P_w(R)$ (and hence $P_b(R)$) in virtue of Lemma \ref{lem7}. Then we have
$$P_w(R)\subset (U_w^0\cap U_b^1\cap U_w^2)^o$$
and
$$P_b(R)\subset \mathbb{P}^1(\mathbb{C})\setminus(U_w^0\cup U_b^1\cup U_w^2).$$
Since $\#(P_w(R)\cup P_b(R))=\#P(R)=k$ is finite, we can choose disjoint smooth closed discs (closed topological discs with smooth boundaries)  $\mathcal{D}_1, \mathcal{D}_2, \cdots, \mathcal{D}_l$ for some $2\leq l\leq k$, such that 
\begin{itemize}
\item $P_w(R)\cup P_b(R)\subset\cup_{j=1}^l (\mathcal{D}_j)^o$.
\item $\cup_{j=1}^l \mathcal{D}_j\subset (U_w^0\cap U_b^1\cap U_w^2)^o\cup \big(\mathbb{P}^1(\mathbb{C})\setminus(U_w^0\cup U_b^1\cup U_w^2)\big)$.
\item $R^{-1}(\Xi)$ is homotopic to $R^{-2}(\Xi)$ in  $\mathbb{P}^1(\mathbb{C})\setminus\cup_{j=1}^l (\mathcal{D}_j)^o$.
\end{itemize}
In the best case $2$ such smooth closed topological discs will be enough while in the worst case we need $k$ such discs. In case that 
$$R^{-1}(\Xi)\cap R^{-2}(\Xi)=\emptyset,$$
a homotopy is obvious. In case that 
$$R^{-1}(\Xi)\cap R^{-2}(\Xi)\neq\emptyset,$$
a homotopy can be obtained by combining series of homotopies between corresponding arcs respectively in $R^{-1}(\Xi)$ and $ R^{-2}(\Xi)$ joining the intersections (which are in fact basepoints of the homotopies) in $\mathbb{P}^1(\mathbb{C})\setminus\cup_{j=1}^l (\mathcal{D}_j)^o$. Such homotopies exist because if there are postcritical points which prevent the corresponding pairs of arcs to be homotopic to each other, then there will also be postcritical points which prevent $\Xi$ to be isotopic to $R^{-1}(\Xi)$ \emph{rel.} $P(R)$.  Now let $M=\mathbb{P}^1(\mathbb{C})\setminus\cup_{j=1}^l (\mathcal{D}_j)^o$, apply the Baer-Epstein-Zieschang Theorem, there exists an isotopy $\{G_t\}_{t\in[0,1]}$ from $R^{-1}(\Xi)$ to $R^{-2}(\Xi)$ \emph{rel.} $\cup_{j=1}^l \mathcal{D}_j$. The isotopy is also \emph{rel.} $P(R)$ since $P(R)\subset\cup_{j=1}^l (\mathcal{D}_j)^o$.  It has to be orientation-reversing considering Lemma \ref{lem7}.  Finally, let
$$H_t(z)=\left\{
\begin{array}{lll}
F_{2t}(z) & for & 0\leq t\leq 1/2,\\
G_{2t-1}(z) & for & 1/2< t\leq 1
\end{array}
\right.$$
be the combination of the two successive orientation-reversing isotopies $\{F_t\}_{t\in[0,1]}$ and $\{G_t\}_{t\in[0,1]}$. Then $\{H_t\}_{t\in[0,1]}$ gives an orientation-preserving isotopy from $\Xi$ to $R^{-2}(\Xi)$ \emph{rel.} $P(R)$. 
\end{proof}

\begin{figure}[ht]
    \centering
    \def\svgwidth{\columnwidth}
    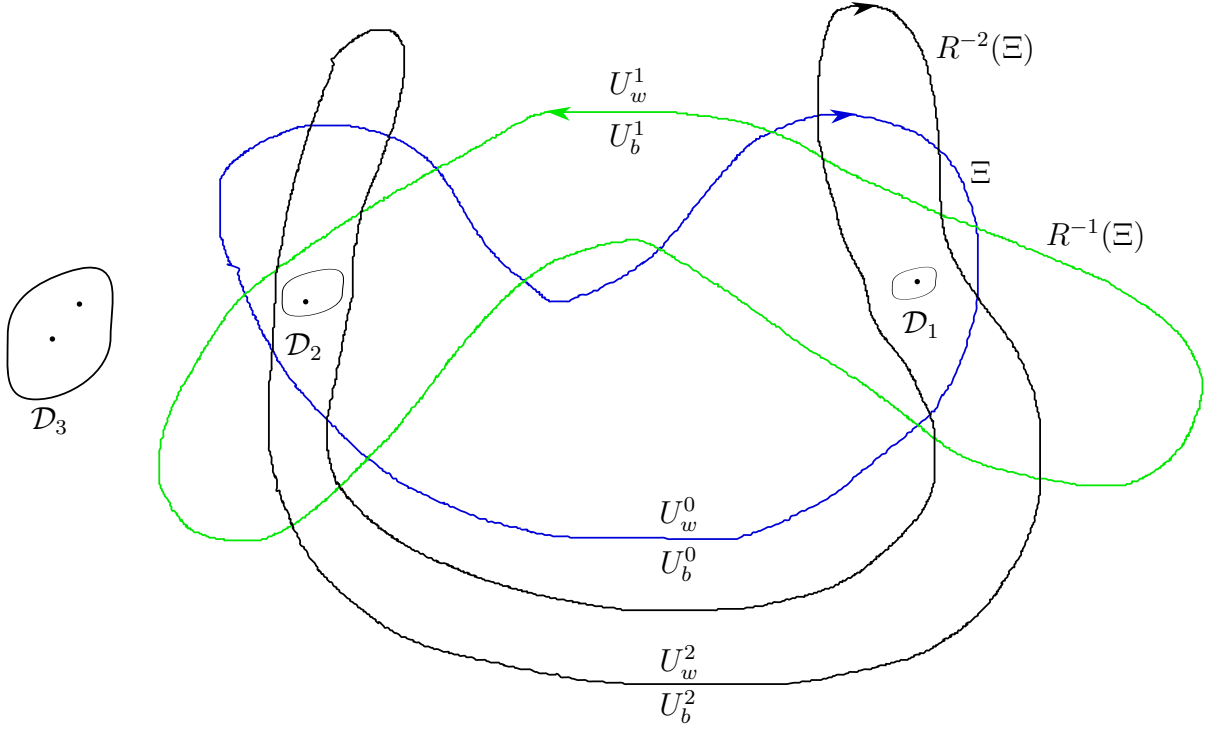
    \caption{Topological dynamics of a sample map $R$ with an OR equator-II}
    \label{fig1}
\end{figure}

In Figure \ref{fig1} we provide the readers a topological illustration of the twice backward iterations of an OR equator of some sample rational map with $\#P(R)=4$. 

\begin{rem}
Lemma \ref{lem4} extends to higher iterates in fact. If $\Xi$ is an OR equator of some rational map $R$, then it is an equator of $R^n$ for even $n\in\mathbb{N}_+$, it is an OR equator of $R^n$ for odd $n\in\mathbb{N}_+$. However, this does not mean that $R$ will arise as an $n$-mating for $n\geq 3$ in case $R$ is hyperbolic, see Corollary \ref{cor7}. 
\end{rem}

A dual result for equators of rational maps on the Riemann sphere holds, following analogous arguments as in proof of Lemma \ref{lem4}. 
\begin{coro}\label{cor1}
Let $R: \mathbb{P}^1(\mathbb{C})\rightarrow\mathbb{P}^1(\mathbb{C})$ be a rational map.  If $\Xi$ is an equator of $R$, then it is an equator of $R^k$ for any $k\in\mathbb{N}_+$.
\end{coro}

Considering  Lemma \ref{lem4} and Corollary \ref{cor1} together, it is interesting to ask whether some inverse conclusion holds for  rational maps whose twice iterate admit an equator. Unfortunately this is wrong. There exists some rational map $R$, such that $R^2$ admits an equator $\Xi$, while $\Xi$ is neither an equator nor an OR equator of $R$, see for example Proposition \ref{pro2}. However, an inverse conclusion for Lemma \ref{lem4} and Corollary \ref{cor1} does hold for some families of rational maps, for example, (Thurston) polynomials.

\begin{pro}\label{pro1}
Let $R:\mathbb{P}^1(\mathbb{C})\rightarrow\mathbb{P}^1(\mathbb{C})$ be a postcritically finite rational map such that $R^2$ is M\"obius conjugate to a polynomial.  If $R^2$ admits an equator $\Xi$, then it is either an equator or an OR equator of $R$.
\end{pro}

Proof is left to the attentive readers.  Equipped with all the above results, now we are in a position to justify Theorem \ref{thm1}. 

\begin{proof}[Proof of Theorem \ref{thm1}]
Let $\Xi$ be an OR equator of $R$. It follows from Lemma \ref{lem4} that $\Xi$ is an equator of $R^2$. Then the hyperbolic map $R^2$ is a $1$-mating in virtue of \cite[Theorem 4.2]{Mey1}.

\end{proof}

The following two results follow naturally from a combination of Theorem \ref{thm1} and \cite[Theorem 4.2]{Mey1}.

\begin{coro}\label{cor7}
If a hyperbolic postcritically finite rational map $R:\mathbb{P}^1(\mathbb{C})\rightarrow\mathbb{P}^1(\mathbb{C})$ admits an OR equator, then $R$ is an $n$-mating with its fold $F(R)=n\leq 2$. 
\end{coro}

\begin{coro}
If a hyperbolic postcritically finite rational map $R:\mathbb{P}^1(\mathbb{C})\rightarrow\mathbb{P}^1(\mathbb{C})$ admits an OR equator but does not admit any equator, then $R$ is a $2$-mating. 
\end{coro}

\subsection{Orientation-reversing matings versus primitive 2-matings}
It has ever been quite alluring to us that an inverse conclusion to Theorem \ref{thm1} would hold, that is,  every hyperbolic postcritically finite rational map arising as a 2-mating  would admit an OR equator. Surprisingly, this is wrong. See Section \ref{sec5} for 2-matings without any OR equator. Thus the hyperbolic 2-matings are classified into two disjoint families in fact: ones with OR equators and ones without OR equators.  

\begin{defn}
A  hyperbolic rational map arising as a 2-mating admitting an OR equator is called an orientation-reversing (OR) mating, while it is called a \emph{primitive 2-mating} if it does not admit any OR equator.
\end{defn}

The former notion has been suggested in \cite[Section 12]{Mey1}.  It would be interesting to try to give criteria on distinguishing the OR matings from the primitive 2-matings, possibly for certain families. Such criteria would also provide inverse conclusions for Lemma \ref{lem4} and Corollary \ref{cor1}, as well as extensions of Proposition \ref{pro1}. 
\begin{rem}
For non-hyperbolic $2$-matings, one can still use the classification to distinguish ones with and without OR equators, however, it does not make much sense as there are non-hyperbolic $1$-matings without any equator (\cite{Mey3, Mil2}). It seems better to consider classification of non-hyperbolic $2$-matings according to whether they admit OR pseudo-equators or not (\cite{Mey1}), at least for non-hyperbolic $2$-matings with their Julia sets being $\mathbb{P}^1(\mathbb{C})$.     
\end{rem}

For some postcritically finite rational map $R$, let $P_w(R)$ and $P_b(R)$ be some non-travail partition of $P(R)$. Let $\Xi$ be a Jordan curve circular to both $P_w(R)$ and $P_b(R)$. Be careful that the non-splitting property of $R^{-1}(\Xi)$  is not enough to guarantee $\Xi$ will evolve into an equator or OR equator of $R$. There is possibility that the two domains $U_w^0$ and $R^{-1}(U_w^0)=U_w^1(R)$ (alternatively, $U_w^0$ and $R^{-1}(U_b^0)=U_b^1$) may confine some postcritical points trickily which hinders the isotopy between $\Xi$ and $R(\Xi)$ \emph{rel.} $P(R)$, as one will see in the following sections.

\section{Matings of pairs of half polynomials}\label{sec4}

In this section we mean to describe the dynamics of hyperbolic postcritically finite rational maps arising as OR matings via matings of half polynomials. For hyperbolic rational maps arising as $1$-matings, the description of their dynamics is explicit through the unmating process. More precisely, if some hyperbolic rational map $R$ arises as a $1$-mating, let $\Xi$ be an equator and $\{H_t\}_{t\in I}$ be an isotopy from $\Xi$ to $R^{-1}(\Xi)$. Then the restrictive maps 
\begin{center}
$(H_1\circ R)|_{U_w^1}: U_w^1\rightarrow U_w^1$ and $(H_1\circ R)|_{U_b^1}: U_b^1\rightarrow U_b^1$
\end{center}  
are Thurston equivalent to two polynomials respectively, moreover, $R$ is topologically conjugate to the $1$-mating of the two polynomials. One can refer to \cite{Mey1, Wit}.

\subsection{Truncations of the deformed OR matings and the induced polynomials} For a degree-$d$ hyperbolic rational map $R$ arising as an OR mating with $d\geq 2$, the above results are still applicable to derive topological description of the Fatou and Julia set of $R$, as well as dynamics of the twice iteration $R^2$. In this case let $\Xi$ be an OR equator of $R$, let $\{H_t\}_{t\in I}$ be an isotopy from $\Xi$ to $R^{-1}(\Xi)$ (orientation-reversing of course). Now we consider the following two restrictive maps
\begin{equation}\label{eq12}
T_{wb}=(H_1\circ R)|_{U_w^1}: U_w^1\rightarrow U_b^1
\end{equation} 
and 
\begin{equation}\label{eq13}
T_{bw}=(H_1\circ R)|_{U_b^1}: U_b^1\rightarrow U_w^1.
\end{equation}
The two maps coincide with each other on $\partial U_w^1=\partial U_b^1=R^{-1}(\Xi)$, that is, 
\begin{equation}\label{eq16}
T_{wb}|_{\partial U_w^1}=T_{bw}|_{\partial U_b^1}=(H_1\circ R)|_{\partial U_w^1}=(H_1\circ R)|_{\partial U_b^1},
\end{equation}
which is a $d$-fold covering map. Be careful that $H_1$ satisfies (\ref{eq14}) instead of (\ref{eq4}) now.  Let 
\begin{equation}
T_w=T_{bw}\circ T_{wb}: U_w^1\rightarrow U_w^1
\end{equation}
and 
\begin{equation}
T_b=T_{wb}\circ T_{bw}: U_b^1\rightarrow U_b^1
\end{equation}
be the compositions. The two maps are both $2d$-branched covering maps. Following Meyer's arguments (possibly with mild modifications), we can obtain the following results.  

\begin{theorem}[Meyer-Wittner]\label{thm4}
For a hyperbolic postcritically finite $2$-mating $R$ with an OR equator $\Xi$, $T_w$ and $T_b$ are Thurston equivalent to two monic polynomials $P_w: \hat{\mathbb{C}}_w\rightarrow \hat{\mathbb{C}}_w$ and $P_b: \hat{\mathbb{C}}_b\rightarrow \hat{\mathbb{C}}_b$ respectively. Moreover, we have
\begin{center}
$R^2\stackrel{Thu}{\approx}P_w\ \mathrm{\rotatebox[origin=c]{90}{$\vDash$}}_\mathrm{F}\ P_b$ 
\end{center}
and
\begin{center}
$R^2\stackrel{top}{\approx} P_w$ \rotatebox[origin=c]{90}{$\vDash$} $P_b$.
\end{center}
\end{theorem}

Since $\mathcal{J}(R)=\mathcal{J}(R^2)$ ($\mathcal{F}(R)=\mathcal{F}(R^2)$), the mating provides topological description of the Julia (Fatou) set of $R$, as well as dynamics of $R^2$. Here we are particularly interested in the topological description of dynamics of $R$, instead of its twice iterate $R^2$. The mating of half polynomials is designed to achieve this goal. 

From now on we will work on the exploding Riemann sphere $\hat{\mathbb{C}}$ for convenience. However, all the results can be laid on the Riemann sphere in fact.  We will distinguish the white exploding Riemann sphere $\hat{\mathbb{C}}_w=\hat{\mathbb{C}}$ from the black one $\hat{\mathbb{C}}_b=\hat{\mathbb{C}}$ in the following. According to Theorem \ref{thm4}, since $T_w$ is Thurston equivalent to the monic polynomial $P_w$, there exist two homeomorphisms $f_0: \hat{\mathbb{C}}_w\rightarrow U_w^1$ and $f_1: \hat{\mathbb{C}}_w\rightarrow U_w^1$ isotopic to each other \emph{rel.} $P(P_w)\cup\infty S^1$, such that 
\begin{equation}\label{eq17}
T_{bw}\circ T_{wb}\circ f_0=T_w\circ f_0=f_1\circ P_w,
\end{equation} 
that is, the following diagram is commutative,
\begin{center}
\begin{tikzcd}
\hat{\mathbb{C}}_w \arrow[r, "f_0"] \arrow[d, "P_w"] & U_w^1 \arrow[d, "T_w"] \\
\hat{\mathbb{C}}_w \arrow[r, "f_1"]  & U_w^1.
\end{tikzcd}
\end{center}

Symmetrically, there exist two homeomorphisms $g_0: \hat{\mathbb{C}}_w\rightarrow U_w^1$ and $g_1: \hat{\mathbb{C}}_w\rightarrow U_w^1$ isotopic to each other \emph{rel.} $P(P_b)\cup\infty S^1$, such that 
\begin{equation}\label{eq18}
T_{wb}\circ T_{bw}\circ g_0=T_b\circ g_0=g_1\circ P_b,
\end{equation} 
that is, the following diagram is commutative,
\begin{center}
\begin{tikzcd}
\hat{\mathbb{C}}_b \arrow[r, "g_0"] \arrow[d, "P_b"] & U_b^1 \arrow[d, "T_b"] \\
\hat{\mathbb{C}}_b \arrow[r, "g_1"]  & U_b^1.
\end{tikzcd}
\end{center}

We require these homeomorphisms $f_0, f_1, g_0, g_1$ bear some consistence on the boundaries of their domains.

\begin{lemma}\label{lem8}
The homeomorphisms $f_0, f_1, g_0, g_1$ can be chosen such that
\begin{equation}\label{eq25}
f_0(\infty e^{2\pi it})=g_0(\infty e^{2\pi it})
\end{equation}
and
\begin{equation}\label{eq26}
f_1(\infty e^{2\pi it})=g_1(\infty e^{2\pi it})
\end{equation}
for any $t\in I$.
\end{lemma}
\begin{proof}
Recall from (\ref{eq15}) that, as degree-$2d$ polynomials on $\hat{\mathbb{C}}$, we have
\begin{center}
$P_w(\infty e^{2\pi it})=\infty e^{4\pi idt}$  
\end{center}
and 
\begin{center}
$P_b(\infty e^{2\pi it})=\infty e^{4\pi idt}$ 
\end{center}
on $\infty S^1$. Note that  $T_w|_{\partial U_w^1}: \partial U_w^1\rightarrow \partial U_w^1$ and $T_b|_{\partial U_b^1}: \partial U_b^1\rightarrow \partial U_b^1$ are both (orientation-preserving) $2d$-fold coverings without branching points. According to (\ref{eq16}), the two restrictive maps coincide with each other,
\begin{equation}\label{eq19}
T_w|_{\partial U_w^1}=T_b|_{\partial U_b^1}.
\end{equation}
Now we consider the restrictions of these maps on  $\partial \hat{\mathbb{C}}_w=\infty S^1$, $\partial \hat{\mathbb{C}}_b=\infty S^1$ and $\partial U_w^1=\partial U_b^1$ respectively. Recall that there exist homeomorphisms $f_0, f_1$ isotopic to each other \emph{rel.} $P(P_w)\cup\infty S^1$ satisfying (\ref{eq17}), so we have
\begin{equation}\label{eq27}
f_0|_{\infty S^1}=f_1|_{\infty S^1}.
\end{equation}
Symmetrically there exist homeomorphisms $g_0, g_1$ satisfying (\ref{eq18}) such that
\begin{equation}\label{eq28}
g_0|_{\infty S^1}=g_1|_{\infty S^1}
\end{equation}
since $g_0$ is isotopic to $g_1$ \emph{rel.} $P(P_b)\cup\infty S^1$. When they are restricted on the boundary $\infty S^1$, considering (\ref{eq27}) and (\ref{eq28}), the two equations (\ref{eq17}) and (\ref{eq18}) degenerate into
\begin{equation}\label{eq20}
T_w\circ f_0=f_0\circ z^{2d}
\end{equation}
and
\begin{equation}\label{eq21}
T_b\circ g_0=g_0\circ z^{2d}
\end{equation} 
respectively. Then (\ref{eq19}), (\ref{eq20}) and (\ref{eq21}) together induce
\begin{equation}\label{eq22}
z^{2d}=(f_0^{-1}\circ g_0)\circ z^{2d}\circ (f_0^{-1}\circ g_0)^{-1}
\end{equation}
on $\infty S^1$. Note that this type of Schr\"oder equation $(\varphi(z))^{2d}=\varphi(z^{2d})$ admits a unique solution up to multiplication by a $(2d-1)$-th root of unity. So
$$f_0^{-1}\circ g_0=e^{2\pi i\frac{j}{2d-1}}\mathds{1},$$  
or
\begin{equation}\label{eq23}
g_0(z)=f_0(e^{2\pi i\frac{j}{2d-1}}z)
\end{equation}
for some $0\leq j\leq 2d-1$ and any $z\in \infty S^1$. Now combining (\ref{eq17}) and (\ref{eq23}) together, we have the following commuting diagram,
\begin{center}
\begin{tikzcd}
\hat{\mathbb{C}}_w \arrow[r, "h_0"] \arrow[d, "P_{*w}"] & \hat{\mathbb{C}}_w \arrow[r, "f_0"] \arrow[d, "P_w"] & U_w^1 \arrow[d, "T_w"] \\
\hat{\mathbb{C}}_w \arrow[r, "h_1"] & \hat{\mathbb{C}}_w \arrow[r, "f_1"]  & U_w^1,
\end{tikzcd}
\end{center}
in which $P_{*w}$ is a monic polynomial M\"obius conjugate to $P_w$ and
\begin{equation}\label{eq24}
h_0(z)=h_1(z)=e^{2\pi i\frac{j}{2d-1}}z.
\end{equation}
Moreover, the maps 
\begin{center}
$f_{*0}=f_0\circ h_0$ 
\end{center}
and 
\begin{center}
$f_{*1}=f_1\circ h_1$ 
\end{center}
are isotopic to each other \emph{rel.} $P(P_{*w})\cup\infty S^1$. Note that 
\begin{center}
$f_{*0}|_{\infty S^1}=g_0|_{\infty S^1}$ 
\end{center}
and 
\begin{center}
$f_{*1}|_{\infty S^1}=g_1|_{\infty S^1}$ 
\end{center}
in virtue of (\ref{eq27}), (\ref{eq28}) and (\ref{eq23}). It is easy to see that $f_0, f_1$ and $P_w$ can be substituted by the new homeomorphisms $f_{*0}, f_{*1}$ and polynomial $P_{*w}$ such that (\ref{eq17}), (\ref{eq18}) still hold, while all the homeomorphisms $f_{*0}, f_{*1}, g_0, g_1$ coincide with each other on the boundary.
\end{proof}

\subsection{Half polynomials and their matings}

The half polynomials are maps between $\hat{\mathbb{C}}_w$ and $\hat{\mathbb{C}}_b$, which always appear in pairs.

\begin{defn}\label{def1}
Let $\mathcal{P}_{wb}: \hat{\mathbb{C}}_w\rightarrow \hat{\mathbb{C}}_b$ and $\mathcal{P}_{bw}: \hat{\mathbb{C}}_b\rightarrow \hat{\mathbb{C}}_w$ be a pair of maps satisfying
\begin{enumerate}[(1).]
\item\label{ite1} $\mathcal{P}_{wb}(\infty S^1)=\mathcal{P}_{bw}(\infty S^1)=\infty S^1=\mathcal{P}_{wb}^{-1}(\infty S^1)=\mathcal{P}_{bw}^{-1}(\infty S^1)$.

\item\label{ite2}  $\mathcal{P}_{wb}|_{\infty S^1}=\mathcal{P}_{bw}|_{\infty S^1}$. 
\end{enumerate}
They are called a pair of (monic) half polynomials, if there exists two (monic) polynomials $P_w: \hat{\mathbb{C}}_w\rightarrow \hat{\mathbb{C}}_w, P_b: \hat{\mathbb{C}}_b\rightarrow \hat{\mathbb{C}}_b$ and two homeomorphisms $f: \hat{\mathbb{C}}_w\rightarrow \hat{\mathbb{C}}_w, g: \hat{\mathbb{C}}_b\rightarrow \hat{\mathbb{C}}_b$ satisfying $f|_{\infty S^1}=\mathds{1}_{\infty S^1}=g|_{\infty S^1}$, such that 
\begin{equation}\label{eq29}
\mathcal{P}_{wb}\circ f\circ \mathcal{P}_{bw}=P_b
\end{equation}
and
\begin{equation}\label{eq30}
\mathcal{P}_{bw}\circ g\circ \mathcal{P}_{wb}=P_w.
\end{equation}
\end{defn}

Be careful that $\hat{\mathbb{C}}_w$ and $\hat{\mathbb{C}}_b$ are viewed as two distinct copies of $\hat{\mathbb{C}}$ now, so the half polynomials are in fact not endomorphisms on $\hat{\mathbb{C}}$. The maps $f$ and $g$ are called \emph{bonding maps} of the half polynomials $\mathcal{P}_{wb}, \mathcal{P}_{bw}$.  This definition fits into our theory in this work, while we suspect it may be of additional interest in some special cases. For example, in case $f=g=\mathds{1}$ on $\hat{\mathbb{C}}_w$ and $\hat{\mathbb{C}}_b$,  (\ref{eq29}) degenerates to
\begin{center}
$\mathcal{P}_{wb}\circ \mathcal{P}_{bw}=P_b$
\end{center}
while (\ref{eq30}) degenerates to
\begin{center}
$\mathcal{P}_{bw}\circ \mathcal{P}_{wb}=P_w.$
\end{center}

Any pair of (monic) polynomials can of course be viewed as a pair of (monic) half polynomials. The half polynomials can be employed to describe the dynamics of $2$-matings directly. To do this we need the following definition. For $\hat{\mathbb{C}}_w$ and $\hat{\mathbb{C}}_b$ being two copies of $\hat{\mathbb{C}}$, let $\sim_+$ be the relationship $\infty e^{2\pi it}\sim_+\infty e^{2\pi it}$ for points in $\hat{\mathbb{C}}_w$ and $\hat{\mathbb{C}}_b$ respectively, $t\in I$. Let $\hat{\mathbb{C}}_w\sqcup\hat{\mathbb{C}}_b$ be the disjoint union.
\begin{defn}\label{def5}
For a pair of half polynomials $\mathcal{P}_{wb}: \hat{\mathbb{C}}_w\rightarrow \hat{\mathbb{C}}_b$ and $\mathcal{P}_{bw}: \hat{\mathbb{C}}_b\rightarrow \hat{\mathbb{C}}_w$, their mating
\begin{center}
 $\mathcal{P}_{wb}\amalg \mathcal{P}_{bw}: \hat{\mathbb{C}}_w\sqcup\hat{\mathbb{C}}_b/\sim_+\rightarrow\hat{\mathbb{C}}_w\sqcup\hat{\mathbb{C}}_b/\sim_+$ 
\end{center}
is defined to be
$$\mathcal{P}_{wb}\amalg \mathcal{P}_{bw}(z)=\left\{
\begin{array}{lll}
\mathcal{P}_{wb}(z) & for & z\in \hat{\mathbb{C}}_w,\\
\mathcal{P}_{bw}(z) & for & z\in \hat{\mathbb{C}}_b.
\end{array}
\right.$$
\end{defn}

The mating is well-defined according to Definition \ref{def1} \ref{ite2}  We use a different symbol so one can easily distinguish it from the mating of polynomials. Now we are well-prepared to justify Theorem \ref{thm5}.

\begin{proof}[Proof of Theorem \ref{thm5}]
For the hyperbolic postcritically finite rational map $R$ admitting some OR equator $\Xi$, let  $P_w: \hat{\mathbb{C}}_w\rightarrow \hat{\mathbb{C}}_w$ and $P_b: \hat{\mathbb{C}}_b\rightarrow \hat{\mathbb{C}}_b$ be polynomials Thurston equivalent to $T_w$ and $T_b$ respectively. In virtue of Lemma \ref{lem8}, assume the homeomorphisms $f_0, f_1, g_0, g_1$ satisfy (\ref{eq25}) and (\ref{eq26}). We rewrite (\ref{eq17}) as
\begin{equation}\label{eq31}
f_1^{-1}\circ T_{bw}\circ g_0\circ g_0^{-1}\circ g_1\circ g_1^{-1}\circ T_{wb}\circ f_0=P_w,
\end{equation} 
rewrite (\ref{eq18}) as
\begin{equation}\label{eq32}
g_1^{-1}\circ T_{wb}\circ f_0\circ f_0^{-1}\circ f_1\circ f_1^{-1}\circ T_{bw}\circ g_0=P_b.
\end{equation} 
Now let 
\begin{equation}\label{eq34}
\mathcal{P}_{bw}=f_1^{-1}\circ T_{bw}\circ g_0
\end{equation}
and 
\begin{equation}\label{eq35}
\mathcal{P}_{wb}=g_1^{-1}\circ T_{wb}\circ f_0,
\end{equation}
let
\begin{center}
$g=g_0^{-1}\circ g_1$ and $f=f_0^{-1}\circ f_1$. 
\end{center}
Then (\ref{eq31}) and (\ref{eq32}) induce (\ref{eq30}) and (\ref{eq29}) respectively. It is easy to check that $\mathcal{P}_{bw}$ and $\mathcal{P}_{wb}$ satisfy Definition \ref{def1} \ref{ite1}  Considering (\ref{eq16}), (\ref{eq25}) and (\ref{eq26}), $\mathcal{P}_{bw}$ and $\mathcal{P}_{wb}$ satisfy Definition \ref{def1} \ref{ite2}  Note that since $f_0|_{\infty S^1}=f_1|_{\infty S^1}$ and $g_0|_{\infty S^1}=g_1|_{\infty S^1}$, we have
\begin{center}
$f|_{\infty S^1}=\mathds{1}=g|_{\infty S^1}$.
\end{center}
These guarantee that $\{\mathcal{P}_{bw}, \mathcal{P}_{wb}\}$ is a pair of half polynomials. In the following we show that (\ref{eq33}) holds. To do this, first note that
\begin{center}
$H_1\circ R|_{U_w^1}= T_{wb}=g_1\circ \mathcal{P}_{wb}\circ f_0^{-1}$
\end{center}
considering (\ref{eq12}) and (\ref{eq35}). So
\begin{equation}\label{eq59}
R|_{U_w^1}=H_1^{-1}\circ g_1\circ \mathcal{P}_{wb}\circ f_0^{-1}: U_w^1\rightarrow U_w^0.
\end{equation}
Following similar lines we can deduce
\begin{equation}\label{eq60}
R|_{U_b^1}=H_1^{-1}\circ f_1\circ \mathcal{P}_{bw}\circ g_0^{-1}: U_b^1\rightarrow U_b^0.
\end{equation}
Now let 
$$h_0(z)=\left\{
\begin{array}{lll}
f_0(z) & for & z\in \hat{\mathbb{C}}_w,\\
g_0(z) & for & z\in \hat{\mathbb{C}}_b,
\end{array}
\right.$$
and
$$h_1(z)=\left\{
\begin{array}{lll}
H_1^{-1}\circ f_1(z) & for & z\in \hat{\mathbb{C}}_w,\\
H_1^{-1}\circ g_1(z) & for & z\in \hat{\mathbb{C}}_b.
\end{array}
\right.$$
The two homeomorphisms are well-defined in virtue of (\ref{eq25}) and (\ref{eq26}). Moreover, $h_0$ is isotopic to $h_1$ \emph{rel.} $P(\mathcal{P}_{wb})\cup P(\mathcal{P}_{bw})$. Now (\ref{eq59}) and (\ref{eq60}) together give
$$R\circ h_0=h_1\circ (\mathcal{P}_{wb}\amalg \mathcal{P}_{bw}),$$
which completes the proof of (\ref{eq33}). 
\end{proof}

For hermaphroditic matings, Theorem \ref{thm5} also provides another interpretation of dynamics of them besides unmating them into pairs of polynomials of the same degrees.

\section{$4$-matings, hermaphroditic matings and OR matings in the bicritical family}\label{sec6}

A  rational map is called \emph{bicritical} if it admits exactly two critical points. Let $\EuScript{M}$ be the collection of  all the bicritical rational maps of degree $k\geq 2$ up to holomorphic conjugacies (we customarily use $d$ to denote the degree of a rational map, however, the symbol is expropriated for another parameter in this section).  Milnor (\cite{Mil3}) initiated the study on the geometric properties of this moduli space according to dynamics of maps in the family. He classified $\EuScript{M}$ into three basic families: the \emph{connectedness locus}, the \emph{hyperbolic shift locus}, and the \emph{parabolic shift locus}, 
$$\EuScript{M}=\EuScript{C}\cup \EuScript{S}_{hyp}\cup \EuScript{S}_{par},$$
according to connectivity of the Julia sets and expanding property of maps on their Julia sets.   He gave geometric description of all the three partitions in the moduli space $\EuScript{M}$, as well as the smooth curves $Per_1(\lambda)$ (the equivalent classes of bicritical maps with a fixed point of multiplier $\lambda$). One can find maps with various interesting properties in the bicritical family, however, here we are particularly interested in the unmating property of bicritical maps in the connectedness locus. 

\subsection{A family of postcritically finite maps in the connectedness locus}
By putting one critical point at $0$ and the other one at infinity, a bicritical rational map  of degree $k$ is always of the form
$$R_{a,b,c,d}=\cfrac{az^k+b}{cz^k+d}$$
for some $(a,b,c,d)\in \mathbb{C}^4$  (\cite[Lemma 1.1]{Mil3}). Without loss of generality one can require $(a,b,c,d)\in \{(a,b,c,d): ad-bc=1\}$ since the map depends only on the ratio of these parameters. For $k\geq 2$, let 

$$\Omega_1^k=\{(a,b,c,d) \in \mathbb{C}^4: ad-bc=1, cb^k+d^{k+1}=0\},$$

$$\Omega_2^k=\{(a,b,c,d) \in \mathbb{C}^4: ad-bc=1, a^{k+1}+bc^k=0\},$$

$$\Omega_3^k=\{(a,b,c,d) \in \mathbb{C}^4: ad-bc=1, ab^{k-1}+d^k=0\},$$
and
$$\Omega_4^k=\{(a,b,c,d) \in \mathbb{C}^4: ad-bc=1, a^k+dc^{k-1}=0\}.$$

$\Omega_j^k$ is a two-dimension complex manifold for any $j\in\{1,2,3,4\}$ and $k\geq 2$.   Now let 

$$\Omega_+^k=(\Omega_1^k\setminus \Omega_3^k)\cap(\Omega_2^k\setminus\Omega_4^k)=(\Omega_1^k\cap\Omega_2^k)\setminus (\Omega_3^k\cup \Omega_4^k)$$  
and
$$\Omega_-^k=(\Omega_3^k\setminus \Omega_1^k)\cap(\Omega_4^k\setminus\Omega_2^k)=(\Omega_3^k\cap\Omega_4^k)\setminus (\Omega_1^k\cup \Omega_2^k)$$
for $k\geq 2$, in which  $\Omega_1^k\cap\Omega_2^k$ and $\Omega_3^k\cap\Omega_4^k$ are all Riemann surfaces. Let $\Omega^k=\Omega_+^k\cup \Omega_-^k$ for $k\geq 2$ and $\Omega=\cup_{k=2}^\infty\Omega^k$.  The following result demonstrates that maps with parameters in $\Omega$ are all in the connectedness locus in fact.
\begin{pro}\label{pro6}
Let $k\geq 2$ be an integer. For any $(a,b,c,d) \in \Omega^k$, we have
\begin{center}
$(R_{a,b,c,d})\in \EuScript{C}$.
\end{center}
\end{pro}
\begin{proof}
For any sequel of parameters $(a,b,c,d) \in \Omega_1^k\setminus\Omega_3^k$, we can guarantee that
\begin{equation}\label{eq38}
R_{a,b,c,d}(\frac{b}{d})=\infty.
\end{equation}
For any sequel of parameters $(a,b,c,d) \in \Omega_2^k\setminus\Omega_4^k$, we can guarantee that
\begin{equation}\label{eq39}
R_{a,b,c,d}(\frac{a}{c})=0.
\end{equation}
Then for any  $(a,b,c,d) \in \Omega_+^k$, combining  (\ref{eq38}) and (\ref{eq39}) together, we see the orbits of the two critical points $0$ and $\infty$ as following,
\begin{equation}\label{eq40}
\begin{tikzcd}
0 \arrow[r, "R_{a,b,c,d}"] & \cfrac{b}{d} \arrow[r, "R_{a,b,c,d}"]  & \infty \arrow[r, "R_{a,b,c,d}"] & \cfrac{a}{c} \arrow[lll, bend left, "R_{a,b,c,d}"].
\end{tikzcd}
\end{equation}
Note that we always have  $\cfrac{b}{d}\neq \cfrac{a}{c}$ for any $(a,b,c,d) \in \Omega_+^k$. 

Now for any sequel of parameters $(a,b,c,d) \in \Omega_3^k\setminus\Omega_1^k$, we have
\begin{equation}\label{eq41}
R_{a,b,c,d}(\frac{b}{d})=0.
\end{equation}
For any sequel of parameters $(a,b,c,d) \in \Omega_4^k\setminus\Omega_2^k$, we have
\begin{equation}\label{eq42}
R_{a,b,c,d}(\frac{a}{c})=\infty.
\end{equation}
Then for any  $(a,b,c,d) \in \Omega_-^k$, combining  (\ref{eq41}) and (\ref{eq42}) together, we see the orbits of the two critical points $0$ and $\infty$ are typically two disjoint cycles,
\begin{equation}\label{eq43}
\begin{array}{c}
\begin{tikzcd}
0 \arrow[r, bend left, "R_{a,b,c,d}"] & \cfrac{b}{d} \arrow[l, bend left, "R_{a,b,c,d}"],
\end{tikzcd}\\
\begin{tikzcd}
\infty \arrow[r, bend left, "R_{a,b,c,d}"] & \cfrac{a}{c} \arrow[l, bend left, "R_{a,b,c,d}"],
\end{tikzcd}
\end{array}
\end{equation}
since we always have  $\cfrac{b}{d}\neq \cfrac{a}{c}$ for any $(a,b,c,d) \in \Omega_-^k$, except the case $k=2$. The two cycles degenerate into one in case $k=2$, see Lemma \ref{lem15}. The above arguments show that for any $(a,b,c,d) \in \Omega^k$, $R_{a,b,c,d}$ is always postcritically finite, which forces $\mathcal{J}(R_{a,b,c,d})$ to be connected.

\end{proof}

In fact both $\Omega_+^k$ and $\Omega_-^k$ can be simplified, see Lemma \ref{lem16}. One can easily see the following result from the proof of Proposition \ref{pro6}. 
\begin{coro}
For any postcritically finite bicritical map $R_{a,b,c,d}$ with $(a,b,c,d) \in \Omega^k$ and $k\geq 2$, we always have 
$$P(R_{a,b,c,d})=\{0, \cfrac{b}{d}, \infty, \cfrac{a}{c}\}.$$
\end{coro}
Be careful that $\#P(R_{a,b,c,d})$ is not always $4$ for any $R_{a,b,c,d}$ with $(a,b,c,d) \in \Omega^k$ and $k\geq 2$. In case of $k=2$ we have $\#P(R_{a,b,c,d})=2$ for $(a,b,c,d) \in \Omega_-^2$, see  Corollary \ref{cor2}. This property leads to some peculiar unmating behaviour for these postcritically finite maps of degree $2$, compared with their cousins of higher degree in the bicritical family.  
 
There are plenty maps of interesting unmating properties in $\EuScript{M}$.  In our concerns, we find $4$-matings, hermaphroditic matings and orientation-reversing matings in the connectedness locus  $\EuScript{C}$.

\subsection{Unmatability of some sample maps of certain degrees with parameters in $\Omega$}

Our final goal in this section would be to prove Theorem \ref{thm8}. However, instead of presenting the shrivelling proof directly,  we start from exploring the  unmatability of some special maps. There are two reasons for us to cope with the problem in this style. The first one is that through explorations on the unmatability of these special maps the abstract concepts and techniques are condensed into concrete ones, which will make the proof more ample and easily understandable. The readers can probably grasp the core ideas through these special cases before the proof of Theorem \ref{thm8}.  The second one is that some maps take some distinguished features individually, not exhibited in the proof of Theorem \ref{thm8}, which may inspire the readers in some other aspects.  

Let $\beta_{3,1}=e^{2\pi i/3}$ be one of the cubic roots of $1$. Let 
$$\omega_{+,2}=(\cfrac{1}{\sqrt{\beta_{3,1}-1}}, -\cfrac{1}{\sqrt{\beta_{3,1}-1}}, -\cfrac{1}{\sqrt{\beta_{3,1}-1}}, \cfrac{\beta_{3,1}}{\sqrt{\beta_{3,1}-1}}).$$
One can check that $\omega_{+,2}\in\Omega_+^2$. There are in fact two conjugate choices for $\sqrt{\beta_{3,1}-1}$, while one can take either of them, it does not affect the results below. The induced bicritical map is
$$R_{\omega_{+,2}}=\cfrac{z^2-1}{-z^2+\beta_{3,1}}.$$
The critical points of the map $R_{\omega_{+,2}}$ go as following,
\begin{center}
\begin{tikzcd}
0 \arrow[r, "R_{\omega_{+,2}}"] & -\beta_{3,1}^{-1} \arrow[r, "R_{\omega_{+,2}}"]  & \infty \arrow[r, "R_{\omega_{+,2}}"] & -1 \arrow[lll, bend left, "R_{\omega_{+,2}}"].
\end{tikzcd}
\end{center}

\begin{figure}[h]
\centering
\includegraphics[scale=1]{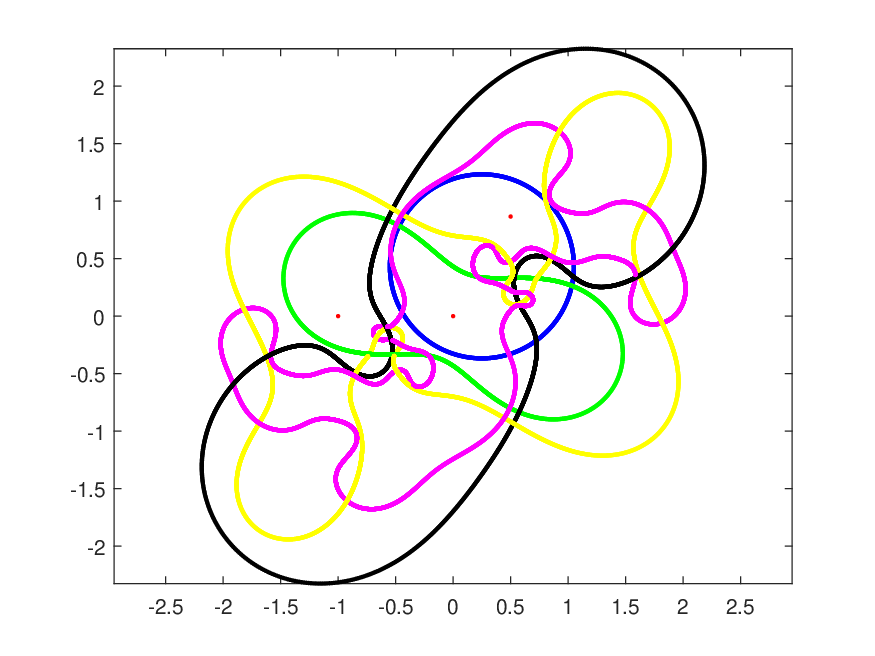}
\caption{$\Xi_{\omega_{+,2},V}$ and $\{R_{\omega_{+,2}}^{-j}(\Xi_{\omega_{+,2},V})\}_{j=1}^4$}
\label{fig9}
\end{figure}

\begin{figure}[h]
\centering
\includegraphics[scale=1]{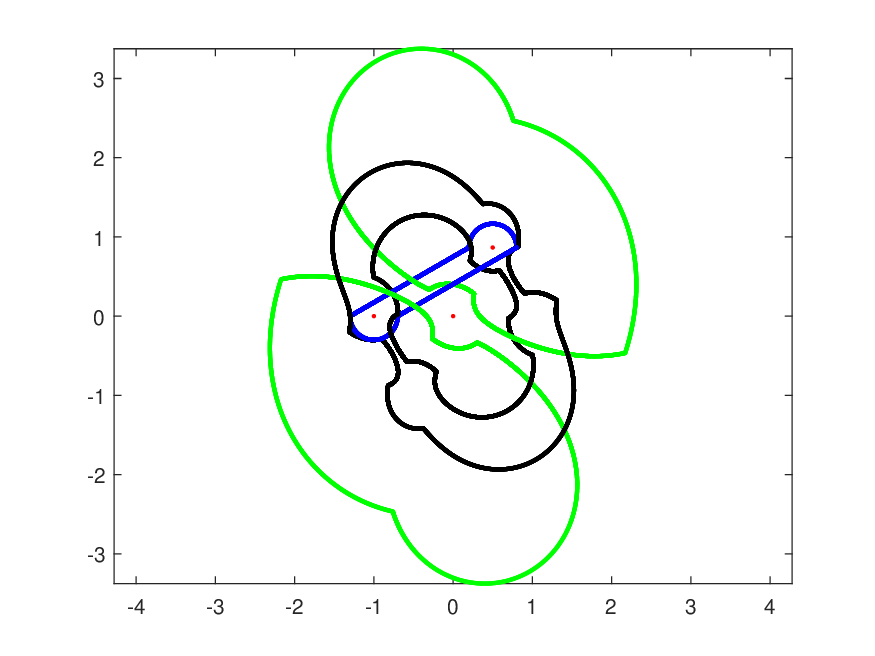}
\caption{$\Xi_{\omega_{+,2},VI}$, $R_{\omega_{+,2}}^{-1}(\Xi_{\omega_{+,2},VI})$ and $R_{\omega_{+,2}}^{-2}(\Xi_{\omega_{+,2},VI}$)}
\label{fig10}
\end{figure}

\begin{figure}[h]
\centering
\includegraphics[scale=1]{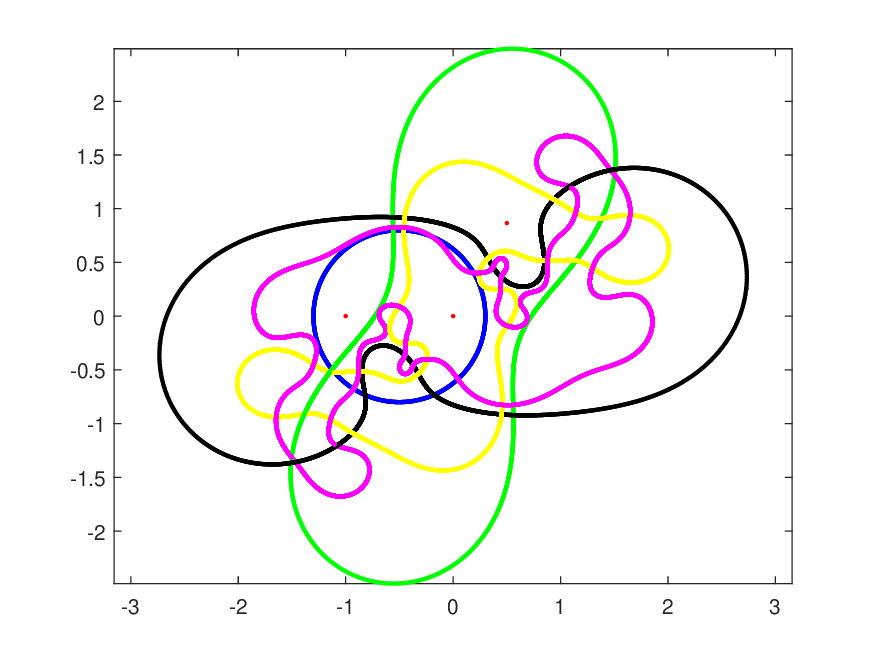}
\caption{$\Xi_{\omega_{+,2},VII}$ and $\{R_{\omega_{+,2}}^{-j}(\Xi_{\omega_{+,2},VII})\}_{j=1}^4$}
\label{fig11}
\end{figure}

Considering the unmatability of the map $R_{\omega_{+,2}}$, we have the following result.
\begin{pro}\label{pro4}
The map $R_{\omega_{+,2}}$ is a $4$-mating. 
\end{pro}
\begin{proof}
First note that for any $j\in\mathbb{N}_+$, $P(R_{\omega_{+,2}}^j)=P(R_{\omega_{+,2}})=\{0, -\beta_{3,1}^{-1}, \infty, -1\}$. Moreover, the orbits of the postcritical points of $R_{\omega_{+,2}}^j$ follow exactly one of the following four styles,

\begin{center}
$\begin{array}{ccc}
\begin{tikzcd}
0 \arrow[r, "R_{\omega_{+,2}}"] & -\beta_{3,1}^{-1} \arrow[r, "R_{\omega_{+,2}}"]  & \infty \arrow[r, "R_{\omega_{+,2}}"] & -1 \arrow[lll, bend left, "R_{\omega_{+,2}}"],
\end{tikzcd}\\
\begin{tikzcd}
\ \ \ 0 \arrow[r, bend left, "R_{\omega_{+,2}}^2"] & \infty \arrow[l, bend left, "R_{\omega_{+,2}}^2"],
\end{tikzcd}
\begin{tikzcd}
-\beta_{3,1}^{-1} \arrow[r, bend left, "R_{\omega_{+,2}}^2"]  & \ \ \ -1, \ \ \  \arrow[l, bend left, "R_{\omega_{+,2}}^2"]
\end{tikzcd}\\
\begin{tikzcd}
0 \arrow[r, "R_{\omega_{+,2}}^3"] & -1 \arrow[r, "R_{\omega_{+,2}}^3"]  & \infty \arrow[r, "R_{\omega_{+,2}}^3"] &  -\beta_{3,1}^{-1} \arrow[lll, bend left, "R_{\omega_{+,2}}^3"],
\end{tikzcd}
\end{array}$
\end{center}
\begin{center}
\begin{tikzpicture}
    \node (a) at (0,0) {0};
    \draw[->] (a) to[in=45, out=-45, looseness=10] node{\ \ \ \ \ \ \  $R_{\omega_{+,2}}^4,$} (a);
    \node (b) at (3,0) {$-\beta_{3,1}^{-1}$};
    \draw[->] (b) to[in=45, out=-45, looseness=5] node{\ \ \ \ \ \ \  $R_{\omega_{+,2}}^4,$} (b);
     \node (c) at (6,0) {$\infty$};
    \draw[->] (c) to[in=45, out=-45, looseness=10] node{\ \ \ \ \ \ \  $R_{\omega_{+,2}}^4,$} (c);
    \node (d) at (9,0) {-1};
    \draw[->] (d) to[in=45, out=-45, looseness=10] node{\ \ \ \ \ \ \  $R_{\omega_{+,2}}^4$} (d);
\end{tikzpicture}
\end{center}
for any $j\in\mathbb{N}_+$. To look for an equator or OR equator of $P(R_{\omega_{+,2}}^j)$, we consider partitions of the postcritical set into two subsets satisfying the immune or swapping property.  We list all the non-travail partitions of $P(R_{\omega_{+,2}}^j)$ into two subsets in Table \ref{tab4}.
\begin{table}[ht]
\caption{Partitions of $P(R_{\omega_{+,2}}^j)$} 
\centering 
\begin{tabular}{c c c} 
\hline 
Type &  $P_w(R_{\omega_{+,2}}^j)$ & $P_b(R_{\omega_{+,2}}^j)$ \\
\hline
I     &   $0$    &    $-\beta_{3,1}^{-1}, \infty, -1$\\

II   &    $-\beta_{3,1}^{-1}$    &    $0, \infty, -1$\\

III     &   $\infty$    &    $0, -\beta_{3,1}^{-1}, -1$\\

IV     &   $-1$    &    $0, -\beta_{3,1}^{-1}, \infty$\\

V     &   $0, -\beta_{3,1}^{-1}$    &    $\infty, -1$\\

VI   &    $0, \infty$    &    $-\beta_{3,1}^{-1},  -1$\\

VII     &   $0, -1$    &    $-\beta_{3,1}^{-1}, \infty$\\
\hline 
\end{tabular}
\label{tab4} 
\end{table} 
One can see that all the listed partitions satisfy the immune property under $R_{\omega_{+,2}}^4$. The Type V and VII partition in Table \ref{tab4} satisfy the swapping  property under $R_{\omega_{+,2}}^2$. The Type VI partition in Table \ref{tab4} satisfies the immune  property under $R_{\omega_{+,2}}^2$, while it satisfies the swapping property under $R_{\omega_{+,2}}$ and $R_{\omega_{+,2}}^3$.  We deal with the unmatability according to the listed types of partitions.
\begin{itemize}
\item  Type I-IV. In these cases there are no equators for $R_{\omega_{+,2}}^j$ for any $j\in\mathbb{N}_+$. Let $\Xi$ be a Jordan curve which partitions the postcritical set into one of the four types mentioned before. One can easily show that $R_{\omega_{+,2}}^{-4}(\Xi)$ splits into more than two disjoint connected components, so $\Xi$ can not serve as an equator (or OR equator) for $R_{\omega_{+,2}}^j$ for any $j\in\mathbb{N}_+$.    

\item  Type V. In this case there are OR equators for $R_{\omega_{+,2}}^2$ and equators for $R_{\omega_{+,2}}^4$ . Let $\Xi$ be any Jordan curve which is circular to both 
\begin{center}
$P_w(R_{\omega_{+,2}}^2)=\{0, -\beta_{3,1}^{-1}\}$ and $P_b(R_{\omega_{+,2}}^2)=\{\infty, -1\}$.
\end{center}
One can show that $R_{\omega_{+,2}}^{-2}(\Xi)$ is always connected  and still circular to both $P_w(R_{\omega_{+,2}}^2)$ and $P_b(R_{\omega_{+,2}}^2)$, so $\Xi$ is an OR equator for $R_{\omega_{+,2}}^2$ considering the swapping property of $P_w(R_{\omega_{+,2}}^2)$ and $P_b(R_{\omega_{+,2}}^2)$ under $R_{\omega_{+,2}}^2$. It is also an equator for $R_{\omega_{+,2}}^4$ in virtue of Lemma \ref{lem4}. One can refer to Figure \ref{fig9} for a concrete curve $\Xi_{\omega_{+,2}, V}$.   

\item Type VI. In this case there are still no equators (or OR equators) for $R_{\omega_{+,2}}^2$ (or $R_{\omega_{+,2}}, R_{\omega_{+,2}}^3$). In fact the pre-image $R_{\omega_{+,2}}^{-1}(\Xi)$ of any Jordan curve $\Xi$ circular to $\{0, \infty\}$  and $\{-\beta_{3,1}^{-1},  -1\}$ splits into two disjoint connected components, which forces $R_{\omega_{+,2}}^{-j}(\Xi)$ to split into at least two disjoint connected components for any $j\geq 1$. One can refer to Figure \ref{fig10} for the splitting process of some concrete Jordan curve $\Xi_{\omega_{+,2}, VI}$ circular to $\{0, \infty\}$  and $\{-\beta_{3,1}^{-1},  -1\}$.

\item  Type VII. In this case there are OR equators for $R_{\omega_{+,2}}^2$ and equators for $R_{\omega_{+,2}}^4$. The case is symmetric to the Type V case. One can refer to Figure \ref{fig11} for a concrete OR equator $\Xi_{\omega_{+,2}, VII}$ for $R_{\omega_{+,2}}^2$ (an equator for $R_{\omega_{+,2}}^4$).
\end{itemize} 
After exhausting all the types of partitions, we summarize that $R_{\omega_{+,2}}$ is a $4$-mating (which may be called an orientation-reversing $4$-mating).
\end{proof}

Let $r=0.3$. The circle 
$$\Xi_{\omega_{+,2}, V}=\{-(2\beta_{3,1})^{-1}+(0.5+r)e^{2\pi ti}\}_{0\leq t< 1}$$
inducing Type V partition (in Table \ref{tab4}) of the postcritical set in Figure \ref{fig9} is an (OR) equator for ($R_{\omega_{+,2}}^2$) $R_{\omega_{+,2}}^4$.   

Let $r=0.3$, $\theta=angle(1-(\beta_{3,1})^{-1})$ and $L=|1-(\beta_{3,1})^{-1}|$. We illustrate the readers how the pre-images of the Jordan curve
\begin{center}
$
\begin{array}{ll}
\Xi_{\omega_{+,2},VI}=& \{-(\beta_{3,1})^{-1}+re^{2\pi it}\}_{0\leq t\leq 0.5}\cup\{\tan\theta(t+r+1)\}_{-1-r\leq t\leq -1-r+L\cos\theta}\vspace{3mm}\\
& \cup\{-1+re^{2\pi it}\}_{0.5\leq t\leq 1}\cup\{\tan\theta(t-r+1)\}_{r-1\leq t\leq r-1+L\cos\theta}
\end{array}
$
\end{center}
inducing Type VI partition (in Table \ref{tab4}) of the postcritical set split in Figure \ref{fig10}.

Let $r=0.3$. The circle 
$$\Xi_{\omega_{+,2}, VII}=\{-0.5+(0.5+r)e^{2\pi ti}\}_{0\leq t< 1}$$
inducing Type VII partition (in Table \ref{tab4}) of the postcritical set in Figure \ref{fig11} is an (OR) equator for ($R_{\omega_{+,2}}^2$) $R_{\omega_{+,2}}^4$.   

It seems from the proof of Proposition \ref{pro4} that there are two different ways in unmating the map $R_{\omega_{+,2}}^4$, however, essentially, it is not a shared $1$-mating. Suppose $\Xi_{V}$ is an equator of $R_{\omega_{+,2}}^4$ induced from  Type V partition (in Table \ref{tab4}) of the postcritical set such that $P_w(R_{\omega_{+,2}}^4)=\{0, -\beta_{3,1}^{-1}\}$ and $P_b(R_{\omega_{+,2}}^4)=\{\infty, -1\}$. Then we get a pair of degree-$16$ polynomials $P_{w,V}$ and $P_{b,V}$, such that
\begin{center}
$P_{w,V}$ \rotatebox[origin=c]{90}{$\vDash$} $P_{b,V}\stackrel{top}{\approx}R_{\omega_{+,2}}^4$.
\end{center} 
Now if $\Xi_{VII}$ is an equator of $R_{\omega_{+,2}}^4$ induced from the Type VII partition (in Table \ref{tab4}) of the postcritical set such that $P_w(R_{\omega_{+,2}}^2))=\{0, -1\}$ and $P_b(R_{\omega_{+,2}}^2)=\{\infty, -\beta_{3,1}^{-1}\}$. Then we get a pair of degree-$16$ polynomials $P_{w,VII}$ and $P_{b,VII}$, such that
\begin{center}
$P_{w,VII}$ \rotatebox[origin=c]{90}{$\vDash$} $P_{b,VII}\stackrel{top}{\approx}R_{\omega_{+,2}}^4$.
\end{center} 
We have 
\begin{center}
$P_{w,V}=P_{b,VII}$ and $P_{b,V}=P_{w,VII}$
\end{center}
up to M\"obius conjugations  by comparing the Hubbard trees of the two pairs of polynomials in fact.

We then consider the map 
$$R_{\omega_{+,3}}=\cfrac{z^3-1}{z^3+1}$$
for $\omega_{+,3}=(\frac{1}{\sqrt{2}}, -\frac{1}{\sqrt{2}}, \frac{1}{\sqrt{2}}, \frac{1}{\sqrt{2}})\in \Omega_+^3$.
The critical points of the map $R_{\omega_{+,3}}$ go as
\begin{center}
\begin{tikzcd}
0 \arrow[r, "R_{\omega_{+,3}}"] & -1 \arrow[r, "R_{\omega_{+,3}}"]  & \infty \arrow[r, "R_{\omega_{+,3}}"] & 1 \arrow[lll, bend left, "R_{\omega_{+,3}}"].
\end{tikzcd}
\end{center}

It turns out that the map $R_{\omega_{+,3}}$ inherits the unmatability of the map $R_{\omega_{+,2}}$.
\begin{pro}\label{pro3}
The map $R_{\omega_{+,3}}$ is a $4$-mating. 
\end{pro}

\begin{proof}
The proof follows almost in line with the proof of Proposition \ref{pro4}, details are left to the readers.
\end{proof}

We again illustrate readers the behaviours of some concrete Jordan curves inducing different partitions of the postcritical set of $R_{\omega_{+,3}}^j$ by some pictures, in case one may discover further phenomenons through them.

\begin{figure}[h]
\centering
\includegraphics[scale=1]{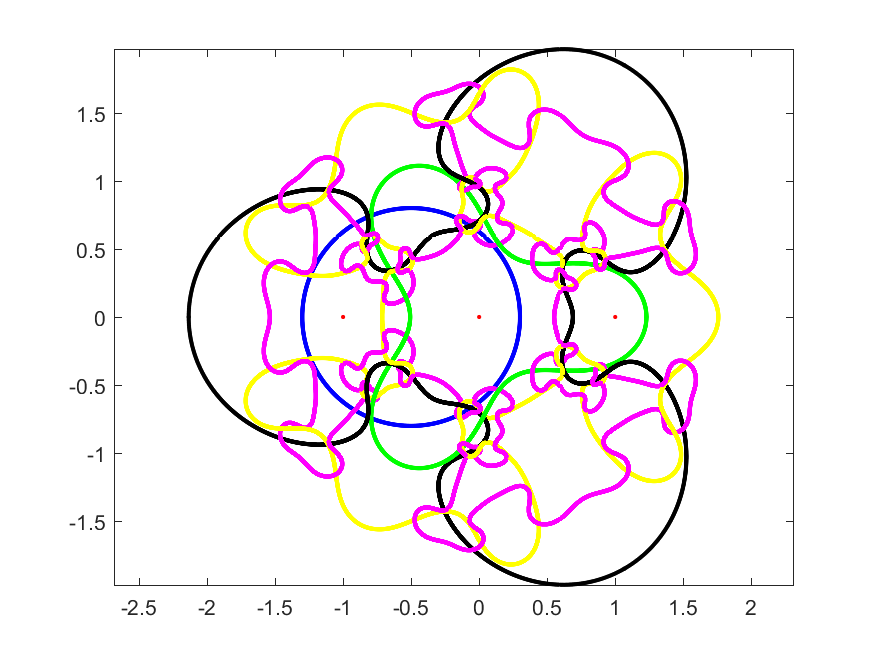}
\caption{$\Xi_{\omega_{+,3},V}$ and $\{R_{\omega_{+,3}}^{-j}(\Xi_{\omega_{+,3},V})\}_{j=1}^4$}
\label{fig12}
\end{figure}

Let $r=0.3$. The circle 
$$\Xi_{\omega_{+,3}, V}=\{-0.5+(0.5+r)e^{2\pi ti}\}_{0\leq t< 1}$$
circular to $\{0, -1\}$ and $\{\infty, 1\}$ in Figure \ref{fig12} is an (OR) equator for ($R_{\omega_{+,3}}^2$) $R_{\omega_{+,3}}^4$.   

\begin{figure}[h]
\centering
\includegraphics[scale=1]{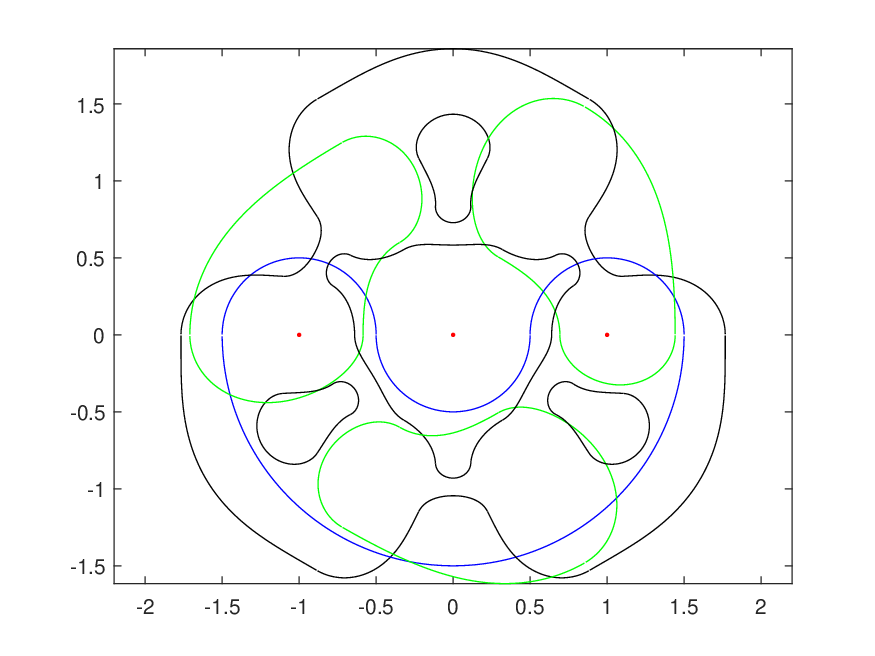}
\caption{$\Xi_{\omega_{+,3},VI}$, $R_{\omega_{+,3}}^{-1}(\Xi_{\omega_{+,3},VI})$ and $R_{\omega_{+,3}}^{-2}(\Xi_{\omega_{+,3},VI}$)}
\label{fig13}
\end{figure}

Let $r=0.5$. We illustrate the readers how the pre-images of the Jordan curve
\begin{center}
$
\begin{array}{ll}
\Xi_{\omega_{+,3},VI}=& \{-1+re^{2\pi it}\}_{0\leq t\leq 0.5}\cup\{(1+r)e^{2\pi it}\}_{0.5\leq t\leq 1}\vspace{3mm}\\
& \cup\{1+re^{2\pi it}\}_{0\leq t\leq 0.5}\cup\{re^{2\pi it}\}_{0.5\leq t\leq 1}
\end{array}
$
\end{center}
circular to $\{0, \infty\}$ and $\{-1, 1\}$ split in Figure \ref{fig13}.

\begin{figure}[h]
\centering
\includegraphics[scale=1]{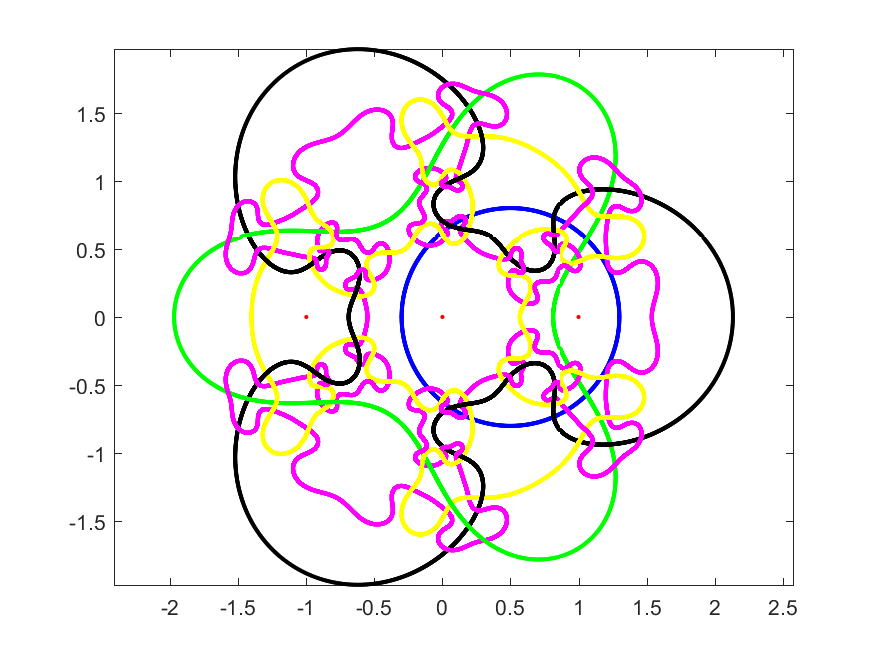}
\caption{$\Xi_{\omega_{+,3},VII}$ and $\{R_{\omega_{+,3}}^{-j}(\Xi_{\omega_{+,3},VII})\}_{j=1}^4$}
\label{fig14}
\end{figure}

Let $r=0.3$. The circle 
$$\Xi_{\omega_{+,3}, VII}=\{0.5+(0.5+r)e^{2\pi ti}\}_{0\leq t< 1}$$
circular to $\{0, 1\}$ and $\{\infty, -1\}$ in Figure \ref{fig14} is an (OR) equator for ($R_{\omega_{+,3}}^2$) $R_{\omega_{+,3}}^4$.   

In the following we turn to the bicritical maps with parametrizations in $\Omega^k_{-}$. We start from the degree $2$ case. In this case these bicritical maps take very simple formulas. 
\begin{lemma}\label{lem15}
We have 
$$\Omega^2_{-}=\{(0,b,c,0) \in \mathbb{C}^4: bc=-1\}.$$
So any bicritical map $R_\omega$ with $\omega\in\Omega^2_{-}$ takes the form
$$R_\omega=-\cfrac{b^2}{z^2}$$
for some $b\in (\mathbb{C}\setminus\{0\})$.
\end{lemma}
\begin{proof}
Let $\omega=(a,b,c,d)\in\Omega^2_{-}$. Then these entries must satisfy
\begin{equation}\label{eq44}
\left\{
\begin{array}{lllll}
ad-bc=1,\\
ab+d^2=0,\\
a^2+dc=0,\\
cb^2+d^3\neq 0,\\
a^3+bc^2\neq 0.
\end{array}
\right.
\end{equation} 
The second and third equality in (\ref{eq44}) imply that if one of the entries $a$ and $d$ equals $0$, then the other one equals $0$. We now show that a solution of (\ref{eq44}) must satisfy 
\begin{equation}\label{eq46}
a=d=0.
\end{equation}
We do this by reduction to absurdity. Suppose $\omega=(a,b,c,d)\in\Omega^2_{-}$ satisfies $a\neq 0$ and $d\neq 0$ simultaneously. Then the second and third equality in (\ref{eq44}) induce
\begin{equation}\label{eq45}
abcd=(ad)^2.
\end{equation}
We get $$ad-bc=0$$ from (\ref{eq45}) since $ad\neq 0$. This contradicts the first equality in (\ref{eq44}), which justifies (\ref{eq46}). Under (\ref{eq46}) we get $bc=-1$ from the first equality of (\ref{eq44}) and $b\neq 0, c\neq 0$ from the last two inequalities of (\ref{eq44}) (or from $bc=-1$ itself).
\end{proof}

Interestingly, the case of $k=2$ is a degenerate case compared with cases of $k\geq 3$ for maps parametrized by $\Omega^k_{-}$. This phenomenon does not repeat in the maps parametrized by $\Omega^k_{+}$ as there are no degenerate maps in them for any $k\geq 2$.  In virtue of Lemma \ref{lem15} the following results are instant.
\begin{coro}\label{cor2}
For any $\omega\in\Omega^2_{-}$, we have 
$$P(R_\omega)=\{0, \infty\}.$$
\end{coro}

\begin{coro}\label{cor3}
Any bicritical map $R_\omega$ with $\omega\in\Omega^2_{-}$ is an OR mating.
\end{coro}

\begin{proof}
One can show that any map $R_{0,b,-b^{-1},0}=-\cfrac{b^2}{z^2}$ with $b\in (\mathbb{C}\setminus\{0\})$ does not admit any equator, while any circle centred at $0$ with small enough radius for fixed $b\in (\mathbb{C}\setminus\{0\})$ can serve as an OR equator for $R_{0,b,-b^{-1},0}$. Then the conclusion follows from Lemma \ref{lem15} instantly.
\end{proof}

We continue to explore the unmatability of the maps parametrised by $\Omega^3_{-}$ in the following. We choose the representative map
$$R_{\omega_{-,3}}=\cfrac{-z^3+i}{iz^3-1}$$
for $\omega_{-,3}=(-\frac{1}{\sqrt{2}}, \frac{i}{\sqrt{2}}, \frac{i}{\sqrt{2}}, -\frac{1}{\sqrt{2}})\in \Omega_-^3$.
The critical points of the map $R_{\omega_{-,3}}$ go as
\begin{center}
\begin{tikzcd}
\ 0\ \arrow[r, bend left, "R_{\omega_{-,3}}"] & -i \arrow[l, bend left, "R_{\omega_{-,3}}"],
\end{tikzcd}
\end{center}
\begin{center}
\begin{tikzcd}
\infty \arrow[r, bend left, "R_{\omega_{-,3}}"] & i \arrow[l, bend left, "R_{\omega_{-,3}}"].
\end{tikzcd}
\end{center}

\begin{figure}[h]
\centering
\includegraphics[scale=1]{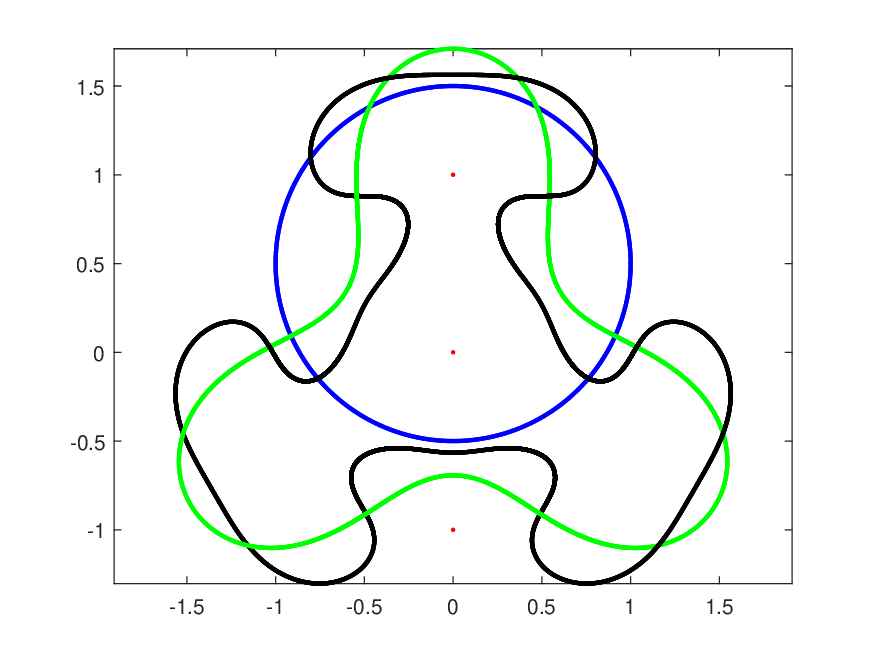}
\caption{$\Xi_{\omega_{-,3},V}$ and $\{R_{\omega_{-,3}}^{-j}(\Xi_{\omega_{-,3},V})\}_{j=1}^2$}
\label{fig15}
\end{figure}

\begin{figure}[h]
\centering
\includegraphics[scale=1]{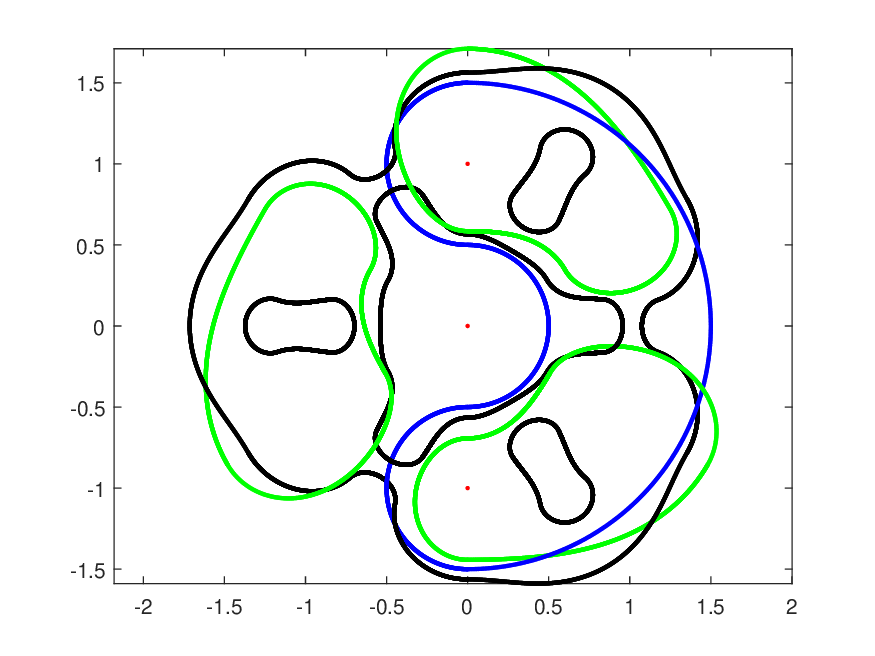}
\caption{$\Xi_{\omega_{-,3},VI}$ and $\{R_{\omega_{-,3}}^{-j}(\Xi_{\omega_{-,3},VI})\}_{j=1}^2$}
\label{fig16}
\end{figure}

\begin{figure}[h]
\centering
\includegraphics[scale=1]{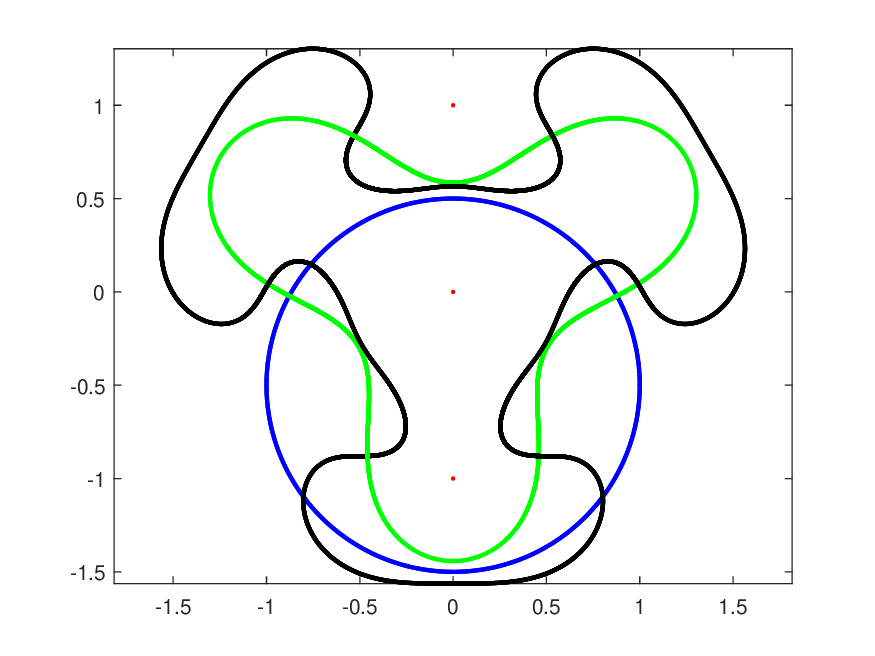}
\caption{$\Xi_{\omega_{-,3},VII}$ and $\{R_{\omega_{-,3}}^{-j}(\Xi_{\omega_{-,3},VII})\}_{j=1}^2$}
\label{fig17}
\end{figure}

\begin{pro}\label{pro5}
The map $R_{\omega_{-,3}}$ is a hermaphroditic mating.
\end{pro}
\begin{proof}
For any $j\in\mathbb{N}_+$, $P(R_{\omega_{-,3}}^j)=P(R_{\omega_{-,3}})=\{0, -i, \infty, i\}$. When $j$ is odd the critical points follow the same orbits of critical points under $R_{\omega_{-,3}}$. When $j$ is even the critical points go as
\begin{center}
\begin{tikzpicture}
    \node (a) at (0,0) {0};
    \draw[->] (a) to[in=45, out=-45, looseness=10] node{\ \ \ \ \ \ \  $R_{\omega_{-,3}}^j,$} (a);
    \node (b) at (3,0) {-i};
    \draw[->] (b) to[in=45, out=-45, looseness=5] node{\ \ \ \ \ \ \  $R_{\omega_{-,3}}^j,$} (b);
     \node (c) at (6,0) {$\infty$};
    \draw[->] (c) to[in=45, out=-45, looseness=10] node{\ \ \ \ \ \ \  $R_{\omega_{-,3}}^j,$} (c);
    \node (d) at (9,0) {i};
    \draw[->] (d) to[in=45, out=-45, looseness=10] node{\ \ \ \ \ \ \  $R_{\omega_{-,3}}^j.$} (d);
\end{tikzpicture}
\end{center}
To look for an equator or OR equator of $R_{\omega_{-,3}}^j$ for $j\in\mathbb{N}_+$, we consider partitions of the postcritical set into two subsets satisfying the immune or swapping property.  We list all the non-travail partitions of $P(R_{\omega_{-,3}}^j)$ into two subsets in Table \ref{tab5}.
\begin{table}[ht]
\caption{Partitions of $P(R_{\omega_{-,3}}^j)$} 
\centering 
\begin{tabular}{c c c} 
\hline 
Type &  $P_w(R_{\omega_{-,3}}^j)$ & $P_b(R_{\omega_{-,3}}^j)$ \\
\hline
I     &   $0$    &    $-i, \infty, i$\\

II   &    $-i$    &    $0, \infty, i$\\

III     &   $\infty$    &    $0, -i, i$\\

IV     &   $i$    &    $0, -i, \infty$\\

V     &   $0, -i$    &    $\infty, i$\\

VI   &    $0, \infty$    &    $-i,  i$\\

VII     &   $0, i$    &    $-i, \infty$\\
\hline 
\end{tabular}
\label{tab5} 
\end{table} 
One can see that all the listed partitions satisfy the immune property under $R_{\omega_{-,3}}^2$. The Type V partition in Table \ref{tab5} satisfies the immune  property under $R_{\omega_{-,3}}$. The Type VI and VII partitions in Table \ref{tab5} satisfy the swapping  property under $R_{\omega_{-,3}}$. We deal with the unmatability of $R_{\omega_{-,3}}^j$ according to the listed types of partitions.
\begin{itemize}
\item  Type I-IV. In these cases there are no equators for $R_{\omega_{-,3}}^j$ for any $j\in\mathbb{N}_+$ since the pre-image of any Jordan curve inducing these types of partitions  splits under $R_{\omega_{-,3}}$.  

\item  Type V. In this case there are equators for $R_{\omega_{-,3}}$. In fact any Jordan curve inducing the Type V partition in Table \ref{tab5} is an equator for  $R_{\omega_{-,3}}$.  One can refer to Figure \ref{fig15} for a concrete one $\Xi_{\omega_{-,3}, V}$.

\item Type VI. In this case there are no equators or OR equators for $R_{\omega_{-,3}}$. In fact the pre-image $R_{\omega_{-,3}}^{-1}(\Xi)$ of any Jordan curve $\Xi$ circular to $\{0, \infty\}$  and $\{-i,  i\}$ splits into three disjoint connected components. One can refer to Figure \ref{fig16} for the splitting process of some concrete Jordan curve $\Xi_{\omega_{-,3}, VI}$.

\item  Type VII. In this case there are OR equators for $R_{\omega_{-,3}}$. In fact any Jordan curve inducing the Type VII partition in Table \ref{tab5} is an OR equator for  $R_{\omega_{-,3}}$. One can refer to Figure \ref{fig17} for a concrete one $\Xi_{\omega_{-,3}, VII}$.
\end{itemize} 
\end{proof}

Let $r=0.5$. The circle 
$$\Xi_{\omega_{-,3}, V}=\{0.5i+re^{2\pi ti}\}_{0\leq t< 1}$$
inducing Type V partition in Table \ref{tab5} in Figure \ref{fig15} is an equator for $R_{\omega_{-,3}}$.   

Let $r=0.5$. We illustrate the readers how the pre-images of the Jordan curve
\begin{center}
$
\begin{array}{ll}
\Xi_{\omega_{-,3},VI}=& \{i+re^{2\pi it}\}_{0.25\leq t\leq 0.75}\cup\{(1-r)e^{2\pi it}\}_{-0.25\leq t\leq 0.25}\vspace{3mm}\\
& \cup\{-i+re^{2\pi it}\}_{0.25\leq t\leq 0.75}\cup\{(1+r)e^{2\pi it}\}_{-0.25\leq t\leq 0.25}\vspace{3mm}
\end{array}
$
\end{center}
inducing Type VI partition in Table \ref{tab5} split in Figure \ref{fig16}.

Let $r=0.5$. The circle 
$$\Xi_{\omega_{-,3}, VII}=\{-0.5i+(0.5+r)e^{2\pi ti}\}_{0\leq t< 1}$$
inducing Type VII partition in Table \ref{tab5} in Figure \ref{fig17} is an OR equator for $R_{\omega_{-,3}}$.   

As the final special case in this section, we consider unmatability of the maps parametrised by $\Omega^4_{-}$. We choose the representative map
$$R_{\omega_{-,4}}=\cfrac{z^4-\beta_{3,1}}{z^4-1}$$
for $\omega_{-,4}=(\frac{1}{\sqrt{\beta_{3,1}-1}}, -\frac{\beta_{3,1}}{\sqrt{\beta_{3,1}-1}}, \frac{1}{\sqrt{\beta_{3,1}-1}}, -\frac{1}{\sqrt{\beta_{3,1}-1}})\in \Omega_-^4$.
The critical points of the map $R_{\omega_{-,4}}$ go as
\begin{center}
\begin{tikzcd}
\ \ \ 0 \ \ \ \arrow[r, bend left, "R_{\omega_{-,4}}"] & \beta_{3,1} \arrow[l, bend left, "R_{\omega_{-,4}}"],
\end{tikzcd}
\end{center}
\begin{center}
\begin{tikzcd}
\infty \arrow[r, bend left, "R_{\omega_{-,4}}"] & 1 \arrow[l, bend left, "R_{\omega_{-,4}}"].
\end{tikzcd}
\end{center}

It turns out that the unmating property of the map  $R_{\omega_{-,4}}$ is the same as that of $R_{\omega_{-,3}}$.

\begin{figure}[h]
\centering
\includegraphics[scale=1]{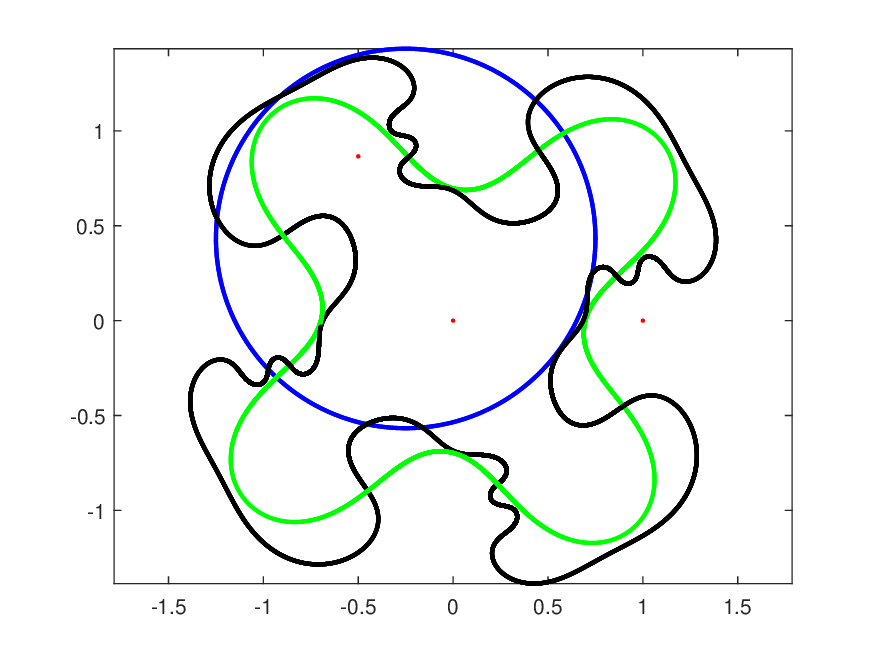}
\caption{$\Xi_{\omega_{-,4},V}$ and $\{R_{\omega_{-,4}}^{-j}(\Xi_{\omega_{-,4},V})\}_{j=1}^2$}
\label{fig18}
\end{figure}

\begin{figure}[h]
\centering
\includegraphics[scale=1]{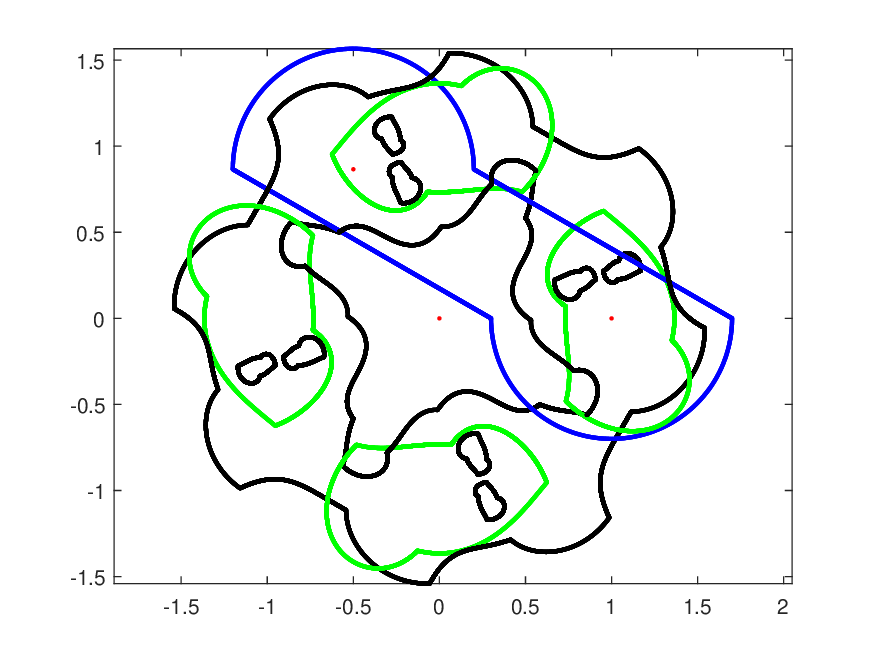}
\caption{$\Xi_{\omega_{-,4},VI}$ and $\{R_{\omega_{-,4}}^{-j}(\Xi_{\omega_{-,4},VI})\}_{j=1}^2$}
\label{fig19}
\end{figure}

\begin{figure}[h]
\centering
\includegraphics[scale=1]{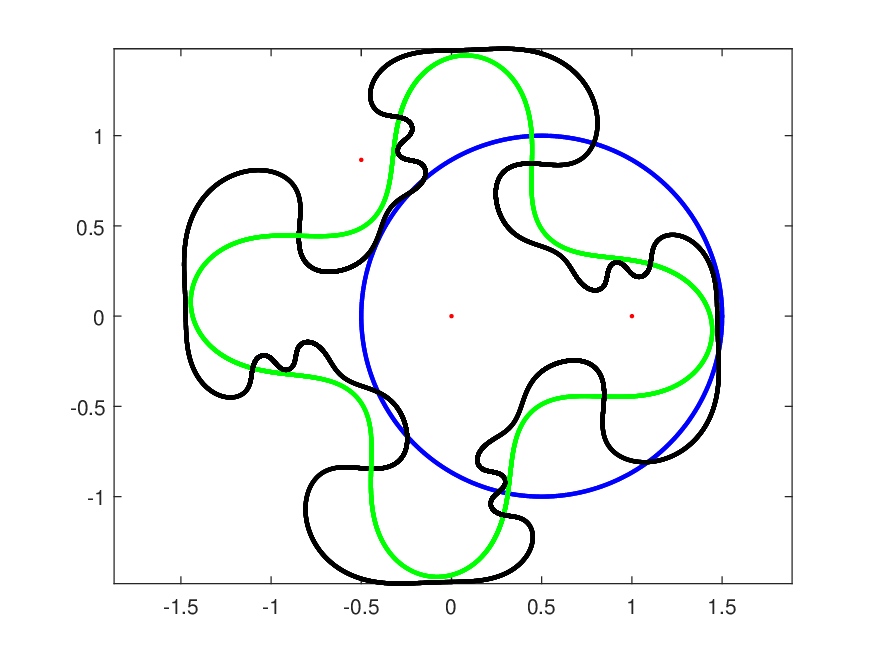}
\caption{$\Xi_{\omega_{-,4},VII}$ and $\{R_{\omega_{-,4}}^{-j}(\Xi_{\omega_{-,4},VII})\}_{j=1}^2$}
\label{fig20}
\end{figure}

\begin{pro}
The map $R_{\omega_{-,4}}$ is a hermaphroditic mating.
\end{pro}
\begin{proof}
The proof follows almost in line with the proof of Proposition \ref{pro5}. One can refer to Figure \ref{fig18} for an equator $\Xi_{\omega_{-,4}, V}$ and Figure \ref{fig20} for an OR equator $\Xi_{\omega_{-,4}, VII}$ of $R_{\omega_{-,4}}$. Any Jordan curve circular to $\{0, \infty\}$ and $\{\beta_{3,1}, 1\}$ splits under $R_{\omega_{-,4}}^{-1}$, see Figure \ref{fig19} for the splitting process of a concrete one $\Xi_{\omega_{-,4}, VI}$.
\end{proof}

Let $r=0.5$. The circle 
$$\Xi_{\omega_{-,4}, V}=\{0.5\beta_{3,1}+(0.5+r)e^{2\pi ti}\}_{0\leq t< 1}$$
circular to $\{0, \beta_{3,1}\}$ and $\{\infty, 1\}$ in Figure \ref{fig18} is an equator for $R_{\omega_{-,4}}$.   

Let $r=0.7$. The pre-images of the Jordan curve
\begin{center}
$
\begin{array}{ll}
\Xi_{\omega_{-,4},VI}=& \{\beta_{3,1}+re^{2\pi it}\}_{0\leq t\leq 0.5}\cup\{(1-r)+t(\beta_{3,1}-1)\}_{0\leq t\leq 1}\vspace{3mm}\\
& \cup\{1+re^{2\pi it}\}_{0.5\leq t\leq 1}\cup\{(1+r)+t(\beta_{3,1}-1)\}_{0\leq t\leq 1}\vspace{3mm}
\end{array}
$
\end{center}
circular to $\{0, \infty\}$ and $\{\beta_{3,1}, 1\}$ in Figure \ref{fig19} split under $R_{\omega_{-,4}}$.

Let $r=0.5$. The circle 
$$\Xi_{\omega_{-,4}, VII}=\{0.5+(0.5+r)e^{2\pi ti}\}_{0\leq t< 1}$$
circular to $\{0, 1\}$ and $\{\infty, \beta_{3,1}\}$ in Figure \ref{fig20} is an OR equator for $R_{\omega_{-,4}}$.   

\subsection{The proof of Theorem \ref{thm8}}

Along with knowledge on the pioneering cases above, we are well prepared to present the proof of Theorem \ref{thm8} now. We first show that $\Omega_+^k$ and $\Omega_-^k$ can be simplified to the following Riemann surfaces for any $k\geq 2$.

\begin{lemma}\label{lem16}
For any integer $k\geq 2$, we have
\begin{equation}\label{eq51}
\Omega_+^k=\Omega_1^k\cap\Omega_2^k
\end{equation}
and 
\begin{equation}\label{eq52}
\Omega_-^k=\Omega_3^k\cap\Omega_4^k.
\end{equation}
\end{lemma}

\begin{proof}
We first claim that the following system of equations 
\begin{equation}\label{eq47}
\left\{
\begin{array}{lllll}
ad-bc=1,\\
cb^k+d^{k+1}= 0,\\
ab^{k-1}+d^k=0.
\end{array}
\right.
\end{equation} 
admits no solutions. We justify the claim according to whether $b=0$ or not.
\begin{itemize}
\item Case $b=0$.  If $b=0$, we get $d=0$ from the second equality of (\ref{eq47}). These cause incompatibility with the first equality of (\ref{eq47}).
\item Case $b\neq 0$. If $b\neq 0$, we get $c=-\cfrac{d^{k+1}}{b^k}$ from the second equality of (\ref{eq47}) and $a=-\cfrac{d^k}{b^{k-1}}$ from the third equality of (\ref{eq47}). Then we have
$$ad-bc=-\cfrac{d^k}{b^{k-1}}d-(-\cfrac{d^{k+1}}{b^k})b=0,$$
which again causes incompatibility with the first equality of (\ref{eq47}).
\end{itemize}

The claim guarantees 
\begin{equation}\label{eq49}
\Omega_1^k\cap \Omega_3^k=\emptyset
\end{equation}
for any $k\geq 2$. By similar arguments we can show that the following system of equations 
\begin{equation}\label{eq48}
\left\{
\begin{array}{lllll}
ad-bc=1,\\
a^{k+1}+bc^k= 0,\\
a^k+dc^{k-1}=0.
\end{array}
\right.
\end{equation} 
admits no solutions. This guarantees 
\begin{equation}\label{eq50}
\Omega_2^k\cap \Omega_4^k=\emptyset
\end{equation}
for any $k\geq 2$. Now combining (\ref{eq49}) and (\ref{eq50}) together we get (\ref{eq51}). An analogous proof applies to (\ref{eq52}), whose details are left to the enthusiastic readers.
\end{proof}

\begin{proof}[Proof of Theorem \ref{thm8}]
We start from proving Theorem \ref{thm8} (A).  For fixed $k\geq 2$, any two maps with parameters in $\Omega_+^k$ ($\Omega_-^k$) are conformally conjugate to each other in virtue of the B\"ottcher's Theorem (extensions from the Fatou sets to the Julia sets are guaranteed since the Julia sets of hyperbolic maps are always locally connected), so it suffices for us to consider an individual map with parameters in $\Omega_+^k$ ($\Omega_-^k$). Let $R_{\omega_{+,k}}$ be some postcritically finite bicritical map with $\omega_{+,k}=(a,b,c,d)\in \Omega_+^k$ for some $k\geq 2$.   Note that we always have
\begin{center}
$\#P(R_{\omega_{+,k}})=4$
\end{center}
in this case. The postcritical points of $R_{\omega_{+,k}}$ go as
\begin{center}
\begin{tikzcd}
0 \arrow[r, "R_{\omega_{+,k}}"] & b/d \arrow[r, "R_{\omega_{+,k}}"]  & \infty \arrow[r, "R_{\omega_{+,k}}"] & a/c \arrow[lll, bend left, "R_{\omega_{+,k}}"].
\end{tikzcd}
\end{center}
So the postcritical points of $R_{\omega_{+,k}}^{4l+1}$ go as
\begin{center}
\begin{tikzcd}
0 \arrow[r, "R_{\omega_{+,k}}^{4l+1}"] & b/d \arrow[r, "R_{\omega_{+,k}}^{4l+1}"]  & \infty \arrow[r, "R_{\omega_{+,k}}^{4l+1}"] & a/c \arrow[lll, bend left, "R_{\omega_{+,k}}^{4l+1}"]
\end{tikzcd}
\end{center}
for any integer $l\geq 0$. The postcritical points of $R_{\omega_{+,k}}^{4l+2}$ go as
\begin{center}
\begin{tikzcd}
\ 0\ \arrow[r, bend left, "R_{\omega_{+,k}}^{4l+2}"] & \infty \arrow[l, bend left, "R_{\omega_{+,k}}^{4l+2}"],
\end{tikzcd}
\begin{tikzcd}
b/d \arrow[r, bend left, "R_{\omega_{+,k}}^{4l+2}"]  &  a/c \arrow[l, bend left, "R_{\omega_{+,k}}^{4l+2}"]
\end{tikzcd}
\end{center}
for any $l\geq 0$. The postcritical points of $R_{\omega_{+,k}}^{4l+3}$ go as
\begin{center}
\begin{tikzcd}
0 \arrow[r, "R_{\omega_{+,k}}^{4l+3}"] & a/c \arrow[r, "R_{\omega_{+,k}}^{4l+3}"]  & \infty \arrow[r, "R_{\omega_{+,k}}^{4l+3}"] &  b/d \arrow[lll, bend left, "R_{\omega_{+,k}}^{4l+3}"],
\end{tikzcd}
\end{center}
for any $l\geq 0$. The postcritical points of $R_{\omega_{+,k}}^{4l+4}$ go as
\begin{center}
\begin{tikzpicture}
    \node (a) at (0,0) {0};
    \draw[->] (a) to[in=45, out=-45, looseness=10] node{\ \ \ \ \ \ \  $R_{\omega_{+,k}}^{4l+4},$} (a);
    \node (b) at (3,0) {$b/d$};
    \draw[->] (b) to[in=45, out=-45, looseness=5] node{\ \ \ \ \ \ \  $R_{\omega_{+,2}}^{4l+4},$} (b);
     \node (c) at (6,0) {$\infty$};
    \draw[->] (c) to[in=45, out=-45, looseness=10] node{\ \ \ \ \ \ \  $R_{\omega_{+,2}}^{4l+4},$} (c);
    \node (d) at (9,0) {$a/c$};
    \draw[->] (d) to[in=45, out=-45, looseness=10] node{\ \ \ \ \ \ \  $R_{\omega_{+,2}}^{4l+4}$} (d);
\end{tikzpicture}
\end{center}
for any $l\geq 0$. To look for an equator or OR equator of $R_{\omega_{+,k}}^j$ for $j\geq 1$, we consider partitions of the postcritical set into two subsets satisfying the immune or swapping property.  We list all the non-travail partitions of $P(R_{\omega_{+,k}}^j)$ into two subsets in Table \ref{tab6}.
\begin{table}[ht]
\caption{Partitions of $P(R_{\omega_{+,k}}^j)$} 
\centering 
\begin{tabular}{c c c} 
\hline 
Type &  $P_w(R_{\omega_{+,k}}^j)$ & $P_b(R_{\omega_{+,k}}^j)$ \\
\hline
I     &   $0$    &    $b/d, \infty, a/c$\\

II   &    $b/d$    &    $0, \infty, a/c$\\

III     &   $\infty$    &    $0, b/d, a/c$\\

IV     &   $a/c$    &    $0, b/d, \infty$\\

V     &   $0,  b/d$    &    $\infty, a/c$\\

VI   &    $0, \infty$    &    $b/d,  a/c$\\

VII     &   $0, a/c$    &    $b/d, \infty$\\
\hline 
\end{tabular}
\label{tab6} 
\end{table} 
One can see that all the listed partitions satisfy the immune property under $R_{\omega_{+,k}}^4$. The Type V and VII partition in Table \ref{tab6} satisfy the swapping  property under $R_{\omega_{+,k}}^2$. The Type VI partition in Table \ref{tab6} satisfies the immune  property under $R_{\omega_{+,k}}^2$, while it satisfies the swapping property under $R_{\omega_{+,k}}$ and $R_{\omega_{+,k}}^3$.  We deal with the unmatability of $R_{\omega_{+,k}}^j$ according to the listed types of partitions in Table \ref{tab6}.
\begin{itemize}
\item  Type I-IV. In these cases there are no equators for $R_{\omega_{+,k}}^j$ for any $j\in\mathbb{N}_+$ dividable by $4$. Let $\Xi$ be a Jordan curve which partitions the postcritical set into one of the four types mentioned before. $R_{\omega_{+,k}}^{-4}(\Xi)$ will split into more than two disjoint connected components, because either of $R_{\omega_{+,k}}^{-1}(0)$ and $R_{\omega_{+,k}}^{-1}(\infty)$ splits into at least two different points. Thus $\Xi$ can not serve as an equator (or OR equator) for $R_{\omega_{+,k}}^j$ for any $j\in\mathbb{N}_+$ dividable by $4$.    

\item  Type V. In this case there are OR equators for $R_{\omega_{+,k}}^2$ and equators for $R_{\omega_{+,k}}^4$ . Let $\Xi$ be any Jordan curve which is circular to both 
\begin{center}
$P_w(R_{\omega_{+,k}}^2)=\{0, b/d\}$ and $P_b(R_{\omega_{+,k}}^2)=\{\infty, a/c\}$,
\end{center}
such that $U_w^0$ is a convex set in $\mathbb{P}^1(\mathbb{C})$.  Note that $R_{\omega_{+,k}}^{-2}(\Xi)$ is always connected  and circular to both $P_w(R_{\omega_{+,k}}^2)$ and $P_b(R_{\omega_{+,k}}^2)$, because $R_{\omega_{+,k}}^{-2}(U_w^0)=U_w^2$ is symmetric with respect to $0$ and $0\in R_{\omega_{+,k}}^{-2}(U_w^0)$. Thus $\Xi$ is an OR equator for $R_{\omega_{+,k}}^2$ considering the swapping property of $P_w(R_{\omega_{+,k}}^2)$ and $P_b(R_{\omega_{+,k}}^2)$. It is also an equator for $R_{\omega_{+,k}}^4$ in virtue of Lemma \ref{lem4}. 

\item Type VI. In this case there are no equators (or OR equators) for $R_{\omega_{+,2}}^{4l+2}$ for any $l\geq 0$. In fact the pre-image $R_{\omega_{+,2}}^{-1}(\Xi)$ of any Jordan curve $\Xi$ circular to $P_w(R_{\omega_{+,k}})=\{0, \infty\}$  and $P_b(R_{\omega_{+,k}})=\{b/d,  a/c\}$ splits into $k$ disjoint connected components. To see this, note that $R_{\omega_{+,k}}^{-1}(U_w^0)=U_w^1$ is symmetric with respect to $0$ and $0\notin U_w^1$. This forces $R_{\omega_{+,k}}^{-1}(U_w^0)$ and hence $R_{\omega_{+,k}}^{-2}(\Xi)$ to split into $k$ disjoint connected components. 

\item  Type VII. In this case there are OR equators for $R_{\omega_{+,k}}^2$ and equators for $R_{\omega_{+,k}}^4$. This case is dual to the Type V case.  
\end{itemize} 
After exhausting all the types of partitions, we can conclude that $R_{\omega_{+,k}}$ is a (OR) $4$-mating for any $\omega_{+,k}\in \Omega_+^k$ and $k\geq 2$.

Then we prove Theorem \ref{thm8} (B). Again it suffices for us to consider an individual map with parametrization in $\Omega_-^k$ for $k\geq 3$. Let $R_{\omega_{-,k}}$ be some postcritically finite bicritical map with $\omega_{-,k}=(a,b,c,d)\in \Omega_-^k$ for some $k\geq 3$. Note that we always have
\begin{center}
$\#P(R_{\omega_{-,k}})=4$
\end{center}
in this case. The postcritical points of $R_{\omega_{-,k}}^{2l+1}$ go as 
\begin{center}
\begin{tikzcd}
\ \ 0\ \  \arrow[r, bend left, "R_{\omega_{-,k}}^{2l+1}"] & b/d \arrow[l, bend left, "R_{\omega_{-,k}}^{2l+1}"],
\end{tikzcd}
\begin{tikzcd}
\ \infty\ \arrow[r, bend left, "R_{\omega_{-,k}}^{2l+1}"]  &  a/c \arrow[l, bend left, "R_{\omega_{-,k}}^{2l+1}"]
\end{tikzcd}
\end{center}
for any integer $l\geq 0$. The postcritical points of $R_{\omega_{-,k}}^{2l+2}$ go as
\begin{center}
\begin{tikzpicture}
    \node (a) at (0,0) {0};
    \draw[->] (a) to[in=45, out=-45, looseness=10] node{\ \ \ \ \ \ \  $R_{\omega_{-,k}}^{2l+2},$} (a);
    \node (b) at (3,0) {$b/d$};
    \draw[->] (b) to[in=45, out=-45, looseness=5] node{\ \ \ \ \ \ \  $R_{\omega_{-,k}}^{2l+2},$} (b);
     \node (c) at (6,0) {$\infty$};
    \draw[->] (c) to[in=45, out=-45, looseness=10] node{\ \ \ \ \ \ \  $R_{\omega_{-,k}}^{2l+2},$} (c);
    \node (d) at (9,0) {$a/c$};
    \draw[->] (d) to[in=45, out=-45, looseness=10] node{\ \ \ \ \ \ \  $R_{\omega_{-,k}}^{2l+2}$} (d);
\end{tikzpicture}
\end{center}
for any $l\geq 0$. To look for an equator or OR equator for $R_{\omega_{-,k}}^j$, we consider partitions of the postcritical set into two subsets satisfying the immune or swapping property.  The list is the same as Table \ref{tab6}. One can see that all the listed partitions in Table \ref{tab6} satisfy the immune property under $R_{\omega_{-,k}}^2$. The Type V partition in Table \ref{tab6} satisfies the immune  property under $R_{\omega_{-,k}}$. The Type VI  and  VII partitions in Table \ref{tab6} satisfy the swapping  property under $R_{\omega_{-,k}}$. We deal with the unmatability of $R_{\omega_{-,k}}^j$ according to the listed types of partitions in Table \ref{tab6}.
\begin{itemize}
\item  Type I-IV. In these cases there are no equators for $R_{\omega_{-,k}}^j$ for any even $j\in\mathbb{N}_+$ since the pre-image of any Jordan curve inducing these types of partitions  splits under $R_{\omega_{-,k}}^2$.  

\item  Type V. In this case there are equators for $R_{\omega_{-,k}}$. Let $\Xi$ be any Jordan curve which is circular to both 
\begin{center}
$P_w(R_{\omega_{-,k}})=\{0, b/d\}$ and $P_b(R_{\omega_{-,k}})=\{\infty, a/c\}$,
\end{center}
such that $U_w^0$ is a convex set in $\mathbb{P}^1(\mathbb{C})$.  Note that $R_{\omega_{-,k}}^{-1}(\Xi)$ is always connected  and circular to both $P_w(R_{\omega_{-,k}})$ and $P_b(R_{\omega_{-,k}})$, because $R_{\omega_{-,k}}^{-1}(U_w^0)=U_w^1$ is symmetric with respect to $0$ ($\infty$) and $0\in U_w^1$. Thus $\Xi$ is an equator for $R_{\omega_{-,k}}$ considering the immune property of $P_w(R_{\omega_{-,k}})$ and $P_b(R_{\omega_{-,k}})$.

\item Type VI. In this case there are no equators or OR equators for $R_{\omega_{-,k}}$. In fact the pre-image $R_{\omega_{-,k}}^{-1}(\Xi)$ of any Jordan curve $\Xi$ 
circular to both 
\begin{center}
$P_w(R_{\omega_{-,k}})=\{0, \infty\}$ and $P_b(R_{\omega_{-,k}})=\{b/d, a/c\}$,
\end{center}
splits into $k$ disjoint connected components. This is because $U_w^1$ is symmetric with respect to $0$ ($\infty$) and $0\notin U_w^1$.

\item  Type VII. In this case there are OR equators for $R_{\omega_{-,k}}$. Let $\Xi$ be any Jordan curve which is circular to both 
\begin{center}
$P_w(R_{\omega_{-,k}})=\{0, a/c\}$ and $P_b(R_{\omega_{-,k}})=\{\infty, b/d\}$,
\end{center}
such that $U_w^0$ is a convex set in $\mathbb{P}^1(\mathbb{C})$.  Note that $R_{\omega_{-,k}}^{-1}(\Xi)$ is always connected  and still circular to both $P_w(R_{\omega_{-,k}})$ and $P_b(R_{\omega_{-,k}})$, because $U_w^1$ is symmetric with respect to $0$ and $0\in U_w^1$. Thus $\Xi$ is an OR equator for $R_{\omega_{-,k}}$ considering the swapping property of $P_w(R_{\omega_{-,k}})$ and $P_b(R_{\omega_{-,k}})$ under $R_{\omega_{-,k}}$. 
\end{itemize} 
After exhausting all the types of partitions, we can conclude that $R_{\omega_{-,k}}$ is a hermaphroditic mating for any $\omega_{-,k}\in \Omega_-^k$ with $k\geq 3$.

Finally, Theorem \ref{thm8} (C) follows directly  from Corollary \ref{cor3}.

\end{proof}

It would be interesting to consider the unmatability of other postcritically finite bicritical  maps in the connectedness locus $\EuScript{C}\subset\EuScript{M}$. There will be other types of matings besides the $4$-matings,  hermaphroditic matings and orientation-reversing matings, in the author's prediction.

\section{Atomic maps and primitive $2$-matings in the captures}\label{sec5}

In this section we explore the $n$-unmatability of the postcritically finite capture maps in the quadratic family $\Big\{R_a=\cfrac{a}{z^2+2z}\Big\}_{a\in \mathbb{C}\setminus\{0\}}$. This is a well-studied family by Aspenberg and Yampolsky (\cite{AY}). Let  $V_j$ be the collection of quadratic rational maps with a super-attracting period-$j$ cycle up to M\"obius conjugacies (\emph{cf.} \cite{Ree3}) for $j\in\mathbb{N}_+$. The family $\{R_a\}_{a\in \mathbb{C}\setminus\{0\}}$ can be viewed as a slice of $V_2$ in fact. All maps in the family admit two critical points $\{-1, \infty\}$. If there exists some $k\in\mathbb{N}_+$ such that 
$$R_a^k(-1)=\infty$$
for some $a\in \mathbb{C}\setminus\{0\}$,  then the map $R_a$ is called a (postcritically finite) \emph{capture} map, with $k$ being its generation.  Parameters of the postcritically finite captures are situated at centres of the parabubbles (capture hyperbolic components) in the parameter slice.  
   
\subsection{Captures are not matings in the classical sense}
The following result can be deduced from  Aspenberg-Yampolsky's work \cite{AY}. However, here we provide a proof of different flavour.
\begin{theorem}[Aspenberg-Yampolsky]\label{lem9}
The postcritically finite captures in the family $\Big\{R_a=\cfrac{a}{z^2+2z}\Big\}_{a\in \mathbb{C}\setminus\{0\}}$ are not $1$-matings.
\end{theorem}
\begin{proof}
Assume  $a_*$ is at the center of a $k$-generation parabubble, such that the orbits of the critical points of $R_{a_*}$ are as following,
\begin{center}
\begin{tikzcd}
p_1 \arrow[r, "R_{a_*}"] & p_2 \arrow[r, "R_{a_*}"] & \cdots \arrow[r, "R_{a_*}"] & p_k \arrow[r, "R_{a_*}"]  & \infty \arrow[r, bend left, "R_{a_*}"] & 0 \arrow[l, bend left, "R_{a_*}"], 
\end{tikzcd}
\end{center}
in which $p_1=-1, p_k=-2$. Then its postcritical set is
\begin{center}
$P(R_{a_*})=\{p_2, \cdots, p_k, \infty, 0\}$.  
\end{center}
We claim that there does not exist any non-travail partition of $P(R_{a_*})$ into two subsets which obeys the immune property. This claim is enough to induce the conclusion. To see this, suppose that there exists some partition $\{P_w(R_{a_*}), P_b(R_{a_*})\}$ such that neither of the two sets is empty. Then according to property of the orbits of the critical points, in order to guarantee the partition satisfies the immune property under  $R_{a_*}$, $\infty$ and $0$ must be in the same subset, without loss of generality we assume $\{\infty, 0\}\subset P_w(R_{a_*})$. Then according to the immune property inductively we have
$$\{p_2,\cdots, p_k\}\subset P_w(R_{a_*}),$$
which forces $P_b(R_{a_*})=\emptyset$. This contradiction justifies the claim.  
\end{proof}

Theorem \ref{lem9} applies to any postcritically finite capture in any family in fact. This illustrates a typical constraint of the classical notion of mating.

\subsection{The unmatability of some captures of low generations}
Similar to the investigation on the $n$-unmatability of maps in the bicritical family, we start from some special cases first. We consider the unmatability of the postcritically finite capture $R_2=\cfrac{2}{z^2+2z}$ of generation $2$ first. This parameter sits at the center of the biggest parabubble in the parameter slice. The orbits of its two critical points $-1$ and $ \infty$ follow the graph
\begin{center}
\begin{tikzcd}
-1 \arrow[r,  "R_2"] & -2 \arrow[r, "R_2"] & \infty \arrow[r, bend left, "R_2"] & 0 \arrow[l, bend left, "R_2"]. 
\end{tikzcd}
\end{center}
Thus we have $P(R_2)=\{-2, \infty, 0\}$. $R_2$ is not a $1$-mating in virtue of Theorem \ref{lem9}. In fact any iterate of the map does not arise as a $1$-mating.

\begin{pro}\label{pro7}
$R_2$ is an atomic rational map.
\end{pro}

\begin{figure}[h]
\centering
\includegraphics[scale=1]{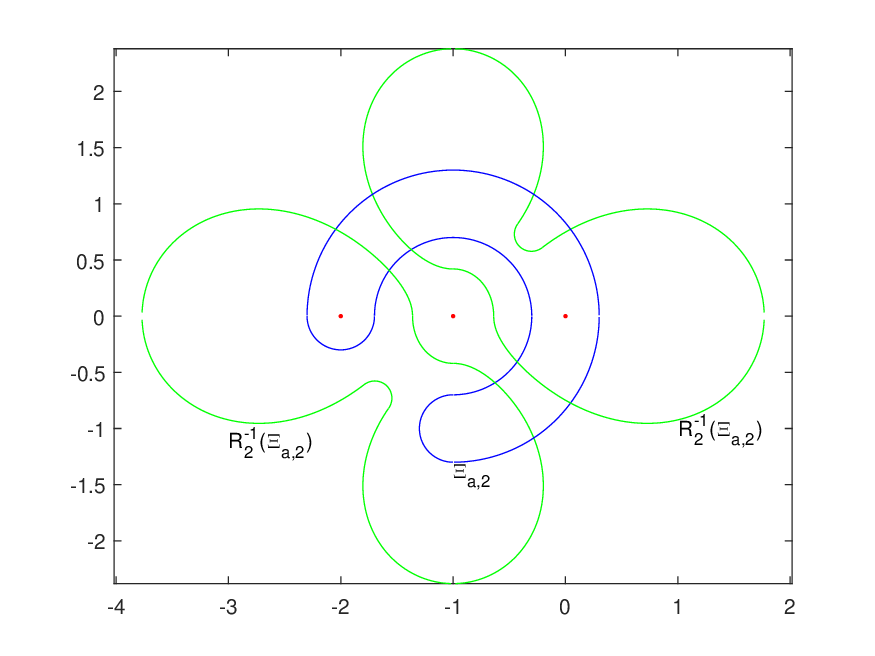}
\caption{$\Xi_{a,2}$ and $R_{2}^{-1}(\Xi_{a,2})$}
\label{fig8}
\end{figure}

\begin{proof}
One can apply \cite[Proposition 4.3]{Mey1} directly here since 
\begin{center}
$\#P(R_2^j)=\#P(R_2)=\#\{-2, \infty, 0\}=3$
\end{center} 
for any $j\in\mathbb{N}_+$. However, we will refine the proof here in order to illustrate the readers the reason that any of its iterations fails to be a $1$-mating.  For $j\in\mathbb{N}_+$, to look for an equator of $R_2^j$, we first search for partitions of $P(R_2^j)$ into two subsets satisfying the immune property. We list all the non-travail partitions of $P(R_2)$ into two subsets in Table \ref{tab1},
\begin{table}[ht]
\caption{Partitions of $P(R_2^j)$} 
\centering 
\begin{tabular}{c c c} 
\hline 
Type &  $P_w(R_2^j)$ & $P_b(R_2^j)$ \\
\hline
I     &   $-2$    &    $0, \infty$\\

II   &    $0$    &    $-2, \infty$\\

III     &   $\infty$    &    $0, -2$\\
\hline 
\end{tabular}
\label{tab1} 
\end{table}

Note that only the Type III partition satisfies the immune property under $R_2^j$ for even $j$. So an equator of $R_2^j$ must be circular to both $P_w(R_2^j)=\{\infty\}$ and $P_b(R_2^j)=\{0, -2\}$ for even $j$. However, any Jordan curve $\Xi$ circular to $\{\infty\}$ and $\{0, -2\}$ is isotopic to a small circle around $\{\infty\}$ \emph{rel.} $P(R_2^j)$, whose pre-image $(R_2^j)^{-1}(\Xi)$ splits into at least two disjoint connected components since $\#R_2^{-j}(\infty)\geq 2$ for any positive even $j$. This denies the possibility for any Jordan curve $\Xi$ to serve as an equator for $R_2^j$ (refer to Figure \ref{fig8} for the splitting process of the Jordan curve $\Xi_{a,2}$ inducing the Type III partition in  Table \ref{tab1} under $R_2^{-1}$).

\end{proof}

We provide the readers a concrete Jordan curve $\Xi_{a,2}$ circular to $P_b(R_2)\cup R_2^{-1}(-1)=\{-2, 0,-1\pm i\}$ and $P_w(R_2)\cup \{-1\}=\{\infty, -1\}$, as well as its pre-image under $R_{2}$ in Figure \ref{fig8}. Let $r=0.3$, the curve is 
\begin{center}
$
\begin{array}{ll}
\Xi_{a,2}=& \{-1+(1-r)e^{2\pi it}\}_{-0.25\leq t\leq 0.5}\cup\{-2+re^{2\pi it}\}_{0.5\leq t\leq 1}\vspace{3mm}\\
& \cup\{-1+(1+r)e^{2\pi it}\}_{-0.25\leq t\leq 0.5}\cup\{-1-i+re^{2\pi it}\}_{0.25\leq t\leq 0.75},
\end{array}
$
\end{center}
such that $-1\pm i\in U_b^0$ but $-1\notin U_b^0$.

In the following we continue to investigate the postcritically finite generation-$3$ captures in the family $\Big\{R_a=\cfrac{a}{z^2+2z}\Big\}_{a\in \mathbb{C}\setminus\{0\}}$. There is only one such map: $f_{3/2}=\cfrac{3/2}{z^2+2z}$. Note that $P(R_{3/2})=\{-3/2, -2, \infty, 0\}$. Its postcritical points go as
\begin{center}
\begin{tikzcd}
-1 \arrow[r,  "R_{3/2}"] & -3/2 \arrow[r, "R_{3/2}"] & -2 \arrow[r, "R_{3/2}"] & \infty \arrow[r, bend left, "R_{3/2}"] & 0 \arrow[l, bend left, "R_{3/2}"]. 
\end{tikzcd}
\end{center}
$R_{3/2}$ is not a $1$-mating in virtue of Theorem \ref{lem9}, however, surprisingly, it is a $2$-mating. More interestingly, it does not admit any OR equator.

\begin{pro}\label{thm6}
$R_{3/2}$ is a primitive $2$-mating.
\end{pro}

We aim to achieve Proposition \ref{thm6} in two steps. We first show that $R_{3/2}$ is not an orientation-reversing mating.

\begin{lemma}\label{lem10}
$R_{3/2}$ admits no OR equators. 
\end{lemma} 
\begin{proof}
To look for an OR equator of $R_{3/2}$, we consider partitions of $P(R_{3/2})$ satisfying the swapping property. We list all the non-travail partitions of $P(R_{3/2})$ into two subsets in Table \ref{tab2}.
\begin{table}[ht]
\caption{Partitions of $P(R_{3/2}^j)$} 
\centering 
\begin{tabular}{c c c} 
\hline 
Type &  $P_w(R_{3/2}^j)$ & $P_b(R_{3/2}^j)$ \\
\hline
I     &   $-3/2$    &    $-2, \infty, 0$\\

II   &    $-2$    &    $-3/2, \infty, 0$\\

III     &   $\infty$    &    $-3/2, -2, 0$\\

IV     &   $0$    &    $-3/2, -2, \infty$\\

V     &   $-3/2, -2$    &    $\infty, 0$\\

VI   &    $-3/2, \infty$    &    $-2,  0$\\

VII     &   $-3/2, 0$    &    $-2, \infty$\\
\hline 
\end{tabular}
\label{tab2} 
\end{table}

One can see that only the Type VI partition satisfies the swapping property under $R_{3/2}$. The Type VI partition is so alluring to induce an OR equator for $R_{3/2}$, unfortunately, this is wrong. Let $\Xi$ be a Jordan curve circular to $\{-3/2, \infty\}$ and $\{-2, 0\}$. $R_{3/2}^{-1}(\Xi)$ is always a single connected Jordan curve.  Let $U_w^0$ be the closure of the connected component of $\mathbb{C}\setminus \Xi$ containing $\{-3/2, \infty\}$ and let $U_b^0$ be the closure of the connected component of $\mathbb{C}\setminus \Xi$ containing $\{-2, 0\}$. Then 
\begin{center}
$\{-2, 0\}\subset R_{3/2}^{-1}(U_w^0)=U_w^1$ while $-3/2\notin R_{3/2}^{-1}(U_w^0)=U_w^1$.
\end{center}
Note that the capture $R_{3/2}$ preserves the real line while it swaps the upper and lower half complex plane. These force $-3/2$ to be confined in a domain enclosed by $U_w^0$ and $R_{3/2}^{-1}(U_b^0)=U_b^1$, which prevents $\Xi$ to be isotopic to $R_{3/2}^{-1}(\Xi)$ \emph{rel.} $P(R_{3/2})$.
\end{proof}

\begin{figure}[h]
\centering
\includegraphics[scale=1]{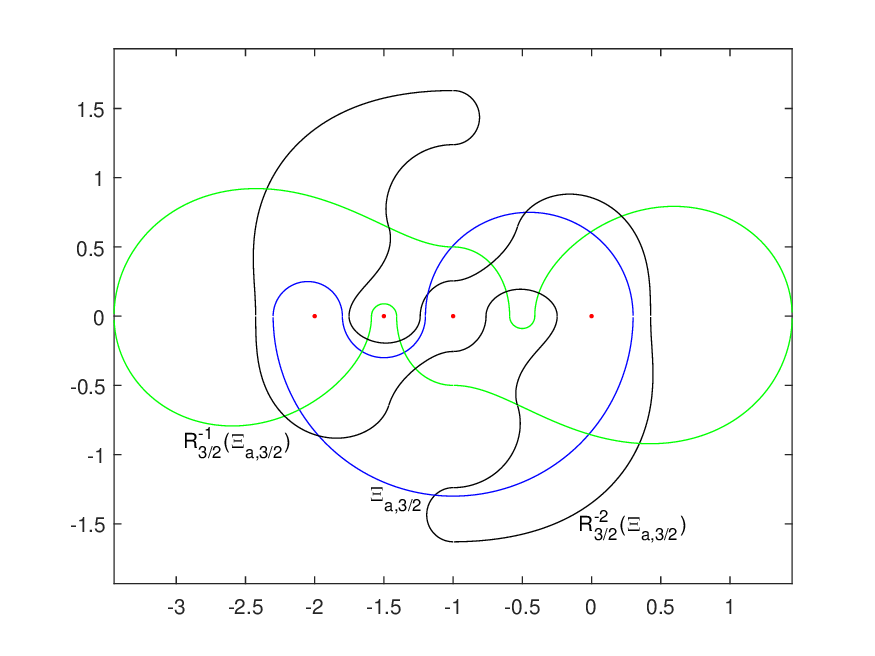}
\caption{$\Xi_{a,3/2}$, $R_{3/2}^{-1}(\Xi_{a,3/2})$ and $R_{3/2}^{-2}(\Xi_{a,3/2})$}
\label{fig4}
\end{figure}

Although $R_{3/2}$ does not admit any equator or OR equator, the unmatability changes after merely once iteration.
\begin{lemma}\label{lem11}
$R_{3/2}^2$ admits an equator. 
\end{lemma}
\begin{proof}
To look for an equator of $R_{3/2}^2$, we consider partitions of $P(R_{3/2}^2)$ satisfying the immune property. Since $P(R_{3/2}^2)=P(R_{3/2})$, we again consider the partitions in Table \ref{tab2}. Be careful that we are considering the iterate $R_{3/2}^2$ instead of $R_{3/2}$ now.  
Obviously  only the Type VI partition is possible to induce an equator of $R_{3/2}^2$. In fact any Jordan curve circular to $\{-2, 0\}$ and $\{-3/2, \infty\}$ serves as an equator for $R_{3/2}^2$, since $R_{3/2}^2$ preserves the upper and lower half complex plane now. One can  refer to Figure \ref{fig4} for a concrete one $\Xi_{a,3/2}$.
\end{proof}

The Jordan curve  
\begin{center}
$
\begin{array}{ll}
\Xi_{a,3/2}=& \{-1.5+0.3e^{2\pi it}\}_{0.5\leq t\leq 1}\cup\{-2.05+0.25e^{2\pi it}\}_{0\leq t\leq 0.5}\vspace{3mm}\\
& \cup\{-1+1.3e^{2\pi it}\}_{0.5\leq t\leq 1}\cup\{-0.45+0.75e^{2\pi it}\}_{0\leq t\leq 0.5}
\end{array}
$
\end{center}
inducing the Type VI partition of $P(R_{3/2}^2)$ in Table \ref{tab2} in Figure \ref{fig4} is an equator for $R_{3/2}^2$. 

\begin{proof}[Proof of Proposition \ref{thm6}]
The result follows from a combination of Lemma \ref{lem10}, Lemma \ref{lem11} and \cite[Theorem 4.2]{Mey1} instantly.
\end{proof}

Now we continue to investigate unmatability of the postcritically finite generation-$4$ captures in the family $\Big\{R_a=\cfrac{a}{z^2+2z}\Big\}_{a\in \mathbb{C}\setminus\{0\}}$. The number of them is more than one.
\begin{pro}
There are totally three postcritically finite generation-$4$ captures in the family 
$$\Big\{R_a=\cfrac{a}{z^2+2z}\Big\}_{a\in \mathbb{C}\setminus\{0\}},$$ 
which are parametrized by
\begin{center}
$a_{4,1}=1.36110308052864737763464656\cdots$,
\end{center}
\begin{center}
$a_{4,2}=1.319448459735676311182676718\cdots-1.633170240915237656118846731\cdots i$,
\end{center}
and
\begin{center}
$a_{4,3}=1.319448459735676311182676718\cdots+1.633170240915237656118846731\cdots i$ 
\end{center}
respectively. The orbits of their postcritical points are approximately
\begin{center}
\begin{tikzcd}
-1 \arrow[r,  "R_{a_{4,1}}"] & -1.3611 \arrow[r, "R_{a_{4,1}}"] & -1.5652 \arrow[r, "R_{a_{4,1}}"] & -2 \arrow[r, "R_{a_{4,1}}"] & \infty \arrow[r, bend left, "R_{a_{4,1}}"] & 0 \arrow[l, bend left, "R_{a_{4,1}}"],
\end{tikzcd}
\end{center}
\begin{center}
\begin{tikzcd}
-1 \arrow[r,  "R_{a_{4,2}}"] & -1.3194+1.6332i \arrow[r, "R_{a_{4,2}}"] & -0.2174 + 0.5217i \arrow[r, "R_{a_{4,2}}"] & -2 \arrow[r, "R_{a_{4,2}}"] & \infty \arrow[r, bend left, "R_{a_{4,2}}"] & 0 \arrow[l, bend left, "R_{a_{4,2}}"],
\end{tikzcd}
\end{center}
and
\begin{center}
\begin{tikzcd}
-1 \arrow[r,  "R_{a_{4,3}}"] & -1.3194-1.6332i \arrow[r, "R_{a_{4,3}}"] & -0.2174 - 0.5217i \arrow[r, "R_{a_{4,3}}"] & -2 \arrow[r, "R_{a_{4,3}}"] & \infty \arrow[r, bend left, "R_{a_{4,3}}"] & \ 0\ \arrow[l, bend left, "R_{a_{4,3}}"]
\end{tikzcd}
\end{center}
respectively.
\end{pro}
\begin{proof}
The forward iterates of the critical point $-1$ of a postcritically finite capture $R_a$ must enter the super attracting cycle 
\begin{tikzcd}
\infty \arrow[r, bend left, "R_a"] & \ 0\ \arrow[l, bend left, "R_a"] 
\end{tikzcd} 
via $-2$. Then one can get all the postcritically finite captures of generation-$4$ by solving the equation
$$R_a^3(-1)=\cfrac{a(a-2)^2}{2a-3}=-2$$
with respect to $a$. It turns out that there are three solutions. 
\end{proof}

Now we count the number of maps on the  postcritical set of $R_{a_{4,l}}^j$ under $R_{a_{4,l}}^j$ for any $l\in\{1,2,3\}$ and $j\in\mathbb{N}_+$. It turns out that there are only finitely many of them. 

\begin{lemma}\label{lem17}
For any $l\in\{1,2,3\}$ and $j\in\mathbb{N}_+$, the orbits of the postcritical points of $R_{a_{4,l}}^j$ follow exactly one of the following three maps,
\begin{center}
\begin{tikzcd}
R_{a_{4,l}}(-1) \arrow[r, "R_{a_{4,l}}"] & R_{a_{4,l}}^2(-1) \arrow[r, "R_{a_{4,l}}"] & -2 \arrow[r, "R_{a_{4,l}}"] & \infty \arrow[r, bend left, "R_{a_{4,l}}"] & 0 \arrow[l, bend left, "R_{a_{4,l}}"],
\end{tikzcd}
\end{center}
\begin{center}
\begin{tikzpicture}
    \node (a) at (0,0) {$R_{a_{4,l}}(-1)$};
    \node (b) at (3,0) {$-2$};
    \node (c) at (5,0) {$0$};
    \node (d) at (9,0) {$R_{a_{4,l}}^2(-1)$};
    \node (e) at (12,0) {$\infty$};
    
    \draw[->] (a) to node [above] {$R_{a_{4,l}}^2$} (b);
    \draw[->] (b) to node [above] {$R_{a_{4,l}}^2$} (c);
    \draw[->] (c) to[in=45, out=-45, looseness=10] node{\ \ \ \ \ \ \  $R_{a_{4,l}}^2,$} (c);
    \draw[->] (d) to node [above] {$R_{a_{4,l}}^2$} (e);
    \draw[->] (e) to[in=45, out=-45, looseness=10] node{\ \ \ \ \ \ \  $R_{a_{4,l}}^2,$} (e);
\end{tikzpicture}
\end{center}
\begin{center}
\begin{tikzpicture}
    \node (a) at (0:1) {$\infty$};
    \node (b) at (0:3) {$0$};
    \node (c) at (0:5) {$R_{a_{4,l}}^2(-1).$};
    \node (d) at (145:2) {$R_{a_{4,l}}(-1)$};
    \node (e) at (210:1.5) {$-2$};
    
    \draw[->] (a) to [above, bend left] node  {$R_{a_{4,l}}^3$} (b);
    \draw[->] (b) to [below, bend left] node  {$R_{a_{4,l}}^3$} (a);
    \draw[->] (c) to node [above] {$R_{a_{4,l}}^3$} (b);
    \draw[->] (d) to node [above] {$R_{a_{4,l}}^3$} (a);
    \draw[->] (e) to node [below] {$R_{a_{4,l}}^3$} (a);
\end{tikzpicture}
\end{center}
\end{lemma}
\begin{proof}
One can check that, for any $j\in\mathbb{N}_+$, the orbits of the postcritical points of $R_{a_{4,l}}^{2j}$ split into two disjoint orbits in the second map under $R_{a_{4,l}}^{2j}$. For any $j\in\mathbb{N}_+$, the postcritical points of $R_{a_{4,l}}^{2j+1}$ go exactly along the third map under $R_{a_{4,l}}^{2j+1}$.
\end{proof}

Interestingly, the three maps $R_{a_{4,1}}, R_{a_{4,2}}, R_{a_{4,3}}$ all inherit the unmatability from that of the map $R_2$, instead of $R_{3/2}$.

\begin{pro}
$R_{a_{4,1}}, R_{a_{4,2}}, R_{a_{4,3}}$ are all atomic rational maps.
\end{pro}

\begin{figure}[h]
\centering
\includegraphics[scale=1]{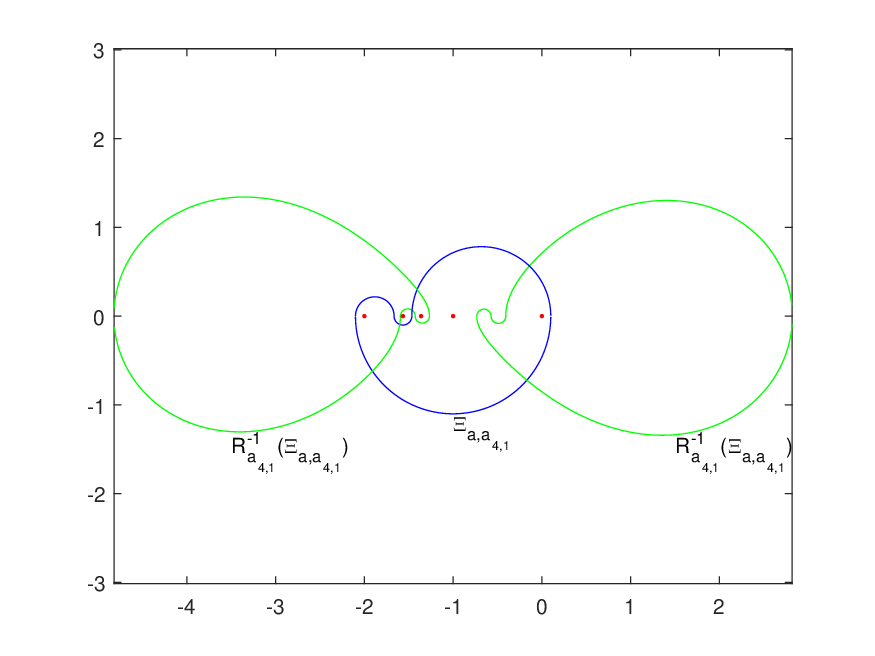}
\caption{$\Xi_{a,a_{4,1}}$ and $R_{a_{4,1}}^{-1}(\Xi_{a,a_{4,1}})$}
\label{fig5}
\end{figure}

\begin{figure}[h]
\centering
\includegraphics[scale=1]{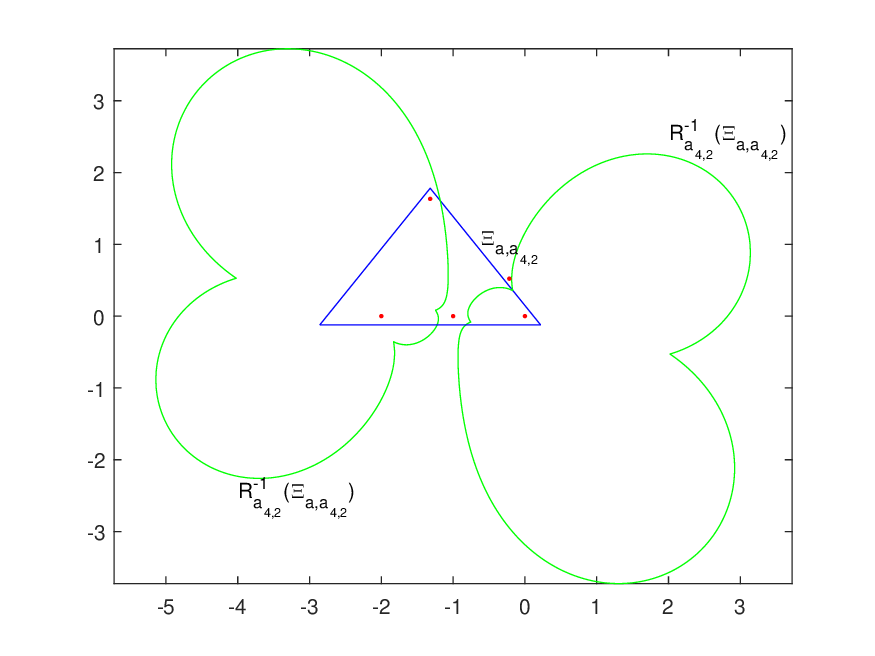}
\caption{$\Xi_{a,a_{4,2}}$ and $R_{a_{4,2}}^{-1}(\Xi_{a,a_{4,2}})$}
\label{fig6}
\end{figure}

\begin{proof}
To look for an equator or an OR equator of $R_{a_{4,1}}^j$ for $j\geq 1$, we consider partitions of $P(R_{a_{4,1}}^j)$ satisfying the immune or swapping property. We list all the non-travail partitions of $P(R_{a_{4,1}}^j)$ into two subsets in Table \ref{tab3}.
\begin{table}[ht]
\caption{Partitions of $P(R_{a_{4,1}}^j)$} 
\centering 
\begin{tabular}{c c c} 
\hline 
Type &  $P_w(R_{a_{4,1}}^j)$ & $P_b(R_{a_{4,1}}^j)$ \\
\hline
I     &   $-1.3611$    &    $-1.5652, -2, \infty, 0$\\

II   &    $-1.5652$    &    $-1.3611, -2, \infty, 0$\\

III   &    $-2$    &    $-1.3611, -1.5652, \infty, 0$\\

IV     &   $\infty$    &    $-1.3611, -1.5652, -2, 0$\\

V     &   $0$    &    $-1.3611, -1.5652, -2, \infty$\\

VI   &    $-1.3611, -1.5652,$    &    $-2,  \infty, 0$\\

VII   &    $-1.3611, -2$    &    $-1.5652, \infty, 0$\\

VIII   &    $-1.3611, \infty$    &    $-1.5652, -2, 0$\\

IX   &    $-1.3611, 0$    &    $-1.5652, -2, \infty$\\

X   &    $-1.3611, -1.5652, -2$    &    $\infty, 0$\\

XI   &    $-1.3611, -1.5652, \infty$    &    $-2, 0$\\

XII   &    $-1.3611, -1.5652, 0$    &    $-2, \infty$\\

XIII   &    $-1.3611, -2, \infty$    &    $-1.5652, 0$\\

XIV   &    $-1.3611, -2, 0$    &    $-1.5652, \infty$\\

XV   &    $-1.3611, \infty, 0$    &    $-1.5652, -2$\\
\hline 
\end{tabular}
\label{tab3} 
\end{table}

Considering Lemma \ref{lem17}, only Type XIV partition satisfies the swapping property under $R_{a_{4,1}}$ (no partition satisfies the immune property under $R_{a_{4,1}}^j$ for any $j\geq 3$ besides the Type XIV partition). So the Type XIV partition satisfies the immune property under $R_{a_{4,1}}^2$. Then an equator of  $R_{a_{4,1}}^2$ must be circular to  $\{-1.3611, -2, 0\}$  and  $\{-1.5652, \infty\}$. Unfortunately, any Jordan curve circular to  $\{-1.3611, -2, 0\}$  and  $\{-1.5652, \infty\}$ splits into two connected components under $R_{a_{4,1}}^{-1}$ (refer to Figure \ref{fig5}), which prevents it from being escalated into an equator of $R_{a_{4,1}}^2$. By further checking one can see that the Type XIV partition satisfies the immune property under $R_{a_{4,1}}^{2j}$ for any $j\geq 2$, moreover, it is the only partition with this property. Since any Jordan curve circular to  $\{-1.3611, -2, 0\}$  and  $\{-1.5652, \infty\}$ splits into at least two connected components under $R_{a_{4,1}}^{-1}$, an equator of $R_{a_{4,1}}^j$ is also impossible for any $j\geq 3$. 

Following similar lines above, we can show that any Jordan curve circular to  $\{-1.3194+1.6332i, -2, 0\}$  and  $\{-0.2174 + 0.5217i, \infty\}$ splits into at least two connected components under backward iterations of $R_{a_{4,2}}$, so an equator of $R_{a_{4,2}}^j$ is also impossible for any $j\geq 2$ (refer to Figure \ref{fig6}). $R_{a_{4,3}}$ is an atomic rational map since its dynamics is symmetric to the dynamics of $R_{a_{4,2}}$.
\end{proof}

The Jordan curve 
\begin{center}
$
\begin{array}{ll}
\Xi_{a,a_{4,1}}=& \{-1.5652+0.1e^{2\pi it}\}_{0.5\leq t\leq 1}\cup\{-1.8826+0.2174e^{2\pi it}\}_{0\leq t\leq 0.5}\vspace{3mm}\\
& \cup\{-1+1.1e^{2\pi it}\}_{0.5\leq t\leq 1}\cup\{-0.6826+0.7826e^{2\pi it}\}_{0\leq t\leq 0.5}
\end{array}
$
\end{center}
inducing the Type XIV partition in Table \ref{tab3} of $P(R_{a_{4,1}})$ and its pre-image under $R_{a_{4,1}}$ are depicted in Figure \ref{fig5}. The critical point $-1\in U_w^0$.

The triangle 
\begin{center}
$
\begin{array}{ll}
\Xi_{a,a_{4,2}}=& \{-1.2378(t-0.12)+ti\}_{-1.3194\leq t\leq 0.2169}\vspace{3mm}\\
& \cup\{1.2378(t+2.7588)+ti\}_{-2.8557\leq t\leq -1.3194}\cup\{-0.12+ti\}_{-2.8557\leq t\leq 0.2169}
\end{array}
$
\end{center}
circular to  
\begin{center}
$P_w(R_{a_{4,2}})=\{-1.3194+1.6332i, -2, 0\}$  and  $P_b(R_{a_{4,2}})=\{-0.2174 + 0.5217i, \infty\}$
\end{center}
as well as its pre-image under $R_{a_{4,2}}$ are depicted in Figure \ref{fig6}. The critical point $-1\in U_w^0$.

As the final special case, we turn to the unmatability of the postcritically finite captures of generation $5$ in the family $\Big\{R_a=\cfrac{a}{z^2+2z}\Big\}_{a\in \mathbb{C}\setminus\{0\}}$. There are also finitely many of them.
\begin{pro}
There are totally five postcritically finite generation-$5$ captures in the family $$\Big\{R_a=\cfrac{a}{z^2+2z}\Big\}_{a\in \mathbb{C}\setminus\{0\}},$$
which are parametrized by
\begin{center}
$a_{5,1}=1.29982357191516455241651733\cdots$,
\end{center}
\begin{center}
$a_{5,2}=1.650533365200074732301663490\cdots-0.4235414469091221968931574311\cdots i$,
\end{center}
\begin{center}
$a_{5,3}=1.650533365200074732301663490\cdots+0.4235414469091221968931574311\cdots i$,
\end{center}
\begin{center}
$a_{5,4}=1.69955484884234299149007785\cdots-1.509347665604408818352587205\cdots i$,
\end{center}
and
\begin{center}
$a_{5,5}=1.69955484884234299149007785\cdots+1.509347665604408818352587205\cdots i$ 
\end{center}
respectively. The forward orbits  of the postcritical points of $R_{a_{5,1}}$ are approximately
\begin{center}
\begin{tikzcd}
-1 \arrow[r,  "R_{a_{5,1}}"] & -1.2998 \arrow[r, "R_{a_{5,1}}"] & -1.4282 \arrow[r, "R_{a_{5,1}}"] & -1.5917 \arrow[r, "R_{a_{5,1}}"] & -2 \arrow[r, "R_{a_{5,1}}"] & \infty \arrow[r, bend left, "R_{a_{5,1}}"] & 0 \arrow[l, bend left, "R_{a_{5,1}}"].
\end{tikzcd}
\end{center}
\end{pro}
\begin{proof}
One can get all the postcritically finite captures of generation $5$ by solving the equation
$$R_a^4(-1)=\cfrac{(2a-2)^2}{a(a-2)^4+2(a-2)^2(2a-3)}=-2$$
with respect to $a$. 
\end{proof}

Interestingly but maybe not surprisingly this time, the maps $\{R_{a_{5,l}}\}_{l=1}^5$ all inherit the unmatability from $R_{3/2}$.

\begin{pro}\label{pro2}
Maps in $\{R_{a_{5,l}}\}_{l=1}^5$ are all primitive $2$-matings.
\end{pro}

\begin{figure}[h]
\centering
\includegraphics[scale=1.15]{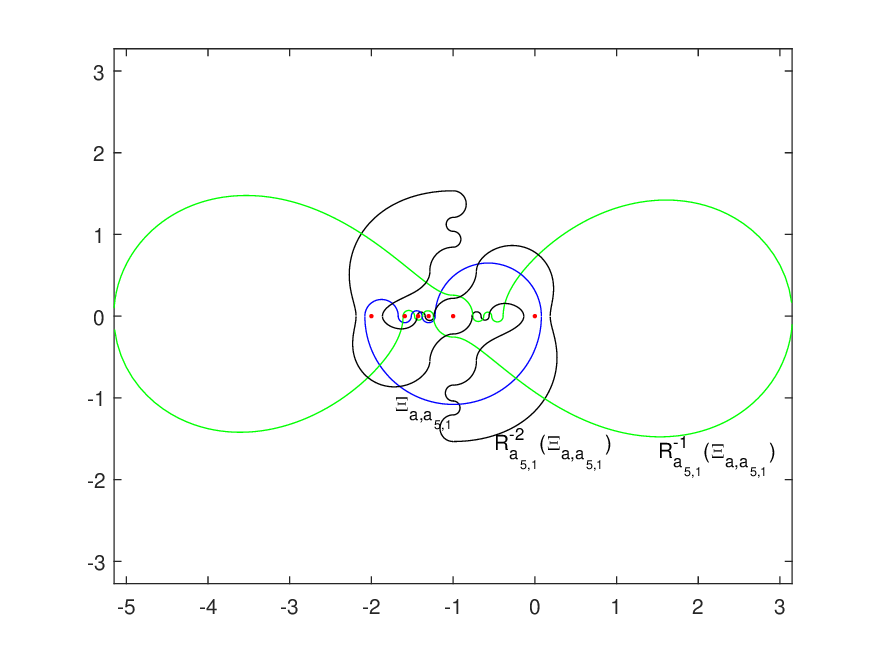}
\caption{$\Xi_{a,a_{5,1}}$, $R_{a_{5,1}}^{-1}(\Xi_{a,a_{5,1}})$and $R_{a_{5,1}}^{-2}(\Xi_{a,a_{5,1}})$}
\label{fig7}
\end{figure}

\begin{proof}
We only justify the conclusion for $R_{a_{5,1}}$, the remaining cases are similar. First note that 
$$P(R_{a_{5,1}})=P(R_{a_{5,1}}^j)=\{-1.2998, -1.4282, -1.5917, -2, \infty, 0\}$$
for any $j\in\mathbb{N}_+$. Of course $R_{a_{5,1}}$ does not admit any equator. Let
\begin{center}
$P_w(R_{a_{5,1}}^j)=\{-1.4282, -2, 0\}$  and  $P_b(R_{a_{5,1}}^j)=\{-1.2998, -1.5917, \infty\}$.
\end{center}
Further arguments show that an OR equator of $R_{a_{5,1}}$ must be circular to both $P_w(R_{a_{5,1}}^j)$ and $P_b(R_{a_{5,1}}^j)$ for any $j\geq 1$. Moreover, an equator of $R_{a_{5,1}}^j$ also must be circular to both $P_w(R_{a_{5,1}}^j)$ and $P_b(R_{a_{5,1}}^j)$ for any $j\geq 2$. Now let $\Xi$ be a Jordan curve circular to both $P_w(R_{a_{5,1}}^j)$ and $P_b(R_{a_{5,1}}^j)$. Let $U_w^0$ be the closure of the connected component of $\mathbb{C}\setminus \Xi$ containing $\{-1.4282, -2,  0\}$ and let $U_b^0$ be the closure of the connected component of $\mathbb{C}\setminus \Xi$ containing $\{-1.2998, -1.5917, \infty\}$. $R_{a_{5,1}}^{-j}(\Xi)$ is always connected for any integer $j\geq 1$. Further analysis shows that the two postcritical points $-1.2998$ and $-1.5917$ must be confined in the finite domains enclosed by $U_w^0$ and $R_{a_{5,1}}^{-1}(U_b^0)=U_b^1$, which prevents $\Xi$ to be isotopic to $R_{a_{5,1}}^{-1}(\Xi)$ \emph{rel.} $P(R_{a_{5,1}})$.  However, they are released from the finite confined domains enclosed by $U_w^0$ and $R_{a_{5,1}}^{-2}(U_w^0)=U_w^2$, so $\Xi$ is isotopic to $R_{a_{5,1}}^{-2}(\Xi)$ \emph{rel.} $P(R_{a_{5,1}})$ (refer to Figure \ref{fig7}). 
\end{proof}

Let $r=0.08$. The Jordan curve 
\begin{center}
$
\begin{array}{ll}
\Xi_{a,a_{5,1}}=& \{-1.5917+re^{2\pi it}\}_{0.5\leq t\leq 1}\cup\{-1.79585-r+0.20415e^{2\pi it}\}_{0\leq t\leq 0.5}\vspace{3mm}\\
& \cup\{-1+(1+r)e^{2\pi it}\}_{0.5\leq t\leq 1}\cup\{-0.6499+r+0.6499e^{2\pi it}\}_{0\leq t\leq 0.5}\vspace{3mm}\\
& \cup\{-1.2998+re^{2\pi it}\}_{0.5\leq t\leq 1}\cup\{-1.44575+(0.14595-r)e^{2\pi it}\}_{0\leq t\leq 0.5}\vspace{3mm}\\
\end{array}
$
\end{center}
circular to $P_w(R_{a_{5,1}}^j)=\{-1.4282, -2, 0\}$  and  $P_b(R_{a_{5,1}}^j)=\{-1.2998, -1.5917, \infty\}$ as well as its pre-images under $R_{a_{5,1}}$ are depicted in Figure \ref{fig7}. 

\subsection{The proof of Theorem \ref{thm7}}

Along with explorations on the unmatability of all the captures of generation $2, 3, 4, 5$, now it is more convenient for us to present the proof of Theorem   \ref{thm7}.
\begin{proof}[Proof of Theorem   \ref{thm7}]
Assume $R_{a_*}$ is a postcritically finite capture of generation $k$ in the family $\Big\{R_a=\cfrac{a}{z^2+2z}\Big\}_{a\in \mathbb{C}\setminus\{0\}}$ such that its critical points go as 
\begin{center}
\begin{tikzcd}
p_1 \arrow[r, "R_{a_*}"] & p_2 \arrow[r, "R_{a_*}"] & \cdots \arrow[r, "R_{a_*}"] & p_k \arrow[r, "R_{a_*}"]  & \infty \arrow[r, bend left, "R_{a_*}"] & 0 \arrow[l, bend left, "R_{a_*}"], 
\end{tikzcd}
\end{center}
in which $p_1=-1, p_k=-2$. Thus $P(R_{a_*})=\{p_2, \cdots, p_k, \infty, 0\}$.  

We first deal with the case $a_*\in\mathbb{R}$ (there is one and only one postcritically finite generation-$k$ capture parametrized by a pure real parameter for every $k\geq 2$). Note that in this case the postcritical points are all reals such that
$$P(R_{a_*})\subset [-2, 0]$$
and we have
$$-1>p_2>p_3 \cdots > p_{k-1} >-2.$$ 
$R_{a_*}$ preserves the real line and swaps the upper and lower half complex plane. Now we consider non-travail partitions of  $P(R_{a_*})$ into two subsets satisfying the immune or swapping property. One can easily figure out now that the only reasonable partition is
$$P_b(R_{a_*})=\{\cdots, p_{k-5}, p_{k-3}, p_{k-1}, \infty\},\ \ P_w(R_{a_*})=\{\cdots, p_{k-4}, p_{k-2}, p_k, 0\}.$$ 
This partition satisfies the swapping property for $R_{a_*}$ while it satisfies the immune property for $R_{a_*}^2$. Let $\Xi$ be a Jordan curve circular to both $P_b(R_{a_*})$ and $P_w(R_{a_*})$. Let $U_b^0$ be the closure of the connected component of $\mathbb{C}\setminus \Xi$ containing $P_b(R_{a_*})$ and let $U_w^0$ be the closure of the connected component of $\mathbb{C}\setminus \Xi$ containing $P_w(R_{a_*})$.  We claim that $R_{a_*}^{-1}(\Xi)$ splits for even $k$, which is enough to force $R_{a_*}$ to be an atomic rational map in this case.   
To see this, note that $R_{a_*}^{-1}(\Xi)$ is symmetric with respect to both the real line and the critical point $-1$. So the connectedness depends on whether $-1\in U_w^1=R_{a_*}^{-1}(U_w^0)$. Now if $k$ is even, the partition is
$$P_b(R_{a_*})=\{p_3, \cdots, p_{k-5}, p_{k-3}, p_{k-1}, \infty\},\ \ P_w(R_{a_*})=\{p_2, \cdots, p_{k-4}, p_{k-2}, p_k, 0\}.$$ 
Obviously we have $-1\in U_w^1$ now, which forces $R_{a_*}^{-1}(U_b^0)=U_b^1$ and hence $R_{a_*}^{-1}(\Xi)$ to split. 

Now if $k$ is odd, the partition is
$$P_b(R_{a_*})=\{p_2, \cdots, p_{k-5}, p_{k-3}, p_{k-1}, \infty\},\ \ P_w(R_{a_*})=\{p_3, \cdots, p_{k-4}, p_{k-2}, p_k, 0\}.$$ 
Obviously we have $-1\notin U_w^1$ now. This guarantees $R_{a_*}^{-1}(U_b^0)=U_b^1$ and hence $R_{a_*}^{-1}(\Xi)$ to be connected. However, since $R_{a_*}$ swaps the upper and lower half plane, the finite members in $P_b(R_{a_*})$ (these are, $p_2, \cdots, p_{k-5}, p_{k-3}, p_{k-1}$) will be confined in the finite domains enclosed by $U_w^0$ and $U_b^1$. This prevents $\Xi$ to be isotopic to $R_{a_*}^{-1}(\Xi)$ \emph{rel.} $P(R_{a_*})$, which denies the possibility for $\Xi$ to be escalated to an OR equator of $R_{a_*}$. So $R_{a_*}$ can not be an OR mating in this case. However, since $R_{a_*}$ swaps the upper and lower half plane,  the finite members in $P_b(R_{a_*})$ are released from the confined domains by $U_w^0$ and $U_w^2=R_{a_*}^{-1}(U_w^1)$ under another backward iteration of $R_{a_*}$. This means that $\Xi$ is isotopic to $R_{a_*}^{-2}(\Xi)$ \emph{rel.} $P(R_{a_*})$. The isotopy has to be orientation-preserving since the partition $P_b(R_{a_*}), P_w(R_{a_*})$ satisfies the immune property under $R_{a_*}^2$.  So $\Xi$ evolves into an equator of $R_{a_*}^2$ in this case, which guarantees that $R_{a_*}^2$ is a $1$-mating.    

Recall from \cite[Section 4]{AY} that the spine of $R_{a_*}$ is the combination of the axes of two bubble rays containing all the critical and postcritical points.  In the remaining cases that $a_*\notin\mathbb{R}$, the spines take the role of the real line, and all the arguments above apply to these cases.   
\end{proof}

\section{The postcritical realization programme of atoms and matings}\label{sec7}

\subsection{The postcritical realization programme for postcritical sets without prescribed maps via the mating and atomic family}

The first level of the postcritical realization programme focuses on \emph{realizing} a (finite) subset of the Riemann sphere by some (postcritically finite) rational map.  

\begin{defn}\label{def6}
For a set $X\subset \mathbb{P}^1(\mathbb{C})$, we say that $X$ is \emph{realized} by a rational map $R: \mathbb{P}^1(\mathbb{C})\rightarrow \mathbb{P}^1(\mathbb{C})$ if $P(R)=X$.
\end{defn}
In this subsection we mean to prove Theorem \ref{thm9}. The proof essentially relies on DeMarco-Koch-McMullen's techniques on realizing the postcritical sets of rational maps in \cite{DKM}. We achieve the goal by two steps.

\begin{proof}[Proof of Theorem \ref{thm9} (A)]
If $\infty\in X$, by \cite[Theorem 2.1]{DKM}, there exists a polynomial $P_*: \mathbb{P}^1(\mathbb{C})\rightarrow\mathbb{P}^1(\mathbb{C})$, such that the postcritical set $P(P_*)=X$. $P_*$ is of course a $1$-mating and rigid. If $\infty\notin X$, we choose arbitrary $x\in X$ and a non-degenerate M\"obius transformation 
$$f(z)=\cfrac{az+b}{cz+d}$$
satisfying $-\cfrac{d}{c}=x$. Now according to \cite[Proof of Theorem 1.1]{DKM}, there exists some polynomial  $P_*$ such that $P(f^{-1}\circ P_*\circ f)=X$. Note that the unmatability of rational maps is preserved under M\"obius conjugations, so the postcritically finite rational map $f^{-1}\circ P_*\circ f$ is still a $1$-mating. It is rigid according to \cite{DH1}.  
\end{proof}

It turns out that the rational maps employed by DeMarco-Koch-McMullen on realizing given postcritical sets are all $1$-matings. An interesting question is that, can these sets be the postcritical sets of some atomic rational maps? We only have answers in case  $\#X=2$ and $3$. 

\begin{proof}[Proof of Theorem \ref{thm9} (B)]
Up to M\"obius conjugations, if a rational map $R$ admits only two postcritical points, then   
\begin{center}
$R(z)=z^k$ or $R(z)=z^{-k}$
\end{center}
for some $k\geq 2$. The former one is a $1$-mating while the latter one is an OR mating. Thus no set $X$ with two points can be the postcritical set of an atomic rational map. In case $\#X=3$, without loss of generality we can assume $X=\{0,1, \infty\}$. The postcritical set of the map $A(z)=1-1/z^2$ in \cite[Fig.3]{DKM} coincides with $X$. It is an atomic rational map.  To see this, note that $C(A)=\{0,\infty\}$ and the critical points go as 
\begin{center}
\begin{tikzcd}
\infty \arrow[r, "A"] & 1 \arrow[r, "A"] & 0 \arrow[ll, bend left, "A"],
\end{tikzcd}
\end{center}
which shows that it is a hyperbolic map. Then $A^j(z)$ (obviously not a Thurston polynomial) does not admit any equator for any $j\geq 1$ in virtue of \cite[Proposition 4.3]{Mey1}. 
\end{proof}

\subsection{The postcritical realization programme for maps via the mating and atomic family}

The further level of the postcritical realization programme focuses on \emph{realizing} a (finite) subset of the Riemann sphere and maps on it simultaneously by some (postcritically finite) rational map. 

\begin{defn}\label{def7}
For a set $X\subset \mathbb{P}^1(\mathbb{C})$ and a map $R_X: X\rightarrow X$ on it, we say that $R_X$ can be \emph{realized} by a rational map  if there exists some rational map $R$ with $P(R)=X$ and the restriction $R|_{X}=R_X$.
\end{defn}

So a map realizing a map $R_X$ must also realize the set $X$.  We are particularly interested in realizing maps on finite sets in $\mathbb{P}^1(\bar{\mathbb{Q}})$, that is, \cite[Question 1.2]{DKM}, via matings and atoms respectively. The setting of the question on $\mathbb{P}^1(\bar{\mathbb{Q}})$ is critical. In case $\#X=3$, there are totally 7 maps on $X$ (in fact 27 types without modulo by M\"obius maps), see \cite[Fig.3. A-G]{DKM}.   DeMarco, Koch and McMullen have realized every map in $\{A,\cdots,G: X\rightarrow X\}$ by an explicit rigid rational map. We strengthen their result by showing that we can achieve their aim via both the mating and atomic family individually, except for the map 
\begin{center}
A:
\begin{tikzcd}
* \arrow[r] & * \arrow[r] & * \arrow[ll, bend left].
\end{tikzcd}
\end{center} 
We also decide the unmatability of every map they construct to realize the prescribed maps on $X$ with  $\#X=3$. 

\begin{theorem}\label{thm10}
Let $X\subset\mathbb{P}^1(\bar{\mathbb{Q}})$ be a finite set satisfying $\#X=3$. Then every map $B,\cdots,G: X\rightarrow X$ can be realized by rigid rational maps arising as $n$-matings and atomic maps respectively.
\end{theorem}

\begin{proof}
In case $\#X=3$, without loss of generality we can assume $X=\{0,1, \infty\}$. DeMarco, Koch and McMullen realize every map B-G by an explicit rational map which is either an $n$-mating or an atomic map. We will supplement either an atomic or an $n$-mating which realizes every map B-G. One can use M\"obius conjugations to adjust the maps on $\{0,1, \infty\}$ for all the maps.

\begin{itemize}

\item map B. The map 
$$B(z)=\cfrac{(z-\alpha)^3}{(z-1+\alpha)^3}$$
with $\alpha^2-\alpha+1=0$ is an $n$-mating for some integer $n\in\mathbb{N}_+$. To see this, note that $C(B(z))=\{\alpha, 1-\alpha\}$. The critical points go as
\begin{center}
\begin{tikzpicture}
    \node (a) at (0:1) {$1$};
    \node (d) at (145:1) {$0$}; 
    \node [left of = d, node distance=0.7in] (b) {$\alpha$};
    \node (e) at (210:1) {$\infty$};
    \node [left of = e, node distance=0.7in] (c) {$1-\alpha$};

    \draw[->] (d) to [above] node  {$B$} (a);
    \draw[->] (b) to [above] node  {$B$} (d);
    \draw[->] (c) to [above] node  {$B$} (e);
    \draw[->] (e) to node [below] {$B$} (a);
    \draw[->] (a) to[in=45, out=-45, looseness=10] node{\ \ \ \ \ \ \  $B.$} (a);
\end{tikzpicture}
\end{center}  
One can see that both critical points of $B$ are strictly pre-periodic, so $\mathcal{J}(B)=\mathbb{P}^1(\mathbb{C})$. Applying \cite[Theorem 1.1]{Mey2} we can see that $B(z)$ is an $n$-mating for some $n\in\mathbb{N}_+$. In the following we provide a hyperbolic atomic rational map also realizing map B. Let 
\begin{center}
$\alpha_1^B=\cfrac{8+4\sqrt{5}i}{9}$  and  $\alpha_2^B=\cfrac{8-4\sqrt{5}i}{9}$.
\end{center}
Let $c_j^B=(4-3\alpha_j^B)/4$ for $j=1,2$. The map 
$$B_{a,\alpha_j^B}(z)=\cfrac{\gamma_j^B(z-\alpha_j^B)^5}{z(z-1)^2}$$ 
with $\gamma_j^B=\cfrac{c_j^B(c_j^B-1)^2}{(c_j^B-\alpha_j^B)^5}$ is  an atomic map for either $j\in\{1,2\}$. In fact one can check that $C(B_{a,\alpha_j^B}(z))=\{1, c_j^B, \alpha_j^B, \infty\}$ for either $j\in\{1,2\}$. The critical points go as
\begin{center}
\begin{tikzpicture}
    \node (a) at (0:1) {$\infty$};
    \node (d) at (145:1) {$1$}; 
    \node [left of = d, node distance=0.7in] (b) {$c_j^B$};
    \node (e) at (210:1) {$0$};
    \node [left of = e, node distance=0.7in] (f) {$\alpha^B_j$};
    
    \draw[->] (d) to [above] node  {$B_{a,\alpha_j^B}$} (a);
    \draw[->] (b) to [above] node  {$B_{a,\alpha_j^B}$} (d);
    \draw[->] (f) to [below] node  {$B_{a,\alpha_j^B}$} (e);
    \draw[->] (e) to node [below] {$B_{a,\alpha_j^B}$} (a);
    \draw[->] (a) to[in=45, out=-45, looseness=10] node{\ \ \ \ \ \ \  $B_{a,\alpha_j^B}.$} (a);
\end{tikzpicture}
\end{center} 
A similar argument as the proof of Proposition \ref{pro7} justifies $B_{a,\alpha_j^B}(z)$ is atomic for 
any $j=1,2$.

\item map C. The map $C(z)=z^2(3-2z)$ is a $1$-mating as a polynomial. The map 
$$C_a(z)=E_a^2(z)$$
is a hyperbolic atomic map also realizing map C (see $E_a(z)$ in map E). Details are left to the readers. 

\item map D. The map $D(z)=(1-\cfrac{2}{z})^2$ is an $n$-mating for some integer $n\in\mathbb{N}_+$. The critical points of $D$ go as
\begin{center}
\begin{tikzpicture}
    \node (a) at (0,0) {$2$};
    \node (b) at (2,0) {$0$};
    \node (c) at (4,0) {$\infty$};
    \node (d) at (6,0) {$1$};

    \draw[->] (a) to node [above] {$D$} (b);
    \draw[->] (b) to node [above] {$D$} (c);
    \draw[->] (c) to node [above] {$D$} (d);
    \draw[->] (d) to[in=45, out=-45, looseness=10] node{\ \ \ \ $D.$} (d);
\end{tikzpicture}
\end{center}
One can see that all critical points are strictly pre-periodic, which means that it is an $n$-mating in virtue of \cite[Theorem 1.1]{Mey2}. The map
$$D_a(z)=-\cfrac{(z+4)(z-1)^2}{9z}$$
is an atomic rational map which also realizes map D. In fact one can check that  $C(D_a(z))=\{1,\infty,-2\}$. The critical points of $D_a(z)$ go as 
\begin{center}
\begin{tikzpicture}
    \node (a) at (0,0) {$-2$};
    \node (b) at (2,0) {$1$};
    \node (c) at (4,0) {$0$};
    \node (d) at (6,0) {$\infty$};

    \draw[->] (a) to node [above] {$D_a$} (b);
    \draw[->] (b) to node [above] {$D_a$} (c);
    \draw[->] (c) to node [above] {$D_a$} (d);
    \draw[->] (d) to[in=45, out=-45, looseness=10] node{\ \ \ \ $D_a.$} (d);
\end{tikzpicture}
\end{center}

\item map E. The map $E(z)=-z^2+1$ is a $1$-mating as a polynomial. The map 
$$E_a(z)=\cfrac{(z-1)^3}{3z-1}$$
is a hyperbolic atomic map also realizing map E. The critical points $\{0,1,\infty\}$ of $E_a(z)$ go as 
\begin{center}
\begin{tikzpicture}
    \node (a) at (0,0) {$\infty$};
    \node (b) at (4,0) {$0$};
    \node (c) at (6,0) {$1.$};
    
    \draw[->] (b) to [above, bend left] node  {$E_a$} (c);
    \draw[->] (c) to [below, bend left] node  {$E_a$} (b);
    \draw[->] (a) to[in=45, out=-45, looseness=10] node{\ \ \ \ $E_a,$} (a);
\end{tikzpicture}
\end{center}

\item map F. The map $F(z)=\cfrac{(2z-1)^2}{4z(z-1)}$ is a hyperbolic atomic map. One can check that $C(F)=\{1/2, \infty\}$ while
\begin{center}
\begin{tikzcd}
1/2 \arrow[r, "F"] & 0 \arrow[r, "F"] & \infty \arrow[r, bend left, "F"] & 1 \arrow[l, bend left, "F"].
\end{tikzcd}
\end{center}

Let $\{\alpha_j^F\}_{j=1,2}$ be the two roots of the polynomial
\begin{center}
$4\alpha^4-19\alpha^3+9\alpha^2+27\alpha+27$
\end{center}
other than $3$ ($3$ is a double root). Let $c_j^F=(\alpha_j^F-3)/2$ for $j=1,2$. Then this type of map can be realized by the non-hyperbolic $n$-mating
\begin{center}
$F_{m,\alpha_j^F}=\cfrac{(\alpha_j^F)^3(z-1)}{(z-\alpha_j^F)^3}$
\end{center}
for either $j=1,2$. In fact one can check that $C(F_{m,\alpha_j^F})=\{\alpha_j^F, c_j^F, \infty\}$ and the critical points go as   
\begin{center}
\begin{tikzcd}
\alpha_j^F \arrow[r, "F_{m,\alpha_j^F}"] & \infty \arrow[r, "F_{m,\alpha_j^F}"] & 0 \arrow[r, bend left, "F_{m,\alpha_j^F}"] & 1 \arrow[l, bend left, "F_{m,\alpha_j^F}"] & c_j^F \arrow[l, "F_{m,\alpha_j^F}"]
\end{tikzcd}
\end{center}
for either $j=1,2$. All the critical points are strictly pre-periodic.

\item map G. The map $G(z)=(2z-1)^2$ is a $1$-mating as a polynomial. The map 
$$G_a(z)=-\cfrac{4(z^2-z)^2}{(2z-1)^2}$$
is a hyperbolic atomic map also realizing map G. To see this, we recall that the generation-$2$ capture $R_2=\cfrac{2}{z^2+2z}$ in the slice of $V_2$ is an atomic rational map. One can check that $P(R_2^2)=\{-2, \infty, 0\}$ while
\begin{center}
\begin{tikzpicture}
    \node (a) at (0,0) {$\infty$};
    \node (b) at (4,0) {$-2$};
    \node (c) at (6,0) {$0$};
    
    \draw[->] (b) to [above] node  {$R_2^2$} (c);
    \draw[->] (c) to[in=45, out=-45, looseness=10] node{\ \ \ \ $R_2^2.$} (c);
    \draw[->] (a) to[in=45, out=-45, looseness=10] node{\ \ \ \ $R_2^2,$} (a);
\end{tikzpicture}
\end{center}
Essentially  $R_2^2$ has realized map G as an atomic map. To be more precise, let us take $f(z)=-2z$, then we get $G_a(z)=f^{-1}\circ R_2^2\circ f(z)$.

\end{itemize}

All the postcritically finite rational maps above are rigid since they all admit only three postcritical points.
\end{proof}

Note that we have decided the unmatability of every map constructed by DeMarco-Koch-McMullen to realize the prescribed maps on $X$ with  $\#X=3$ in the Proof of Theorem \ref{thm9} (B) and Theorem \ref{thm10}  in fact.

\begin{coro}
The maps $B(z), C(z), D(z), E(z), G(z)$ are all $n$-matings for some integer $n\in\mathbb{N}_+$, while the maps $A(z), F(z)$ are both atomic maps. 
\end{coro}





We do not know whether map A can be realized by some non-hyperbolic $n$-matings (it certainly can not be realized by a hyperbolic $n$-mating). It would be interesting to try to decide the fold of the non-hyperbolic $n$-matings in $\{B(z), D(z), F_{m,\alpha_1^F}(z), F_{m,\alpha_2^F}(z)\}$ which realize certain maps in the Proof of Theorem \ref{thm10}. Although most of the maps on $X$ with $\#X=3$ can be realized by some rational maps arising as $n$-matings, there is some distinction in fact. The maps $C, E, G$ are realized by polynomials (which are hyperbolic $1$-matings), while the maps $ B, D, F$ are realized by non-hyperbolic $n$-matings in our results. Although we are not sure whether we can take the fold $n=1$, the maps $A, B, D, F$ can not be realized by polynomials (in fact these can not be realized by hyperbolic $n$-matings). This demonstrates some characteristic of certain maps (maps on $X$), which may deserve special attention in the postcritical realization programme.

\section{The postcritical realization programme for maps on sets of four points}

We mean to establish algebraic structure on collections of maps on finite sets in $\mathbb{P}^1(\mathbb{C})$ in this section. For simplicity of presentation the emphasis is laid on sets of four points, except the first subsection. The structure allows us to simplify the postcritical realization programme for maps on finite sets.

\subsection{The compositive trick}

We have implicitly employed the trick in constructing $C_a(z)$ and $G_a(z)$ in the proof of Theorem \ref{thm10} in fact. 

\begin{defn}\label{def10}
We say a map $R: \mathbb{P}^1(\mathbb{C})\rightarrow\mathbb{P}^1(\mathbb{C})$ preserves a set $X\subset\mathbb{P}^1(\mathbb{C})$ if $R(X)=X$, it sub-preserves a set $X\subset\mathbb{P}^1(\mathbb{C})$ if $R(X)\subset X$, it sup-preserves a set $X\subset\mathbb{P}^1(\mathbb{C})$ if $R(X)\supset X$.
\end{defn}

We do not require $\#X$ to be finite in Definition \ref{def10}. Of course any map can not strictly sup-preserve a finite set.  The compositive trick relies on the following simple observation.

\begin{lemma}\label{lem18}
For two rational maps $R_1, R_2: \mathbb{P}^1(\mathbb{C})\rightarrow\mathbb{P}^1(\mathbb{C})$, if $R_1$ and  $R_2$ both sub-preserve $P(R_1)\cup P(R_2)$, then
\begin{equation}\label{eq56}
P(R_1\circ R_2)\subset P(R_1)\cup P(R_2).
\end{equation} 
\end{lemma}
\begin{proof}
Note that $C(R_1\circ R_2)=C(R_2)\cup R_2^{-1}(C(R_1))$. So we have  
\begin{equation}\label{eq54}
V(R_1\circ R_2)=R_1\circ R_2(C(R_2))\cup R_1(C(R_1))\subset P(R_2)\cup P(R_1)
\end{equation} 
since $R_1(V(R_2))\subset P(R_2)\cup P(R_1)$. Again since $R_1$ and  $R_2$ sub-preserve $P(R_1)\cup P(R_2)$ we have
\begin{equation}\label{eq53}
R_1\circ R_2(V(R_1\circ R_2))\subset P(R_2)\cup P(R_1).
\end{equation}  
Inductively (\ref{eq53}) induces
\begin{equation}\label{eq55}
(R_1\circ R_2)^k(V(R_1\circ R_2))\subset P(R_2)\cup P(R_1)
\end{equation}
for any $k\geq 1$. Finally, (\ref{eq56}) follows from a combination of (\ref{eq54}) and (\ref{eq55}).
\end{proof}

\begin{rem}
Be careful that $R_1$ preserves $P(R_2)$ and  $R_2$ preserves $P(R_1)$ are not enough to guarantee
\begin{equation}\label{eq57}
P(R_1\circ R_2)= P(R_1)\cup P(R_2)
\end{equation}
in Lemma \ref{lem18}.  Some additional restrictions on 
\begin{center}
$R_1|_{C(R_1)\cup C(R_2)\cup P(R_1)\cup P(R_2)}$ and $R_2|_{C(R_1)\cup C(R_2)\cup P(R_1)\cup P(R_2)}$ 
\end{center}
can be imposed to achieve (\ref{eq57}). 
\end{rem}

The importance of Lemma \ref{lem18} in the postcritical realization programme is that it provides a way to build new postcritically finite rational maps by compositions of  established postcritically finite rational maps: the \emph{compositive trick}. 

\begin{coro}\label{cor4}
For two postcritically finite rational maps $R_1, R_2: \mathbb{P}^1(\mathbb{C})\rightarrow\mathbb{P}^1(\mathbb{C})$, if $R_1$ and  $R_2$ both sub-preserve $P(R_1)\cup P(R_2)$, then $R_1\circ R_2$ is a postcritically finite rational map satisfying (\ref{eq56}).
\end{coro}
\begin{proof}
The result follows from Lemma \ref{lem18} instantly.
\end{proof}

Corollary \ref{cor4} takes a more attractive form in the following case.
\begin{coro}\label{cor5}
If two rational maps $R_1, R_2: \mathbb{P}^1(\mathbb{C})\rightarrow\mathbb{P}^1(\mathbb{C})$ are both postcritically finite with 
\begin{center}
$P(R_2)\subset P(R_1)=V(R_1)$,
\end{center}
such that $R_2$ sub-preserves $P(R_1)$, then $R_1\circ R_2$ is a postcritically finite rational map with 
\begin{equation}\label{eq62}
P(R_1\circ R_2)=V(R_1\circ R_2)=P(R_1).
\end{equation}
\end{coro}
\begin{proof}
Since $P(R_1)\cup P(R_2)=P(R_1)$, both $R_1$ and $R_2$ sub-preserve $P(R_1)$. Then $R_1\circ R_2$ is a postcritically finite rational map satisfying
\begin{equation}\label{eq61}
P(R_1\circ R_2)\subset P(R_1)
\end{equation} 
according to Corollary \ref{cor4}. Now since $P(R_1\circ R_2)\supset V(R_1\circ R_2)\supset V(R_1)=P(R_1)$, considering  (\ref{eq61}), we get (\ref{eq62}). 
\end{proof}

The condition $P(R_1)=V(R_1)$ in Corollary \ref{cor5} can be relaxed to require all the initial points in the orbits of the map $R_1\circ R_2|_{P(R_1)}$ are contained in $V(R_1)\cup R_1(V(R_2))$, at the risk of losing (\ref{eq62}). One can also refer to \cite{Koc} for some method on constructing postcritically finite rational maps from the periodic points of maps on the \emph{moduli space} $\mathcal{M}_X$ with $\#X<\infty$.

There are obvious distinctions between \emph{period maps} and \emph{strictly pre-periodic maps} on sets of finite points.

\begin{defn}\label{def9}
Let $X\subset \mathbb{P}^1(\mathbb{C})$ with finitely many points.  A map on $X$ is said to be periodic if all points in $X$ are periodic under the map, it is said to be strictly pre-periodic if at least one point in $X$ is strictly pre-periodic under the map.
\end{defn} 

A periodic map on $X$ preserves $X$ while a strictly pre-periodic one only sub-preserves $X$. The composite operation is closed in the periodic family and the strictly pre-periodic family respectively, however, composition between a periodic map and a strictly pre-periodic one always results in a strictly pre-periodic map. We will investigate realizing the periodic and the strictly pre-periodic maps respectively in the next two subsections, then we combine the results together in the final subsection. Attention will be put on sets with four points in the following, while the techniques developed apply to sets with finitely many points in fact.

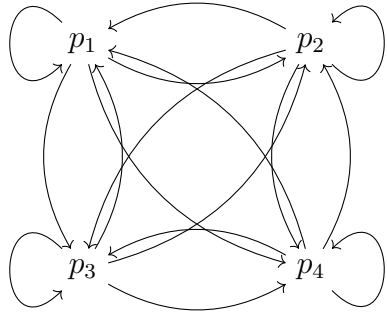
\begin{figure}
\centering
\begin{tikzpicture}
    \node (a) at (0,0) {$p_1$};
    \node (b) at (3,0) {$p_2$};
    \node (c) at (0,-3) {$p_3$};
    \node (d) at (3,-3) {$p_4$};

    \draw[->] (a) to [bend right] node {} (b);
    \draw[->] (b) to [bend right] node {} (a);
    \draw[->] (a) to [bend right] node {} (c);
    \draw[->] (c) to [bend right] node {} (a);
    \draw[->] (a) to [bend right] node {} (d);
    \draw[->] (d) to [bend right] node {} (a);
    \draw[->] (b) to [bend right] node {} (c);
    \draw[->] (c) to [bend right] node {} (b);
    \draw[->] (b) to [bend right] node {} (d);
    \draw[->] (d) to [bend right] node {} (b);
    \draw[->] (c) to [bend right] node {} (d);
    \draw[->] (d) to [bend right] node {} (c);
     \draw[->] (a) to[in=-135, out=135, looseness=6] node {} (a);
     \draw[->] (b) to[in=45, out=-45, looseness=6] node {} (b);
     \draw[->] (c) to[in=-135, out=135, looseness=6] node {} (c);
    \draw[->] (d) to[in=45, out=-45, looseness=6] node {} (d);
\end{tikzpicture}
\caption{maps on $X$ with  $\#X=4$}
\label{fig22} 
\end{figure}

Let $X=\{p_1, p_2, p_3, p_4\}\subset \mathbb{P}^1(\mathbb{C})$ such that $\#X=4$, which means that the four points are always mutually distinct from one another. There are totally $4^4=256$ types of maps on $X$ (refer to Figure \ref{fig22}, note that M\"obius conjugations can only manipulate at most three points simultaneously), which makes the realization of them a tremendous work. However, there are close relationships between these maps, exploiting which we only need to realize few of them in order to achieve them all. 

\subsection{The periodic family on sets of four points}

There are totally 24 types of periodic maps on $X=\{p_1, p_2, p_3, p_4\}$ of four points, we first give a classification of them according to the number of periodic cycles in them. There are 6 types of periodic maps on $X$ with only one periodic cycle. The "$\mathfrak{P}$"  stands for "periodic", "$|$" stands for the number of periodic cycles. There are 11 types of periodic maps on $X$ with two periodic cycles. Additional letters are employed to distinguish the two series in their names. There are 6 types of periodic maps on $X$ with three periodic cycles. One can add an additional integer to indicate distinction of the various types of maps in name of the same series, for example, let
\begin{center}
$\mathfrak{P}|1$: \begin{tikzpicture}
    \node (a) at (0,0) {$p_1$};
    \node (b) at (2,0) {$p_2$};
    \node (c) at (4,0) {$p_3$};
    \node (d) at (6,0) {$p_4$};

    \draw[->] (a) to node {} (b);
    \draw[->] (b) to node {} (c);
    \draw[->] (c) to node {} (d);
    \draw[->] (d) to [in=-10, out=190, looseness=1] node {} (a);
\end{tikzpicture},
\end{center}

\begin{center}
$\mathfrak{P}|6$: \begin{tikzpicture}
    \node (a) at (0,0) {$p_1$};
    \node (b) at (2,0) {$p_2$};
    \node (c) at (4,0) {$p_4$};
    \node (d) at (6,0) {$p_3$};

    \draw[->] (a) to node {} (b);
    \draw[->] (b) to node {} (c);
    \draw[->] (c) to node {} (d);
    \draw[->] (d) to [in=-10, out=190, looseness=1] node {} (a);
\end{tikzpicture}
\end{center}
and
\begin{center}
$\mathfrak{P}||A1$: \begin{tikzpicture}
    \node (a) at (0,0) {$p_1$};
    \node (b) at (2,0) {$p_2$};
    \node (c) at (4,0) {$p_3,$};
    \node (d) at (5,0) {$p_4$};

    \draw[->] (a) to node  {} (b);
    \draw[->] (b) to node {} (c);
    \draw[->] (c) to [in=-10, out=190, looseness=1] node {} (a);
    \draw (d) edge[loop right] node {.} (d);
\end{tikzpicture} 
\end{center}
The sole map $\mathfrak{P}||||$ on $X$ with four periodic cycles is the identity map. All the series of periodic maps and numbers of types in them on $X$ are listed in Table \ref{tab7}.

\begin{table}[ht]
\caption{Number of periodic maps} 
\centering 
\begin{tabular}{c c c} 
\hline 
Series &  map & number of types \\
\hline
$\mathfrak{P}|$     &   
\begin{tikzpicture}
    \node (a) at (0,0) {$*$};
    \node (b) at (2,0) {$*$};
    \node (c) at (4,0) {$*$};
    \node (d) at (6,0) {$*$};

    \draw[->] (a) to node {} (b);
    \draw[->] (b) to node {} (c);
    \draw[->] (c) to node {} (d);
    \draw[->] (d) to [in=-10, out=190, looseness=1] node {} (a);
\end{tikzpicture}
& 6\\

$\mathfrak{P}||A$     &   
\begin{tikzpicture}
    \node (a) at (0,0) {$*$};
    \node (b) at (2,0) {$*$};
    \node (c) at (4,0) {$*,$};
    \node (d) at (5,0) {$*$};

    \draw[->] (a) to node  {} (b);
    \draw[->] (b) to node {} (c);
    \draw[->] (c) to [in=-10, out=190, looseness=1] node {} (a);
    \draw (d) edge[loop right] node {} (d);
\end{tikzpicture}   &  8\\

$\mathfrak{P}||B$     &   
\begin{tikzcd}
* \arrow[r, bend left] & *, \arrow[l, bend left] & * \arrow[r, bend left] & * \arrow[l, bend left]
\end{tikzcd}  &  3\\

$\mathfrak{P}|||$     &   
\begin{tikzpicture}
    \node (a) at (0,0) {$*$};
    \node (b) at (2,0) {$*$};
    \node (c) at (4,0) {$*$};
    \node (d) at (6,0) {$*$};

    \draw (a) edge[loop right] node {,} (a);
    \draw (b) edge[loop right] node {,} (b);
    \draw[->] (c) to [in=170, out=10, looseness=1] node {} (d);
    \draw[->] (d) to [in=-10, out=190, looseness=1] node {} (c);
\end{tikzpicture} & 6\\

$\mathfrak{P}||||$     &   
\begin{tikzpicture}
    \node (a) at (0,0) {$*$};
    \node (b) at (2,0) {$*$};
    \node (c) at (4,0) {$*$};
    \node (d) at (6,0) {$*$};

    \draw (a) edge[loop right] node {,} (a);
    \draw (b) edge[loop right] node {,} (b);
    \draw (c) edge[loop right] node {,} (c);
    \draw (d) edge[loop right] node {} (d);
    
\end{tikzpicture} & 1\\

\hline 
\end{tabular}
\label{tab7} 
\end{table}

It turns out that all the period maps on sets of four points can be generated by merely two ones in them.

\begin{theorem}\label{thm12}
Assuming we can realize the maps $\mathfrak{P}|1$ and $\mathfrak{P}||A1$ on $X=\{p_1, p_2, p_3, p_4\}$  by some rational maps 
\begin{center}
$R_{\mathfrak{P}|1}: \mathbb{P}^1(\mathbb{C})\rightarrow\mathbb{P}^1(\mathbb{C})$ and $R_{\mathfrak{P}||A1}: \mathbb{P}^1(\mathbb{C})\rightarrow\mathbb{P}^1(\mathbb{C})$
\end{center}
respectively with $V(R_{\mathfrak{P}|1})=V(R_{\mathfrak{P}||A1})=X$, then all the periodic maps on $X$ are realizable by rational maps on $\mathbb{P}^1(\mathbb{C})$.
\end{theorem} 

\begin{proof}
We can realize every periodic map on $\{p_1, p_2, p_3, p_4\}$ by compositions of the two rational maps 
$R_{\mathfrak{P}|1}$ and $R_{\mathfrak{P}||A1}$ in fact. For simplicity we only realize one type in every series in the periodic family, enthusiastic readers may supplement realization of the rest types. 
\begin{itemize}
\item Map \begin{tikzpicture}
    \node (a) at (0,0) {$p_1$};
    \node (b) at (2,0) {$p_4$};
    \node (c) at (4,0) {$p_3$};
    \node (d) at (6,0) {$p_2$};

    \draw[->] (a) to node {} (b);
    \draw[->] (b) to node {} (c);
    \draw[->] (c) to node {} (d);
    \draw[->] (d) to [in=-10, out=190, looseness=1] node {} (a);
\end{tikzpicture}. 

This map can be realized by the rational map $R_{\mathfrak{P}|1}^3$. 

\item Map \begin{tikzpicture}
    \node (a) at (0,0) {$p_1$};
    \node (b) at (2,0) {$p_4$};
    \node (c) at (4,0) {$p_2,$};
    \node (d) at (5,0) {$p_3$};

    \draw[->] (a) to node  {} (b);
    \draw[->] (b) to node {} (c);
    \draw[->] (c) to [in=-10, out=190, looseness=1] node {} (a);
    \draw (d) edge[loop right] node {} (d);
\end{tikzpicture}.

This map can be realized by the rational map $R_{\mathfrak{P}|1}^2\circ R_{\mathfrak{P}||A1}$. 

\item Map \begin{tikzcd}
p_1 \arrow[r, bend left] & p_2, \arrow[l, bend left] & p_3 \arrow[r, bend left] & p_4 \arrow[l, bend left]
\end{tikzcd}. 

This map can be realized by the rational map $(R_{\mathfrak{P}|1}\circ R_{\mathfrak{P}||A1})^2$.

\item Map \begin{tikzpicture}
    \node (a) at (0,0) {$p_2$};
    \node (b) at (2,0) {$p_4$};
    \node (c) at (4,0) {$p_1$};
    \node (d) at (6,0) {$p_3$};

    \draw (a) edge[loop right] node {,} (a);
    \draw (b) edge[loop right] node {,} (b);
    \draw[->] (c) to [in=170, out=10, looseness=1] node {} (d);
    \draw[->] (d) to [in=-10, out=190, looseness=1] node {} (c);
\end{tikzpicture}. 

This map can be realized by the rational map $R_{\mathfrak{P}|1}^2\circ R_{\mathfrak{P}||A1}^2 \circ (R_{\mathfrak{P}||A1}\circ R_{\mathfrak{P}|1})^3$.

\item Map $\mathfrak{P}||||$. This map can be realized by the rational map $R_{\mathfrak{P}|1}^4$.
\end{itemize}

It is easy to see that any of the critical-value sets of the above compositive maps is $X$ in virtue of Corollary \ref{cor5}. 
\end{proof}

Obviously there are other combinations of periodic maps which can also generate all the maps on $X$, however, one needs at least two of them.

\begin{coro}\label{cor8} 
Assuming we can realize the maps $\mathfrak{P}|1$ and $\mathfrak{P}|6$ on $X=\{p_1, p_2, p_3, p_4\}$  by some rational maps 
\begin{center}
$R_{\mathfrak{P}|1}: \mathbb{P}^1(\mathbb{C})\rightarrow\mathbb{P}^1(\mathbb{C})$ and $R_{\mathfrak{P}|6}: \mathbb{P}^1(\mathbb{C})\rightarrow\mathbb{P}^1(\mathbb{C})$ 
\end{center}
respectively with $V(R_{\mathfrak{P}|1})=X$, then all the periodic maps on $X$ are realizable by rational maps on $\mathbb{P}^1(\mathbb{C})$.
\end{coro}
\begin{proof}
It suffices for us to realize $\mathfrak{P}||A1$ on $X$ by some rational map on $\mathbb{P}^1(\mathbb{C})$ whose critical-value set coincides with $X$, in virtue of Theorem \ref{thm12}. In fact one can check that the map $\mathfrak{P}||A1$ can be realized by the rational map $R_{\mathfrak{P}|1}^2\circ R_{\mathfrak{P}|6}\circ R_{\mathfrak{P}|1}$ with $V(R_{\mathfrak{P}|1}^2\circ R_{\mathfrak{P}|6}\circ R_{\mathfrak{P}|1})=X$, under assumption of the corollary. 
\end{proof}

In fact the condition $V(R_{\mathfrak{P}|1})=V(R_{\mathfrak{P}||A1})=X$ in Theorem \ref{thm12} can also be relaxed to $V(R_{\mathfrak{P}|1})=X$.  The collection of all the periodic maps on $X$ can be viewed as a non-commutative finite group under the composition operation. We suspect there might be some shortcut to find the generating sets of the group, however, we do not know such method.

\subsection{The strictly pre-periodic family on sets of four points}

There are totally 232 types of strictly pre-periodic maps on $X=\{p_1, p_2, p_3, p_4\}$ of four points. These maps constitute a finite semi-group under the composition operation (which is a monoid by including the identity map $\mathfrak{P}||||$ in them).  We also give a classification of them according to the number of disjoint orbits in them. We first count the number of types of the series of strictly pre-periodic maps with only one orbit in Table \ref{tab10}, then we count the number of types of the series of the strictly pre-periodic maps with two orbits in Table \ref{tab13}. There is only one series ($\mathfrak{S}|||$) of strictly pre-periodic maps with three orbits (12 types). 
\begin{center}
$\mathfrak{S}|||$: 
\begin{tikzpicture}
    \node (a) at (0,0) {$*$};
    \node (b) at (2,0) {$*$};
    \node (c) at (4,0) {$*$};
    \node (d) at (6,0) {$*$};

    \draw (a) edge[loop right] node {,} (a);
    \draw (b) edge[loop right] node {,} (b);
    \draw[->] (c) to node {} (d);
    \draw (d) edge[loop right] node {.} (d);
\end{tikzpicture}
\end{center}

One can add an additional integer after the name of the same series to distinguish various type of maps in the series. For example, let $\{\mathfrak{S}|Aj\}_{j=1}^{24}$ be the 24 types of maps in the $\mathfrak{S}|A$ series with   
\begin{center}
$\mathfrak{S}|A1$:        
\begin{tikzpicture}
    \node (a) at (0,0) {$p_1$};
    \node (b) at (2,0) {$p_2$};
    \node (c) at (4,0) {$p_3$};
    \node (d) at (6,0) {$p_4$};

    \draw[->] (a) to node {} (b);
    \draw[->] (b) to node {} (c);
    \draw[->] (c) to node {} (d);
    \draw[->] (d) to [in=-10, out=190, looseness=1] node {} (b);
\end{tikzpicture}.
\end{center}
Note that twice-iterate operation of the maps in the $\mathfrak{S}|A$ series is closed, so the series admits a partition of two subsets
\begin{center}
$\{\mathfrak{S}|Aj\}_{j=1}^{24}=\{\mathfrak{S}|Aj\}_{j=1}^{12}\cup\{\mathfrak{S}|Aj\}_{j=13}^{24}$
\end{center} 
such that $(\mathfrak{S}|Aj)^2\notin \{\mathfrak{S}|Aj\}_{j=1}^{12}$ for any $1\leq j\leq 12$. Let

\begin{center}
$\mathfrak{S}|F1$:   
\begin{tikzpicture}
    \node (a) at (0,0) {$p_1$};
    \node (b) at (2,0) {$p_2$};
    \node (c) at (4,0) {$p_3$};
    \node (d) at (6,0) {$p_4$};

    \draw[->] (a) to node {} (b);
    \draw[->] (b) to node {} (c);
    \draw[->] (d) to node {} (c);
    \draw[->] (c) to [in=25, out=155, looseness=2] node {} (c);  
\end{tikzpicture}.
\end{center}

\begin{table}[ht]
\caption{Number of strictly pre-periodic maps with one orbit} 
\centering 
\begin{tabular}[ht]{ccc}
\hline 
Series &  map & number of types\\
\hline
$\mathfrak{S}|A$     &   
\begin{tikzpicture}
    \node (a) at (0,0) {$*$};
    \node (b) at (2,0) {$*$};
    \node (c) at (4,0) {$*$};
    \node (d) at (6,0) {$*$};

    \draw[->] (a) to node {} (b);
    \draw[->] (b) to node {} (c);
    \draw[->] (c) to node {} (d);
    \draw[->] (d) to [in=-10, out=190, looseness=1] node {} (b);
\end{tikzpicture}
& 24\\

$\mathfrak{S}|B$     &   
\begin{tikzpicture}
    \node (a) at (0,0) {$*$};
    \node (b) at (2,0) {$*$};
    \node (c) at (4,0) {$*$};
    \node (d) at (6,0) {$*$};

    \draw[->] (a) to node {} (b);
    \draw[->] (b) to node {} (c);
    \draw[->] (c) to node {} (d);
    \draw[->] (d) to [in=-10, out=190, looseness=1] node {} (c);
\end{tikzpicture}
& 24\\

$\mathfrak{S}|C$     &   
\begin{tikzpicture}
    \node (a) at (0,0) {$*$};
    \node (b) at (2,0) {$*$};
    \node (c) at (4,0) {$*$};
    \node (d) at (6,0) {$*$};

    \draw[->] (a) to node {} (b);
    \draw[->] (b) to [in=170, out=10, looseness=1] node {} (c);
    \draw[->] (c) to [in=-10, out=190, looseness=1] node {} (b);
    \draw[->] (d) to node {} (c);
\end{tikzpicture}
& 12\\

$\mathfrak{S}|D$     &   
\begin{tikzpicture}
    \node (a) at (0:1) {$*$};
    \node (d) at (170:1) {$*$}; 
    \node (e) at (190:1) {$*$};
    \node [right of = a, node distance=0.8in] (b) {$*$};

    \draw[->] (d) to  node  {} (a);
    \draw[->] (e) to node  {} (a);
    \draw[->] (b) to [in=10, out=170, looseness=1] node  {} (a);
    \draw[->] (a) to [in=190, out=-10, looseness=1] node  {} (b);
\end{tikzpicture}
& 12\\

$\mathfrak{S}|E$     &   
\begin{tikzpicture}
    \node (a) at (0,0) {$*$};
    \node (b) at (2,0) {$*$};
    \node (c) at (4,0) {$*$};
    \node (d) at (6,0) {$*$};

    \draw[->] (a) to node {} (b);
    \draw[->] (b) to node {} (c);
    \draw[->] (c) to node {} (d);
    \draw (d) edge[loop right] node {} (d);
\end{tikzpicture}
& 24\\

$\mathfrak{S}|F$     &   
\begin{tikzpicture}
    \node (a) at (0,0) {$*$};
    \node (b) at (2,0) {$*$};
    \node (c) at (4,0) {$*$};
    \node (d) at (6,0) {$*$};

    \draw[->] (a) to node {} (b);
    \draw[->] (b) to node {} (c);
    \draw[->] (d) to node {} (c);
    \draw[->] (c) to [in=25, out=155, looseness=2] node {} (c);
\end{tikzpicture}
& 24\\

$\mathfrak{S}|G$     &   
\begin{tikzpicture}
    \node (a) at (0:1) {$*$};
    \node (d) at (170:1) {$*$}; 
    \node (e) at (190:1) {$*$};
    \node [right of = a, node distance=0.8in] (b) {$*$};

    \draw[->] (d) to  node  {} (a);
    \draw[->] (e) to node  {} (a);
    \draw[->] (a) to node  {} (b);
    \draw (b) edge[loop right] node {} (b);
\end{tikzpicture}
& 12\\

$\mathfrak{S}|H$     &   
\begin{tikzpicture}
    \node (a) at (0:1) {$*$};
    \node (d) at (170:1) {$*$}; 
    \node (e) at (190:1) {$*$};
    \node (b) at (180:1) {$*$};

    \draw[->] (d) to  node  {} (a);
    \draw[->] (e) to node  {} (a);
    \draw[->] (b) to node  {} (a);
    \draw (a) edge[loop right] node {} (a);
\end{tikzpicture}
& 4\\

\hline 
\end{tabular}
\label{tab10} 
\end{table}

\begin{table}[ht]
\caption{Number of strictly pre-periodic maps with two orbits} 
\centering 
\begin{tabular}[t]{ccc}
\hline 
Series &  map & number of types\\
\hline
$\mathfrak{S}||A$     &   
\begin{tikzpicture}
    \node (a) at (0,0) {$*$};
    \node (b) at (2,0) {$*$};
    \node (c) at (4,0) {$*$};
    \node (d) at (6,0) {$*$};

    \draw (a) edge[loop right] node {,} (a);
    \draw[->] (b) to  node {} (c);
    \draw[->] (c) to [in=170, out=10, looseness=1] node {} (d);
    \draw[->] (d) to [in=-10, out=-170, looseness=1] node {} (c);
\end{tikzpicture}
& 24\\

$\mathfrak{S}||B$     &   
\begin{tikzpicture}
    \node (a) at (0,0) {$*$};
    \node (b) at (2,0) {$*$};
    \node (c) at (4,0) {$*$};
    \node (d) at (6,0) {$*$};

    \draw (a) edge[loop right] node {,} (a);
    \draw[->] (b) to  node {} (c);
    \draw[->] (c) to node {} (d);
    \draw (d) edge[loop right] node {} (d);
\end{tikzpicture}
& 24\\

$\mathfrak{S}||C$     &   
\begin{tikzpicture}
    \node (a) at (0,0) {$*$};
    \node (b) at (2,0) {$*$};
    \node (c) at (4,0) {$*$};
    \node (d) at (6,0) {$*$};

    \draw (a) edge[loop right] node {,} (a);
    \draw[->] (b) to node {} (c);
    \draw[->] (d) to node {} (c);
    \draw[->] (c) to [in=25, out=155, looseness=2] node {} (c);
\end{tikzpicture}
& 12\\

$\mathfrak{S}||D$     &   
\begin{tikzpicture}
    \node (a) at (0,0) {$*$};
    \node (b) at (2,0) {$*$};
    \node (c) at (4,0) {$*$};
    \node (d) at (6,0) {$*$};

    \draw (b) edge[loop right] node {,} (b);
    \draw[->] (a) to node {} (b);
    \draw[->] (c) to node {} (d);
    \draw (d) edge[loop right] node {} (d);
\end{tikzpicture}
& 12\\

$\mathfrak{S}||E$     &   
\begin{tikzpicture}
    \node (a) at (0,0) {$*$};
    \node (b) at (2,0) {$*,$};
    \node (c) at (4,0) {$*$};
    \node (d) at (6,0) {$*$};

    \draw[->] (a) to [in=170, out=10, looseness=1] node {} (b);
    \draw[->] (b) to [in=-10, out=-170, looseness=1] node {} (a);
    \draw[->] (c) to node {} (d);
    \draw (d) edge[loop right] node {} (d);
\end{tikzpicture}
& 12\\

\hline 
\end{tabular}
\label{tab13} 
\end{table}

It turns out that the semi-group of strictly pre-periodic maps can be generated by the 12 ones in $\{\mathfrak{S}|Aj\}_{j=1}^{12}$ (or $\{\mathfrak{S}|Aj\}_{j=13}^{24}$).

\begin{theorem}\label{thm13}
All the strictly pre-periodic maps on $X=\{p_1, p_2, p_3, p_4\}$ can be given by finite compositions of maps in 
\begin{center}
$\{\mathfrak{S}|Aj\}_{j=1}^{12}$
\end{center}
in the $\mathfrak{S}|A$ series.
\end{theorem}

Theorem \ref{thm13} induces the following result on realizing the strictly pre-periodic maps on $X=\{p_1, p_2, p_3, p_4\}$ instantly. 

\begin{theorem}\label{thm14}
Assuming we can realize the maps in $\{\mathfrak{S}|Aj\}_{j=1}^{12}$ on $X=\{p_1, p_2, p_3, p_4\}$  by some rational maps 
\begin{center}
$\{R_{\mathfrak{S}|Aj}: \mathbb{P}^1(\mathbb{C})\rightarrow\mathbb{P}^1(\mathbb{C})\}_{1\leq j\leq 12}$ 
\end{center}
respectively with $V(R_{\mathfrak{S}|Aj})=X$ for any $1\leq j\leq 12$, then all the strictly pre-periodic maps on $X$ are realizable by rational maps on $\mathbb{P}^1(\mathbb{C})$.
\end{theorem}

\begin{proof}
The result follows from a combination of Theorem \ref{thm13} and Corollary \ref{cor5} instantly.
\end{proof}

In the following we gear towards the proof of Theorem \ref{thm13}. We start from the following result.

\begin{lemma}\label{lem21}
All the maps in the $\mathfrak{S}|A$ series can be given by finite compositions of maps in $\{\mathfrak{S}|Aj\}_{j=1}^{12}$.
\end{lemma}
\begin{proof}
This is simply because any map in $\{\mathfrak{S}|Aj\}_{j=13}^{24}$ arises as twice iterate of some map in $\{\mathfrak{S}|Aj\}_{j=1}^{12}$.
\end{proof}

Then we show that the maps in the $\mathfrak{S}|F$ series can be given by compositions of maps in the $\mathfrak{S}|A$ series in fact.

\begin{lemma}\label{lem19}
All the maps in the $\mathfrak{S}|F$ series can be given by finite compositions of maps in the $\mathfrak{S}|A$ series.
\end{lemma}
\begin{proof}
First note that 
\begin{equation}\label{eq63}
\mathfrak{S}|F1= \begin{tikzpicture}
    \node (a) at (5,0) {$p_1$};
    \node (b) at (6,0) {$p_3$};
    \node (c) at (7,0) {$p_2$};
    \node (d) at (8,0) {$p_4$};

    \draw[->] (a) to node {} (b);
    \draw[->] (b) to node {} (c);
    \draw[->] (c) to node {} (d);
    \draw[->] (d) to [in=-20, out=200, looseness=1] node {} (b);
\end{tikzpicture}
\circ
\begin{tikzpicture}
    \node (a) at (9,0) {$p_2$};
    \node (b) at (10,0) {$p_1$};
    \node (c) at (11,0) {$p_3$};
    \node (d) at (12,0) {$p_4$};

    \draw[->] (a) to node {} (b);
    \draw[->] (b) to node {} (c);
    \draw[->] (c) to node {} (d);
    \draw[->] (d) to [in=-20, out=200, looseness=1] node {} (b);
\end{tikzpicture}.
\end{equation}

The rest 23 maps in the $\mathfrak{S}|F$ series can be given by compositions of two maps in the $\mathfrak{S}|A$ series by alternating subscripts of points in (\ref{eq63}).  
\end{proof}

One can probably grasp the idea that all types of maps in the same series can be given by compositions of maps in certain series as long as one of them can be done.  On the inspiration of Lemma \ref{lem19} we can show the following result.

\begin{lemma}\label{lem20}
All the strictly pre-periodic maps can be given by finite compositions of maps in the $\mathfrak{S}|A$ series.
\end{lemma}
\begin{proof}
We justify the result according to the series of the strictly pre-periodic maps inductively. Some details are left to the readers.
\begin{itemize}
\item Series $\mathfrak{S}|||$. Maps in this series can be given by the composition $\mathfrak{S}|Aj\circ \mathfrak{S}|Ak$ for some $j,k\in\{1,2,\cdots,24\}$. For example, 

\begin{center}
\begin{tikzpicture}
    \node (a) at (0,0) {$p_1$};
    \node (b) at (1,0) {$p_3$};
    \node (c) at (2.5,0) {$p_2$};
    \node (d) at (4,0) {$p_4$};

    \draw[->] (a) to node {} (b);
    \draw (b) edge[loop right] node {,} (b);
    \draw (c) edge[loop right] node {,} (c);
    \draw (d) edge[loop right] node {} (d);     
\end{tikzpicture}
=
\begin{tikzpicture}
    \node (a) at (5,0) {$p_1$};
    \node (b) at (6,0) {$p_2$};
    \node (c) at (7,0) {$p_3$};
    \node (d) at (8,0) {$p_4$};

    \draw[->] (a) to node {} (b);
    \draw[->] (b) to node {} (c);
    \draw[->] (c) to node {} (d);
    \draw[->] (d) to [in=-20, out=200, looseness=1] node {} (b);
\end{tikzpicture}
$\circ$
\begin{tikzpicture}
    \node (a) at (9,0) {$p_1$};
    \node (b) at (10,0) {$p_2$};
    \node (c) at (11,0) {$p_4$};
    \node (d) at (12,0) {$p_3$};

    \draw[->] (a) to node {} (b);
    \draw[->] (b) to node {} (c);
    \draw[->] (c) to node {} (d);
    \draw[->] (d) to [in=-20, out=200, looseness=1] node {} (b);
\end{tikzpicture}
\end{center}

\item Series $\mathfrak{S}|F$. Maps in this series can be given by the composition $\mathfrak{S}|Aj\circ \mathfrak{S}|Ak$ for some $j,k\in\{1,2,\cdots,24\}$, in virtue of Lemma \ref{lem19}.

\item Series $\mathfrak{S}|H$. Maps in this series can be given by the composition $\mathfrak{S}|Fj\circ \mathfrak{S}|Ak$ for some $j,k\in\{1,2,\cdots,24\}$. For example, 

\begin{center}
\begin{tikzpicture}
    \node (a) at (0:1) {$p_2$};
    \node (d) at (160:1) {$p_1$}; 
    \node (e) at (200:1) {$p_3$};
    \node (b) at (180:1) {$p_4$};

    \draw[->] (d) to  node  {} (a);
    \draw[->] (e) to node  {} (a);
    \draw[->] (b) to node  {} (a);
    \draw (a) edge[loop right] node {} (a);
\end{tikzpicture}
=
\begin{tikzpicture}
    \node (a) at (5,0) {$p_1$};
    \node (b) at (6,0) {$p_4$};
    \node (c) at (7,0) {$p_2$};
    \node (d) at (8,0) {$p_3$};

    \draw[->] (a) to node {} (b);
    \draw[->] (b) to node {} (c);
    \draw[->] (d) to node {} (c);
    \draw[->] (c) to [in=25, out=155, looseness=2] node {} (c);    
\end{tikzpicture}
$\circ$
\begin{tikzpicture}
    \node (a) at (9,0) {$p_1$};
    \node (b) at (10,0) {$p_2$};
    \node (c) at (11,0) {$p_3$};
    \node (d) at (12,0) {$p_4$};

    \draw[->] (a) to node {} (b);
    \draw[->] (b) to node {} (c);
    \draw[->] (c) to node {} (d);
    \draw[->] (d) to [in=-20, out=200, looseness=1] node {} (b);
\end{tikzpicture}
\end{center}

\item Series $\mathfrak{S}|G$. Maps in this series can be given by the composition $\mathfrak{S}|Fj\circ \mathfrak{S}|Ak$ for some $j,k\in\{1,2,\cdots,24\}$. For example, 

\begin{center}
\begin{tikzpicture}
    \node (a) at (0:1) {$p_2$};
    \node (d) at (160:0.5) {$p_1$}; 
    \node (e) at (200:0.5) {$p_3$};
    \node [right of = a, node distance=0.5in] (b) {$p_4$};

    \draw[->] (d) to  node  {} (a);
    \draw[->] (e) to node  {} (a);
    \draw[->] (a) to node  {} (b);
    \draw (b) edge[loop right] node {} (b);
\end{tikzpicture}
=
\begin{tikzpicture}
    \node (a) at (5,0) {$p_1$};
    \node (b) at (6,0) {$p_4$};
    \node (c) at (7,0) {$p_2$};
    \node (d) at (8,0) {$p_3$};

    \draw[->] (a) to node {} (b);
    \draw[->] (b) to node {} (c);
    \draw[->] (d) to node {} (c);
    \draw[->] (c) to [in=25, out=155, looseness=2] node {} (c);    
\end{tikzpicture}
$\circ$
\begin{tikzpicture}
    \node (a) at (9,0) {$p_2$};
    \node (b) at (10,0) {$p_1$};
    \node (c) at (11,0) {$p_3$};
    \node (d) at (12,0) {$p_4$};

    \draw[->] (a) to node {} (b);
    \draw[->] (b) to node {} (c);
    \draw[->] (c) to node {} (d);
    \draw[->] (d) to [in=-20, out=200, looseness=1] node {} (b);
\end{tikzpicture}
\end{center}

\item Series $\mathfrak{S}||D$. Maps in this series can be given by the composition $\mathfrak{S}|Fj\circ \mathfrak{S}|Ak$ for some $j,k\in\{1,2,\cdots,24\}$. For example, 
\begin{center}
\begin{tikzpicture}
    \node (a) at (0,0) {$p_1$};
    \node (b) at (1,0) {$p_3$};
    \node (c) at (2.5,0) {$p_4$};
    \node (d) at (3.5,0) {$p_2$};

    \draw (b) edge[loop right] node {,} (b);
    \draw[->] (a) to node {} (b);
    \draw[->] (c) to node {} (d);
    \draw (d) edge[loop right] node {} (d);
\end{tikzpicture}
=
\begin{tikzpicture}
    \node (a) at (5,0) {$p_1$};
    \node (b) at (6,0) {$p_2$};
    \node (c) at (7,0) {$p_3$};
    \node (d) at (8,0) {$p_4$};

    \draw[->] (a) to node {} (b);
    \draw[->] (b) to node {} (c);
    \draw[->] (d) to node {} (c);
    \draw[->] (c) to [in=25, out=155, looseness=2] node {} (c);    
\end{tikzpicture}
$\circ$
\begin{tikzpicture}
    \node (a) at (9,0) {$p_2$};
    \node (b) at (10,0) {$p_1$};
    \node (c) at (11,0) {$p_3$};
    \node (d) at (12,0) {$p_4$};

    \draw[->] (a) to node {} (b);
    \draw[->] (b) to node {} (c);
    \draw[->] (c) to node {} (d);
    \draw[->] (d) to [in=-20, out=200, looseness=1] node {} (b);
\end{tikzpicture}
\end{center}

\item Series $\mathfrak{S}||C$. Maps in this series can be given by the composition $\mathfrak{S}|Gj\circ \mathfrak{S}|Ak$ for some $j\in\{1,2,\cdots,12\}, k\in\{1,2,\cdots,24\}$. For example, 
\begin{center}
\begin{tikzpicture}
    \node (a) at (0,0) {$p_2$};
    \node (b) at (1.5,0) {$p_1$};
    \node (c) at (2.5,0) {$p_4$};
    \node (d) at (3.5,0) {$p_3$};

    \draw (a) edge[loop right] node {,} (a);
    \draw[->] (b) to node {} (c);
    \draw[->] (d) to node {} (c);
    \draw[->] (c) to [in=25, out=155, looseness=2] node {} (c);
\end{tikzpicture}
=
\begin{tikzpicture}
    \node (a) at (5:1) {$p_2$};
    \node (d) at (160:0.5) {$p_1$}; 
    \node (e) at (200:0.5) {$p_3$};
    \node [right of = a, node distance=0.4in] (b) {$p_4$};

    \draw[->] (d) to  node  {} (a);
    \draw[->] (e) to node  {} (a);
    \draw[->] (a) to node  {} (b);
    \draw (b) edge[loop right] node {} (b);
\end{tikzpicture}
$\circ$
\begin{tikzpicture}
    \node (a) at (9,0) {$p_1$};
    \node (b) at (10,0) {$p_2$};
    \node (c) at (11,0) {$p_3$};
    \node (d) at (12,0) {$p_4$};

    \draw[->] (a) to node {} (b);
    \draw[->] (b) to node {} (c);
    \draw[->] (c) to node {} (d);
    \draw[->] (d) to [in=-20, out=200, looseness=1] node {} (b);
\end{tikzpicture}
\end{center}

\item Series $\mathfrak{S}|C$. Maps in this series can be given by the composition $\mathfrak{S}|Fj\circ \mathfrak{S}|Ak$ for some $j,k\in\{1,2,\cdots,24\}$. For example, 
\begin{center}
\begin{tikzpicture}
    \node (a) at (0,0) {$p_1$};
    \node (b) at (1,0) {$p_4$};
    \node (c) at (2,0) {$p_3$};
    \node (d) at (3,0) {$p_2$};

    \draw[->] (a) to node {} (b);
    \draw[->] (b) to [in=170, out=10, looseness=1] node {} (c);
    \draw[->] (c) to [in=-10, out=190, looseness=1] node {} (b);
    \draw[->] (d) to node {} (c);
\end{tikzpicture}
=
\begin{tikzpicture}
    \node (a) at (5,0) {$p_1$};
    \node (b) at (6,0) {$p_3$};
    \node (c) at (7,0) {$p_4$};
    \node (d) at (8,0) {$p_2$};

    \draw[->] (a) to node {} (b);
    \draw[->] (b) to node {} (c);
    \draw[->] (d) to node {} (c);
    \draw[->] (c) to [in=25, out=155, looseness=2] node {} (c);    
\end{tikzpicture}
$\circ$
\begin{tikzpicture}
    \node (a) at (9,0) {$p_2$};
    \node (b) at (10,0) {$p_1$};
    \node (c) at (11,0) {$p_3$};
    \node (d) at (12,0) {$p_4$};

    \draw[->] (a) to node {} (b);
    \draw[->] (b) to node {} (c);
    \draw[->] (c) to node {} (d);
    \draw[->] (d) to [in=-20, out=200, looseness=1] node {} (b);
\end{tikzpicture}
\end{center}

\item Series $\mathfrak{S}|D$. Maps in this series can be given by the composition $\mathfrak{S}||Dj\circ \mathfrak{S}|Ak$ for some $j\in\{1,2,\cdots,12\}, k\in\{1,2,\cdots,24\}$. For example, 
\begin{center}
\begin{tikzpicture}
    \node (a) at (0:1) {$p_2$};
    \node (d) at (160:0.5) {$p_1$}; 
    \node (e) at (200:0.5) {$p_4$};
    \node [right of = a, node distance=0.4in] (b) {$p_3$};

    \draw[->] (d) to  node  {} (a);
    \draw[->] (e) to node  {} (a);
    \draw[->] (b) to [in=10, out=170, looseness=1] node  {} (a);
    \draw[->] (a) to [in=190, out=-10, looseness=1] node  {} (b);
\end{tikzpicture}
=
\begin{tikzpicture}
    \node (a) at (5,0) {$p_1$};
    \node (b) at (6,0) {$p_3$};
    \node (c) at (7.5,0) {$p_4$};
    \node (d) at (8.5,0) {$p_2$};

    \draw (b) edge[loop right] node {,} (b);
    \draw[->] (a) to node {} (b);
    \draw[->] (c) to node {} (d);
    \draw (d) edge[loop right] node {} (d);
\end{tikzpicture}
$\circ$
\begin{tikzpicture}
    \node (a) at (9,0) {$p_1$};
    \node (b) at (10,0) {$p_2$};
    \node (c) at (11,0) {$p_3$};
    \node (d) at (12,0) {$p_4$};

    \draw[->] (a) to node {} (b);
    \draw[->] (b) to node {} (c);
    \draw[->] (c) to node {} (d);
    \draw[->] (d) to [in=-20, out=200, looseness=1] node {} (b);
\end{tikzpicture}
\end{center}

\item Series $\mathfrak{S}|E$. Maps in this series can be given by the composition $\mathfrak{S}|Aj\circ \mathfrak{S}|Ak$ for some $j, k\in\{1,2,\cdots,24\}$. For example, 
\begin{center}
\begin{tikzpicture}
    \node (a) at (0,0) {$p_3$};
    \node (b) at (1,0) {$p_1$};
    \node (c) at (2,0) {$p_4$};
    \node (d) at (3,0) {$p_2$};

    \draw[->] (a) to node {} (b);
    \draw[->] (b) to node {} (c);
    \draw[->] (c) to node {} (d);
    \draw (d) edge[loop right] node {} (d);
\end{tikzpicture}
=
\begin{tikzpicture}
    \node (a) at (5,0) {$p_3$};
    \node (b) at (6,0) {$p_4$};
    \node (c) at (7,0) {$p_1$};
    \node (d) at (8,0) {$p_2$};

    \draw[->] (a) to node {} (b);
    \draw[->] (b) to node {} (c);
    \draw[->] (c) to node {} (d);
    \draw[->] (d) to [in=-20, out=200, looseness=1] node {} (b);
\end{tikzpicture}
$\circ$
\begin{tikzpicture}
    \node (a) at (9,0) {$p_2$};
    \node (b) at (10,0) {$p_1$};
    \node (c) at (11,0) {$p_3$};
    \node (d) at (12,0) {$p_4$};

    \draw[->] (a) to node {} (b);
    \draw[->] (b) to node {} (c);
    \draw[->] (c) to node {} (d);
    \draw[->] (d) to [in=-20, out=200, looseness=1] node {} (b);
\end{tikzpicture}
\end{center}

\item Series $\mathfrak{S}||B$. Maps in this series can be given by the composition $\mathfrak{S}|Aj\circ \mathfrak{S}|Ak$ for some $j, k\in\{1,2,\cdots,24\}$. For example, 
\begin{center}
\begin{tikzpicture}
    \node (a) at (0,0) {$p_1$};
    \node (b) at (1.5,0) {$p_3$};
    \node (c) at (2.5,0) {$p_2$};
    \node (d) at (3.5,0) {$p_4$};

    \draw (a) edge[loop right] node {,} (a);
    \draw[->] (b) to  node {} (c);
    \draw[->] (c) to node {} (d);
    \draw (d) edge[loop right] node {} (d);
\end{tikzpicture}
=
\begin{tikzpicture}
    \node (a) at (5,0) {$p_3$};
    \node (b) at (6,0) {$p_1$};
    \node (c) at (7,0) {$p_4$};
    \node (d) at (8,0) {$p_2$};

    \draw[->] (a) to node {} (b);
    \draw[->] (b) to node {} (c);
    \draw[->] (c) to node {} (d);
    \draw[->] (d) to [in=-20, out=200, looseness=1] node {} (b);
\end{tikzpicture}
$\circ$
\begin{tikzpicture}
    \node (a) at (9,0) {$p_2$};
    \node (b) at (10,0) {$p_1$};
    \node (c) at (11,0) {$p_3$};
    \node (d) at (12,0) {$p_4$};

    \draw[->] (a) to node {} (b);
    \draw[->] (b) to node {} (c);
    \draw[->] (c) to node {} (d);
    \draw[->] (d) to [in=-20, out=200, looseness=1] node {} (b);
\end{tikzpicture}
\end{center}

\item Series $\mathfrak{S}|B$. Maps in this series can be given by the composition $\mathfrak{S}|||j\circ \mathfrak{S}|Ak$ for some $j\in\{1,2,\cdots,12\}, k\in\{1,2,\cdots,24\}$. For example, 
\begin{center}
\begin{tikzpicture}
    \node (a) at (0,0) {$p_3$};
    \node (b) at (1,0) {$p_4$};
    \node (c) at (2,0) {$p_1$};
    \node (d) at (3.5,0) {$p_2$};

    \draw[->] (a) to node {} (b);
    \draw[->] (b) to node {} (c);
    \draw[->] (c) to node {} (d);
    \draw[->] (d) to [in=-10, out=190, looseness=1] node {} (c);
\end{tikzpicture}
=
\begin{tikzpicture}
    \node (a) at (4,0) {$p_3$};
    \node (b) at (5,0) {$p_2$};
    \node (c) at (6.5,0) {$p_1$};
    \node (d) at (8,0) {$p_4$};

    \draw[->] (a) to node {} (b);
    \draw (b) edge[loop right] node {,} (b);
    \draw (c) edge[loop right] node {,} (c);
    \draw (d) edge[loop right] node {} (d);     
\end{tikzpicture}
$\circ$
\begin{tikzpicture}
    \node (a) at (9,0) {$p_2$};
    \node (b) at (10,0) {$p_1$};
    \node (c) at (11,0) {$p_3$};
    \node (d) at (12,0) {$p_4$};

    \draw[->] (a) to node {} (b);
    \draw[->] (b) to node {} (c);
    \draw[->] (c) to node {} (d);
    \draw[->] (d) to [in=-20, out=200, looseness=1] node {} (b);
\end{tikzpicture}
\end{center}

\item Series $\mathfrak{S}||A$. Maps in this series can be given by the composition $\mathfrak{S}|Aj\circ \mathfrak{S}|Bk$ for some $j, k\in\{1,2,\cdots,24\}$. For example, 
\begin{center}
\begin{tikzpicture}
    \node (a) at (0,0) {$p_3$};
    \node (b) at (1.5,0) {$p_1$};
    \node (c) at (2.5,0) {$p_4$};
    \node (d) at (3.5,0) {$p_2$};

    \draw (a) edge[loop right] node {,} (a);
    \draw[->] (b) to  node {} (c);
    \draw[->] (c) to [in=170, out=10, looseness=1] node {} (d);
    \draw[->] (d) to [in=-10, out=-170, looseness=1] node {} (c);
\end{tikzpicture}
=
\begin{tikzpicture}
    \node (a) at (5,0) {$p_1$};
    \node (b) at (6,0) {$p_2$};
    \node (c) at (7,0) {$p_3$};
    \node (d) at (8,0) {$p_4$};

    \draw[->] (a) to node {} (b);
    \draw[->] (b) to node {} (c);
    \draw[->] (c) to node {} (d);
    \draw[->] (d) to [in=-20, out=200, looseness=1] node {} (b);
\end{tikzpicture}
$\circ$
\begin{tikzpicture}
    \node (a) at (9,0) {$p_4$};
    \node (b) at (10,0) {$p_1$};
    \node (c) at (11,0) {$p_3$};
    \node (d) at (12.5,0) {$p_2$};

    \draw[->] (a) to node {} (b);
    \draw[->] (b) to node {} (c);
    \draw[->] (c) to node {} (d);
    \draw[->] (d) to [in=-10, out=190, looseness=1] node {} (c);
\end{tikzpicture}
\end{center}

\item Series $\mathfrak{S}||E$. Maps in this series can be given by the composition $\mathfrak{S}|Aj\circ \mathfrak{S}|Bk$ for some $j, k\in\{1,2,\cdots,24\}$. For example, 
\begin{center}
\begin{tikzpicture}
    \node (a) at (0,0) {$p_4$};
    \node (b) at (1,0) {$p_2,$};
    \node (c) at (2,0) {$p_1$};
    \node (d) at (3,0) {$p_3$};

    \draw[->] (a) to [in=170, out=10, looseness=1] node {} (b);
    \draw[->] (b) to [in=-10, out=-170, looseness=1] node {} (a);
    \draw[->] (c) to node {} (d);
    \draw (d) edge[loop right] node {} (d);
\end{tikzpicture}
=
\begin{tikzpicture}
    \node (a) at (5,0) {$p_1$};
    \node (b) at (6,0) {$p_2$};
    \node (c) at (7,0) {$p_3$};
    \node (d) at (8,0) {$p_4$};

    \draw[->] (a) to node {} (b);
    \draw[->] (b) to node {} (c);
    \draw[->] (c) to node {} (d);
    \draw[->] (d) to [in=-20, out=200, looseness=1] node {} (b);
\end{tikzpicture}
$\circ$
\begin{tikzpicture}
    \node (a) at (9,0) {$p_4$};
    \node (b) at (10,0) {$p_1$};
    \node (c) at (11,0) {$p_2$};
    \node (d) at (12.5,0) {$p_3.$};

    \draw[->] (a) to node {} (b);
    \draw[->] (b) to node {} (c);
    \draw[->] (c) to node {} (d);
    \draw[->] (d) to [in=-10, out=190, looseness=1] node {} (c);
\end{tikzpicture}
\end{center}

\end{itemize}

\end{proof}

\begin{proof}[Proof of Theorem \ref{thm13}]
The result follows from a combination of Lemma \ref{lem21}, \ref{lem19} and \ref{lem20} instantly. 
\end{proof}

There are other generating sets of strictly pre-periodic maps whose compositions can induce all the strictly pre-periodic maps on sets of four points, however, one needs at least $12$ of them.

\subsection{The simplification of \cite[Question 1.2]{DKM}}

Upon the algebraic structure of the periodic family and the strictly pre-periodic family of maps demonstrated in the previous subsections, we can easily achieve the following result. 

\begin{theorem}\label{thm15}
All the  maps  can be given by finite compositions of maps in 
\begin{center}
$\{\mathfrak{P}|1, \mathfrak{P}||A1, \mathfrak{S}|A1\}$
\end{center}
on $X=\{p_1, p_2, p_3, p_4\}$.
\end{theorem}
\begin{proof}
It suffices for us to show that all the maps in $\{\mathfrak{S}|Aj\}_{1\leq j\leq 12}$ can be given by the finite compositions of maps in the periodic family plus $\{\mathfrak{S}|A1\}$, in virtue of Theorem \ref{thm12} and Theorem \ref{thm14}. We only express two of them by compositions of the three given maps, with expressions of rest of them left to the readers as exercises.

\begin{itemize}
\item Map \begin{tikzpicture}
    \node (a) at (5,0) {$p_1$};
    \node (b) at (6,0) {$p_2$};
    \node (c) at (7,0) {$p_4$};
    \node (d) at (8,0) {$p_3$};

    \draw[->] (a) to node {} (b);
    \draw[->] (b) to node {} (c);
    \draw[->] (c) to node {} (d);
    \draw[->] (d) to [in=-20, out=200, looseness=1] node {} (b);
\end{tikzpicture}. 

This map can be given by the composition $\mathfrak{S}|A1\circ (\mathfrak{P}|1\circ\mathfrak{P}||A1)^3$.

\item Map \begin{tikzpicture}
    \node (a) at (5,0) {$p_2$};
    \node (b) at (6,0) {$p_3$};
    \node (c) at (7,0) {$p_4$};
    \node (d) at (8,0) {$p_1$};

    \draw[->] (a) to node {} (b);
    \draw[->] (b) to node {} (c);
    \draw[->] (c) to node {} (d);
    \draw[->] (d) to [in=-20, out=200, looseness=1] node {} (b);
\end{tikzpicture}.

This map can be given by the composition $(\mathfrak{P}||A1\circ\mathfrak{P}|1)^3\circ(\mathfrak{S}|A1\circ (\mathfrak{P}|1)^3)^2$.

\end{itemize}

\end{proof}

In virtue of Theorem \ref{thm15} a simplification of \cite[Question 1.2]{DKM} is instant.

\begin{proof}[Proof of Theorem \ref{thm16}]
The result follows from a combination of Theorem \ref{thm15} and Corollary \ref{cor5} instantly.
\end{proof}

Obviously there are other simplifications of \cite[Question 1.2]{DKM}, as long as one finds a generating set of the monoid constituted by all maps on $X$, however, one needs at least $3$ of them. Moreover, these algebraic structures on collections of maps seem to be of individual interest themselves.

\end{document}